\documentclass[12pt]{article}
\usepackage[a4paper]{geometry}
\usepackage{amsthm}
\usepackage{amsmath}
\usepackage{amssymb}
\usepackage{amsfonts}
\usepackage{graphicx}
\usepackage{graphics}
\usepackage[bf,SL,BF]{subfigure}

 \numberwithin{equation}{section}
\def\eref#1{(\ref{#1})}
\def\N{\mathbb{N}}

\def\R{\mathbb{R}}

\def\Rh{\mathcal{R}_h}
\def\Rht{\mathcal{R}_{h,t}}

\def\A{\mathcal{A}}

\def\Q{\mathcal{Q}}
\def\<{\big\langle}
\def\>{\big\rangle}

\def\diiv{\operatorname{div}}

\def\Tr{\operatorname{Trace}}

\def\det{\operatorname{det}}

\def\esssup{{\operatorname{esssup}}}

\newtheorem{Lemma}{Lemma}[section]

\newtheorem{Theorem}{Theorem}[section]
\newtheorem{Proposition}{Proposition}[section]

\newtheorem{Condition}{Condition}[section]

\theoremstyle{remark}
\newtheorem{Remark}{Remark}[section]

\theoremstyle{definition}

\newtheorem{Definition}{Definition}[section]

 \theoremstyle{definition}
  \newtheorem{example}{Example}

\newcommand{\figbox}[1]{%
  \fbox{%
    \vbox to 1in{%
    \vfil
    \hbox to 2in{%
      \hfil
      #1%
      \hfil}%
    \vfil}}}

\newcommand{\goodgap}{%
  \hspace{\subfigcapskip}}

\begin{document}

\title{Homogenization of parabolic equations with a continuum of space and time scales.\protect\footnotetext{AMS 1991 {\it{Subject Classification}}. Primary 34E13,35B27 ; secondary
35B26, 68P30,  60F05, 35B05.} \protect\footnotetext{{\it{Key words
and phrases}}. Multi scale
problem, compensation, homogenization,  up-scaling, compression.}}          

\author{Houman Owhadi\footnote{
 California Institute of Technology
Applied \& Computational Mathematics, Control \& Dynamical systems,
MC 217-50 Pasadena , CA 91125, owhadi@caltech.edu}  and Lei
Zhang\footnote{
 California Institute of Technology
Applied \& Computational Mathematics MC 217-50 Pasadena , CA 91125,
zhanglei@acm.caltech.edu}}
\date{December 21, 2005}
\maketitle \abstract{This paper addresses the issue of
homogenization of linear divergence form parabolic operators in
situations where no ergodicity and no scale separation in time or
space are available. Namely, we consider divergence form linear
parabolic operators in $\Omega \subset \R^n$ with $L^\infty(\Omega
\times (0,T))$-coefficients. It appears that the inverse operator
maps the unit ball  of $L^2(\Omega\times (0,T))$ into a space of
functions which at small (time and space) scales are close in
$H^1$-norm to a functional space of dimension $n$. It follows that
once one has solved these equations at least $n$-times it is
possible to homogenize them both in space and in time, reducing the
number of operations counts necessary to obtain further solutions.
In practice we show that under a Cordes type condition that the
first order time derivatives and second order space derivatives of
the solution of these operators with respect to harmonic coordinates
are in $L^2$ (instead of $H^{-1}$ with Euclidean coordinates). If
the medium is time independent then it is sufficient to solve $n$
times the associated elliptic equation in order to homogenize the
parabolic equation.
 }

\section{Introduction and main results}
Let $\Omega$ be a bounded and convex domain of class $C^2$ of
$\R^n$. Let $T>0$. Consider the following parabolic PDE
\begin{equation}\label{ghjh52}
\begin{cases}
\partial_t u=\diiv\big(a(x,t)\nabla u(x,t)\big)+g\quad \text{in}\quad \Omega\times (0,T) \\
u(x,t)=0 \quad \text{for}\quad (x,t) \in \big(\partial \Omega \times
(0,T) \big) \cup \big(\Omega \cup \{t=0\}\big).
\end{cases}
\end{equation}
Write $\Omega_T:=\Omega \times (0,T)$. $g$ is a function in
$L^{2}(\Omega_T)$.  $(x,t) \rightarrow a(x,t)$ is a mapping from
$\Omega_T$ into the space of  symmetric positive definite matrices
with entries in $L^\infty(\Omega_T)$. Assume $a$ to be uniformly
elliptic on the closure of $\Omega_T$. This paper  addresses the
issue of the homogenization of \eref{ghjh52} in space and time in
situations where scale separation and ergodicity  at small scales
are not available (see   \cite{BeLiPa78}, \cite{JiKoOl91} and
\cite{Al01} for an introduction to classical homogenization theory).
For that purpose, we will introduce in subsection
\ref{kjskjsjshj871} theorems establishing under Cordes type
conditions the increase of regularity of solutions of \eref{ghjh52}
when derivatives are taken with respect to harmonic coordinates
instead of Euclidean coordinates. In subsections \ref{sub2},
\ref{jksjhs89b} these results will be used to homogenize
\eref{ghjh52} in space and in time. More precisely, assume $a$ to be
written on a fine tessellation with $N$ degrees of freedom. If $a$
is time independent, then by solving $n$-times an elliptic boundary
value-problem associated to \eref{ghjh52} (at a cost of $O(N (\ln
N)^{n+3})$ operations using the Hierarchical matrix method
\cite{Beben05}) it is possible approximate the solutions of
\eref{ghjh52} by solving an homogenized operator with $N^\alpha$
degrees of freedom ($\alpha <1$, $\alpha=0.2$ for instance) or with
a fixed number $M$ of degrees of freedom (numerical experiments
given at the end of this paper have been conducted with $N=16641$
and $M=9$), this problem is of practical importance for oil
extraction and reservoir modeling in geophysics. If $a$ is also
characterized by a continuum of time scales, then the method
presented here does not reduce the number of operation counts
necessary to solve \eref{ghjh52} only one time. However if one needs
to solve \eref{ghjh52} $K$ ($K>n$) times (with different right hand
sides) then by solving \eref{ghjh52} $n$ times it is possible to
obtain an approximation of the solutions of \eref{ghjh52} by solving
an homogenized (in space and time) parabolic equation written on a
coarse tessellation with coarse time steps.

\subsection{Compensation phenomenon}\label{kjskjsjshj871}
 Let $F$ be the solution of
the following parabolic equation
\begin{equation}\label{ghjagsash52}
\begin{cases}
\partial_t F=\diiv\big(a(x,t)\nabla F(x,t)\big)\quad \text{in}\quad \Omega_T \\
F(x,t)=x \quad \text{for}\quad (x,t) \in \big(\partial \Omega \times
(0,T) \big)\\
\diiv\big(a(x,0)\nabla F(x,0)\big)=0\quad \text{in}\quad \Omega.
\end{cases}
\end{equation}
By \eref{ghjagsash52} we mean that $F:=(F_1,\ldots, F_n)$ is a
$n$-dimensional vector field such that each of its entries satisfies
\begin{equation}\label{ghjasgsash52}
\begin{cases}
\partial_t F_i=\diiv\big(a(x,t)\nabla F_i(x,t)\big)\quad \text{in}\quad \Omega_T \\
F_i(x,t)=x_i \quad \text{for}\quad (x,t) \in \big(\partial \Omega
\times (0,T) \big) \\
\diiv\big(a(x,0)\nabla F_i(x,0)\big)=0\quad \text{in}\quad \Omega.
\end{cases}
\end{equation}
Observe that if $a$ is time independent then $F$ is the solution of
an elliptic boundary value problem.

\begin{Definition}
Write
\begin{equation}\label{jhscczxhjd}
\sigma:={^t\nabla F}a\nabla F.
\end{equation}
\end{Definition}

Write $\beta_\sigma$ the Cordes parameter associated to $\sigma$
defined by
\begin{equation}\label{sshgdd7641}
\beta_{\sigma}:=\esssup_{(x,t)\in \Omega_T}\Big(
n-\frac{\big(\Tr[\sigma]\big)^2}{\Tr[^t\sigma \sigma] }\Big).
\end{equation}
Observe that since
\begin{equation}\label{sshgdmmmd7641}
\beta_{\sigma}=\esssup_{(x,t)\in \Omega_T}\Big(
n-\frac{\big(\sum_{i=1}^n
\lambda_{i,\sigma(x,t)}\big)^2}{\sum_{i=1}^n
\lambda_{i,\sigma(x,t)}^2}\Big).
\end{equation}
where $(\lambda_{i,M})$ denotes the eigenvalues of $M$,
$\beta_{\sigma}$ is a measure of the anisotropy of $\sigma$.

\subsubsection{Time independent medium.}
In this subsection we assume that $a$ does not depend on time $t$.
 Write for $p\geq 2$, $W^{2,p}_{D}$ ($D$ for Dirichlet boundary
condition) the Banach space $W^{2,p}_D(\Omega)\cap
W^{1,p}_0(\Omega)$. Equip $W^{2,p}_D(\Omega)$ with the norm
\begin{equation}
\|v\|_{W^{2,p}_D(\Omega)}^2:=\int_{\Omega}\big(\sum_{i,j}(\partial_i
\partial_j v)^2\big)^\frac{p}{2}.
\end{equation}
Equip the space  $L^p(0,T,W^{2,p}_D(\Omega))$ with the norm
\begin{equation}
\|v\|_{L^p(0,T,W^{2,p}_D(\Omega))}^p=\int_0^T\int_{\Omega}
\big(\sum_{i,j} (\partial_i
\partial_j v)^2\big)^\frac{p}{2} \,dx\,dt.
\end{equation}

\begin{Theorem}\label{ksjshjsxxsd8721}
Assume that $\partial_t a\equiv 0$, $g\in L^2(\Omega_T)$, $\Omega$
is convex, $\beta_\sigma<1$ and $(\Tr[\sigma])^{\frac{n}{4}-1}\in
L^\infty(\Omega)$ then $u\circ F^{-1}\in L^2(0,T,W^{2,2}_D(\Omega))$
and
\begin{equation}
\|u\circ F^{-1}\|_{L^2(0,T,W^{2,2}_D(\Omega))}\leq \frac{C
}{1-\beta_\sigma^\frac{1}{2}} \|g\|_{L^2(\Omega_T)}.
\end{equation}
\end{Theorem}
\begin{Remark}
The constant $C$ can be written
$$C=\frac{C_n}{(\lambda_{\min}(a))^\frac{n}{4}}
\big\|(\Tr[\sigma])^{\frac{n}{4}-1}\big\|_{L^\infty(\Omega)}.$$
Through this paper, we write
\begin{equation}
\lambda_{\min}(a):=\inf_{(x,t)\in\Omega_T}\inf_{l\in \R^n,
|l|=1}{^tl. a(x,t).l}.
\end{equation}
\end{Remark}
\begin{Remark}
According to theorem \ref{ksjshjsxxsd8721} although the second order
derivatives of $u$ with respect to Euclidean coordinates are only in
$L^2(0,T,H^{-1}(\Omega))$, they are in $L^2(\Omega_T)$ with respect
to harmonic coordinates.
\end{Remark}
\begin{Remark}
Observe that if $a$ is time independent then $F$ and $\sigma$ are
time independent and $F$ is the solution of the following elliptic
problem:
\begin{equation}\label{dgdgfsghsxzf62}
\begin{cases}
\diiv a \nabla F=0 \quad \text{in}\quad \Omega\\
F(x)=x \quad \text{on}\quad \partial \Omega.
\end{cases}
\end{equation}

 In dimension one $F$ is trivially an
homeomorphism. In dimension $2$ this property follows from
topological constraints \cite{MR2001070} (even with $a_{i,j}\in
L^\infty(\Omega)$), \cite{MR1892102} (one can also deduce from
\cite{MR1892102} that for $n=2$, if $a$ is smooth then the
conditions $\beta_\sigma<1$ and $(\Tr[\sigma])^{-1}\in
L^\infty(\Omega)$ are satisfied). In dimension three and higher $F$
can be non-bijective even if $a$ is smooth, we refer to
\cite{MR1892102} and \cite{MR2073507}, however in dimension $3$ the
assumption $(\Tr[\sigma])^{\frac{n}{4}-1}\in L^\infty(\Omega_T)$
implies that $F$ is an homeomorphism. If $n\geq 4$ we need to assume
that $F$ is an homeomorphism to prove the theorem.
\end{Remark}

\begin{Remark}
In fact the condition $(\Tr(\sigma))^{-1}\in L^p(\Omega_T)$ for
$p<\infty$ depending on $n$ is sufficient to obtain theorem
\ref{ksjshjsxxsd8721} and the following compensation theorems. For
the sake of clarity this paper has been restricted to
$(\Tr(\sigma))^{-1} \in L^\infty(\Omega_T)$.
\end{Remark}

\begin{Remark}
Write
\begin{equation}\label{sjjdhd27}
\mu_\sigma:=\esssup_{\Omega_T}\frac{\lambda_{\max}(\sigma)}{\lambda_{\min}(\sigma)}.
\end{equation}
It is easy to check that $\mu_{\sigma}$ is bounded by an increasing
function of $(1-\beta_{\sigma})^{-1}$ and in dimension two
$\beta_{\sigma}<1$ is equivalent to $\mu_{\sigma}<\infty$.
\end{Remark}
\begin{Remark}
Theorem \ref{ksjshjsxxsd8721} has been called compensation
phenomenon because the composition by $F^{-1}$ increases the
regularity of $u \in L^2(0,T, H^1_0(\Omega))$. The choice of this
name has been motivated by
 F. Murat and L. Tartar's work on
H-convergence \cite{MR1493039}  which is also based on on a
regularization property called compensated compactness or div-curl
lemma introduced in the 70's by Murat and Tartar \cite{MR506997},
\cite{MR584398} (we also refer to \cite{MR1225511} for refinements
of the div-curl lemma).
\end{Remark}

\begin{figure}[httb]
\begin{center}
\includegraphics[%
  scale=0.3]{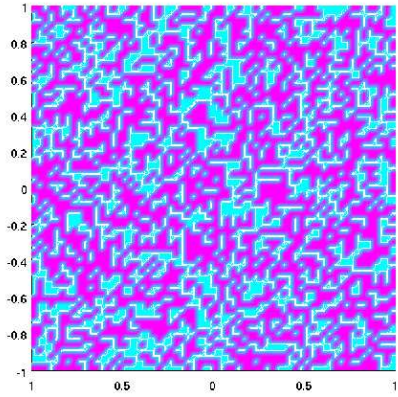}
\caption{Site Percolation} \label{ap7}
\end{center}
\end{figure}
\begin{figure}[httb]
  \begin{center}
    \subfigure[$u$.]
    {\includegraphics[width=0.45\textwidth,height= 0.4\textwidth]{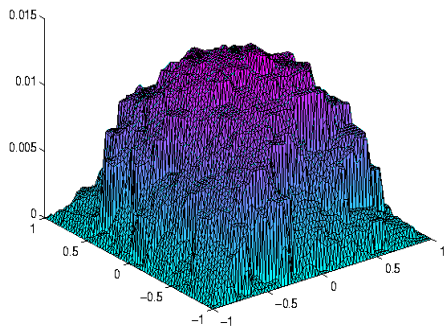}}
    \goodgap
    \subfigure[$u\circ F^{-1}$.]
    {\includegraphics[width=0.45\textwidth,height= 0.4\textwidth]{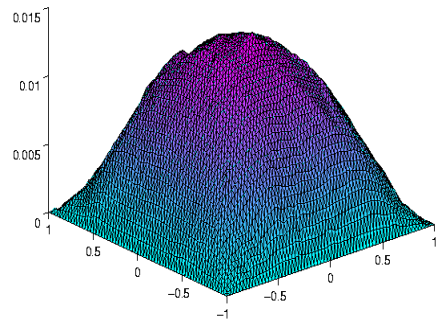}}\\
    \subfigure[$\partial_x u$.]
    {\includegraphics[width=0.45\textwidth,height= 0.4\textwidth]{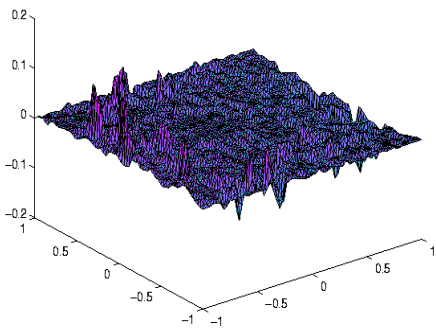}}
    \goodgap
    \subfigure[$\partial x (u\circ F^{-1})$.]
    {\includegraphics[width=0.45\textwidth,height= 0.4\textwidth]{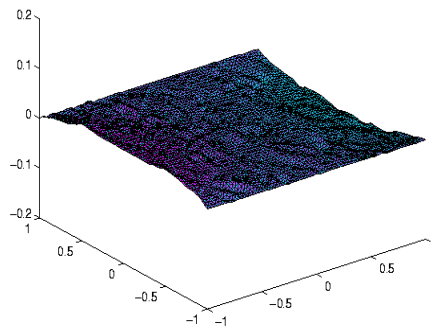}}\\
    \caption{$u$, $u\circ F^{-1}$, $\partial_x u$ and $\partial x (u\circ F^{-1})$ at time $t=1$ for the time independent site percolating medium.}
    \label{unut100tip7}
\end{center}
\end{figure}

The compensation phenomena presented in this subsection can be
observed numerically. In figure \ref{ap7}, the value of $a$ is set
to be equal to $1$ or $100$ with probability $1/2$ on each triangle
of a fine mesh characterized by $16641$ nodes and $32768$ triangles.
 \eref{ghjh52} has been solved numerically on that mesh with $g=1$. $u$, $u\circ F^{-1}$, $\partial_x u$ and $\partial x (u\circ F^{-1})$ have been plotted
 at time $t=1$ in figure \ref{unut100tip7}.

\clearpage

In situations where $g\in L^\infty(0,T,L^2(\Omega))$,  $\partial_t
g\in L^2(0,T, H^{-1}(\Omega))$ or $g\in L^p(\Omega_T)$ with $p>2$,
one can obtain a higher regularity for $u\circ F^{-1}$. This is the
object of the following theorems.

\begin{Theorem}\label{ksjshszjscxsxsdded8721}
Assume that $\Omega$ is convex, $g\in L^\infty(0,T,L^2(\Omega))$, \\
$\partial_t g\in L^2(0,T, H^{-1}(\Omega))$, $\partial_t a\equiv 0$,
$\beta_\sigma<1$ and $(\Tr[\sigma])^{\frac{n}{4}-1}\in
L^\infty(\Omega_T)$ then for all $t\in [0,T]$, $u\circ F^{-1}(.,t)
\in W^{2,2}_D(\Omega)$ and
\begin{equation}\label{equzalundi}
 \|u\circ F^{-1}(.,t)\|_{W^{2,2}_D(\Omega)}\leq \frac{C
}{1-\beta_\sigma^\frac{1}{2}}
\Big(\big\|g\big\|_{L^\infty(0,T,L^2(\Omega))}+\|\partial_t
g\|_{L^2(0,T,H^{-1}(\Omega))}\Big).
\end{equation}
\end{Theorem}
\begin{Remark}
The constant $C$ can be written
$$C=\frac{C_{n,\Omega}}{(\lambda_{\min}(a))^\frac{n}{4}}
\big\|(\Tr[\sigma])^{\frac{n}{4}-1}\big\|_{L^\infty(\Omega_T)}(1+\frac{1}{\lambda_{\min}(a)})^\frac{1}{2}.$$
\end{Remark}

\begin{Theorem}\label{ksjshjsxgghxsd8721}
Assume that $\Omega$ is convex, $g(.,0)\in L^2(\Omega)$,
\\$\partial_t g\in L^2(0,T, H^{-1}(\Omega))$,  $g\in L^p(\Omega_T)$,
$\partial_t a\equiv 0$, $\beta_\sigma<1$ and
$(\Tr[\sigma])^{\frac{n}{4}-1}\in L^\infty(\Omega_T)$ then there
 exists a real number $p_0>2$ depending only on $n,\Omega$
and $\beta_\sigma$ such that for each $p$ such that $2\leq p<p_0$
one has
\begin{equation}
\begin{split}
\|u\circ F^{-1}\|_{L^p(0,T,W^{2,p}_D(\Omega))} \leq &\frac{C
}{1-\beta_\sigma^\frac{1}{2}}
\big(\|g\|_{L^p(\Omega_T)}\\&+\big\|g(.,0)\big\|_{L^2(\Omega)}+\|\partial_t
g\|_{L^2(0,T,H^{-1}(\Omega))}\big).
\end{split}
\end{equation}
\end{Theorem}
\begin{Remark}
The constant $C$ can written
$$C=\frac{C_{n,\Omega,p}}{(\lambda_{\min}(a))^\frac{n}{4}}
\big\|(\Tr[\sigma])^{\frac{n}{4}-1}\big\|_{L^\infty(\Omega_T)}\big(1+\frac{1}{\lambda_{\min}(a)}\big)^\frac{1}{2}.$$
\end{Remark}
Write
\begin{equation}\label{hdaassz7}
\begin{split}
\|v\|_{C^{\gamma}(\Omega)}:=\sup_{x,y\in \Omega,
x\not=y}\frac{|v(x)-v(y)|}{|x-y|^\gamma}.
\end{split}
\end{equation}

\begin{Theorem}\label{hdgwsssawwsd7}
Assume that $n\leq 2$, $\Omega$ is convex, $g(.,0)\in L^2(\Omega)$,
$\partial_t g\in L^2(0,T, H^{-1}(\Omega))$,  $g\in L^p(\Omega_T)$,
$\partial_t a\equiv 0$, $\beta_\sigma<1$, $(\Tr[\sigma])^{-1}\in
L^\infty(\Omega_T)$ and $g\in L^2[0,T;L^{p^*}(\Omega)]$ with
$2<p^*$. Then there exists $p\in (2,p^*]$ and $\gamma(p)>0$ such
that
\begin{equation}\label{hdhsxdzssaqssgc7}
\begin{split}
\big(\int_0^T \big\|\nabla (u\circ
F^{-1})(.,t)\big\|_{C^{\gamma}(\Omega)}^2\,dt\big)^\frac{1}{2} \leq
&\frac{C }{1-\beta_\sigma^\frac{1}{2}}
\big(\|g\|_{L^p(\Omega_T)}\\&+\big\|g(.,0)\big\|_{L^2(\Omega)}+\|\partial_t
g\|_{L^2(0,T,H^{-1}(\Omega))}\big).
\end{split}
\end{equation}
\end{Theorem}
\begin{Remark}
The constant $C$ in \eref{hdhsxdzssaqssgc7} depends on $n$, $p$,
$\Omega$, $\lambda_{\min}(a)$ and
$\big\|(\Tr(\sigma))^{-1}\big\|_{L^\infty(\Omega_T)}$. It is easy to
check that if $n=1$ then the theorem is valid with $\gamma=1/2$.
\end{Remark}

In the following theorems  $\Omega$ is not assumed to be convex.

\begin{Theorem}\label{hdgwsswdswsd7}
Assume $n\geq 2$ and $\partial_t a\equiv 0$. Let $p>2$. There exist
a constant $C^*=C^*(n,\partial \Omega)>0$ a real number $\gamma>0$
depending only on $n,\Omega$ and $p$ such that if
$\beta_{\sigma}<C^*$ then
\begin{equation}\label{hdhsxdznnnnsgxc7}
\begin{split}
\big(\int_0^T \big\|\nabla (u\circ
F^{-1})(.,t)\big\|_{C^{\gamma}(\Omega)}^2\,dt\big)^\frac{1}{2} \leq
 C
\big(&\|g\|_{L^p(\Omega_T)}+\big\|g(.,0)\big\|_{L^2(\Omega)}\\&+\|\partial_t
g\|_{L^2(0,T,H^{-1}(\Omega))}\big).
\end{split}
\end{equation}
\end{Theorem}
\begin{Remark}
The constant $C$ in \eref{hdhsxdznnnnsgxc7} depends on $n$,
$\gamma$, $\Omega$, $C^*$, $\lambda_{\min}(a)$ and
$\big\|(\Tr(\sigma))^{\frac{n}{2p}-1}\big\|_{L^\infty(\Omega_T)}$.
\end{Remark}

It is easy to check that if $a=e(x)S(x,t)$ where $e$ is a time
independent symmetric uniformly elliptic matrix with
$L^\infty(\Omega)$ entries and $S$ is a regular uniformly positive
function then the results given in this sub-section and the
homogenization schemes of sub-section \ref{sub2} remain valid with
the time independent harmonic coordinates associated to $e$, i.e.
solution of $-\diiv e \nabla F=0$.

\subsubsection{Medium with a continuum of time scales.}
In this subsection the entries of $a$ are merely in
$L^\infty(\Omega_T)$.
 We need to introduce the following Cordes type condition.
\begin{Condition}\label{slkssjk88271}
We say that condition \ref{slkssjk88271} is satisfied if and only if
there exists $\delta\in (0,\infty)$ and $\epsilon>0$ such that
\begin{equation}
\esssup_{\Omega_T}
\frac{\delta^2\Tr[{^t\sigma\sigma}]+1}{\big(\delta\Tr[\sigma]+1\big)^2}\leq
\frac{1}{n+\epsilon}.
\end{equation}
\end{Condition}

Write
\begin{equation}
z_\sigma:=\esssup_{\Omega_T} n \frac{\Tr[{^t\sigma
\sigma}]}{(\Tr[\sigma])^2}.
\end{equation}
Observe that $z_\sigma$ is a measure of anisotropy of $\sigma$, in
particular $1\leq z_\sigma \leq n$ and $z_\sigma=1$ if $\sigma$ is
isotropic.
 Write
\begin{equation}
y_\sigma:=\|\Tr[\sigma]\|_{L^\infty(\Omega_T)}\big\|(\Tr[\sigma])^{-1}\big\|_{L^\infty(\Omega_T)}.
\end{equation}

\begin{Proposition}\label{kwjdjwj7}
If $\|\Tr[\sigma]\|_{L^\infty(\Omega_T)}<\infty$ and
$\big\|(\Tr[\sigma])^{-1}\big\|_{L^\infty(\Omega_T)}<\infty$ then
condition \ref{slkssjk88271} is satisfied with
\begin{equation}
\delta:=n \big\|(\Tr[\sigma])^{-1}\big\|_{L^\infty(\Omega_T)}
\end{equation}
and with $\epsilon:=\frac{2ny_\sigma-n}{2n y^2_\sigma}$ provided
that $z_\sigma \leq 1+\frac{\epsilon}{n}$.
\end{Proposition}
\begin{Remark}
Observe that in dimension one $z_\sigma=1$, thus for $n=1$ condition
\ref{slkssjk88271} is satisfied is $\Tr[\sigma]\in
L^\infty(\Omega_T)$ and $(\Tr[\sigma])^{-1}\in L^\infty(\Omega_T)$.
\end{Remark}
We have the following theorems
\begin{Theorem}\label{skjhsj823}
Assume that $\Omega$ is convex, and condition \ref{slkssjk88271} is
satisfied then $u\circ F^{-1}\in
L^2\big(0,T,W^{2,2}_D(\Omega)\big)$, $\partial_t (u\circ F^{-1})\in
L^2(\Omega_T)$ and
\begin{equation}
\|u\circ F^{-1}\|_{L^2(0,T,W^{2,2}_D(\Omega))}+\|\partial_t(u\circ
F^{-1})\|_{L^2(\Omega_T)}\leq C \|g\|_{L^2(\Omega_T)}
\end{equation}
where $C$ depends on $\Omega$, $n$, $\delta$ and $\epsilon$.
\end{Theorem}
\begin{Remark}
According to theorem \ref{skjhsj823} although the second order space
derivatives and first order time derivatives of $u$ with respect to
Euclidean coordinates are only in $L^2(0,T,H^{-1}(\Omega))$, they
are in $L^2(\Omega_T)$ with respect to harmonic coordinates.
\end{Remark}

Similarly we obtain the following theorems in situations where $g\in
L^p(\Omega_T)$ with $p>2$.
\begin{Theorem}\label{skjhsjd823}
Assume that $\Omega$ is convex, and condition \ref{slkssjk88271} is
satisfied then  there exists a number $p_0>2$ depending on
$n,\Omega,\epsilon$ such that for $p\in (2,p_0)$, $u\circ F^{-1}\in
L^p\big(0,T,W^{2,p}_D(\Omega)\big)$, $\partial_t (u\circ F^{-1})\in
L^p(\Omega_T)$ and
\begin{equation}
\|u\circ F^{-1}\|_{L^p(0,T,W^{2,p}_D(\Omega))}+\|\partial_t(u\circ
F^{-1})\|_{L^p(\Omega_T)}\leq C \|g\|_{L^p(\Omega_T)}
\end{equation}
where $C$ depends on $\Omega$, $n$, $\delta$ and $\epsilon$.
\end{Theorem}

\begin{Theorem}\label{skjhsssej823}
Assume that $\Omega$ is convex, and condition \ref{slkssjk88271} is
satisfied then there exists a number $\alpha_0>2$ depending on
$n,\Omega,\epsilon$ such that for $\alpha \in (0,\alpha_0)$, $\nabla
(u\circ F^{-1})\in L^2(0,T,C^\alpha(\Omega))$ and
\begin{equation}
\|\nabla (u\circ F^{-1})(.,t)\|_{L^2(0,T,C^\alpha(\Omega))}\leq C
\|g\|_{L^p(\Omega_T)}
\end{equation}
where $C$ depends on $\Omega$, $\delta$, $n$, and $\epsilon$.
\end{Theorem}

\begin{figure}[httb]
\begin{center}
\includegraphics[%
  scale=0.3]{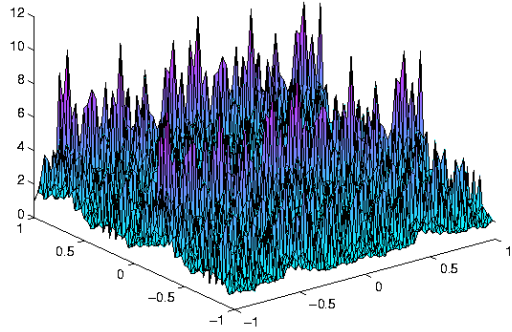}
\caption{$a$ at time $0$.} \label{mediap4}
\end{center}
\end{figure}

These compensation phenomena can be observed numerically. We
consider in dimension $n=2$,
\begin{equation}\label{ksjhskjhdhkdjh872y2}
\begin{split}
a(x,y,t)=\frac{1}{6}(\sum_{i=1}^5 \frac{1.1+\sin(2\pi
x'/\epsilon_{i})}{1.1+\sin(2\pi
y'/\epsilon_{i})}+\sin(4x'^{2}y'^{2})+1)
\end{split}
\end{equation}
with $x'=x+\sqrt{2}t$, $y'=y-\sqrt{2}t$, $\epsilon_{1}=\frac{1}{5}$,
$\epsilon_{2}=\frac{1}{13}$, $\epsilon_{3}=\frac{1}{17}$,
$\epsilon_{4}=\frac{1}{31}$ and $\epsilon_{5}=\frac{1}{65}$. This
medium has been plotted in figure \ref{mediap4} at time $0$ (observe
that $\lambda_{\max}(a)/\lambda_{\min}(a)\sim 100$).

\begin{figure}[httb]
  \begin{center}
    \subfigure[$\partial_x u$.]
    {\includegraphics[width=0.35\textwidth,height= 0.3\textwidth]{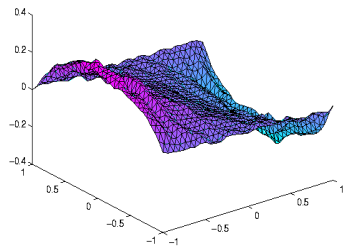}}
    \goodgap
    \subfigure[$\partial_x (u\circ F^{-1})$.]
    {\includegraphics[width=0.35\textwidth,height= 0.3\textwidth]{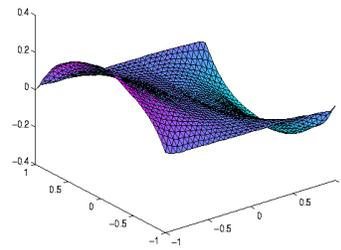}}\\
    \caption{$\partial_x u$ and $\partial_x (u\circ F^{-1})$ at time $t=0.3$ for the Multi-scale Trigonometric time dependent Medium.}
    \label{unut03p4}
\end{center}
\end{figure}

\eref{ghjh52} has been solved numerically on that mesh with $g\equiv
1$ on the fine mesh characterized by $16641$ nodes and $32768$
triangles. Figure \ref{unut03p4} shows $\partial_x u$ and
$\partial_x (u\circ F^{-1})$ at time $0.3$.

\begin{figure}[httb]
\begin{center}
    \subfigure[$u$ and $u\circ F^{-1}$.]
    {\includegraphics[width=0.5\textwidth,height= 0.4\textwidth]{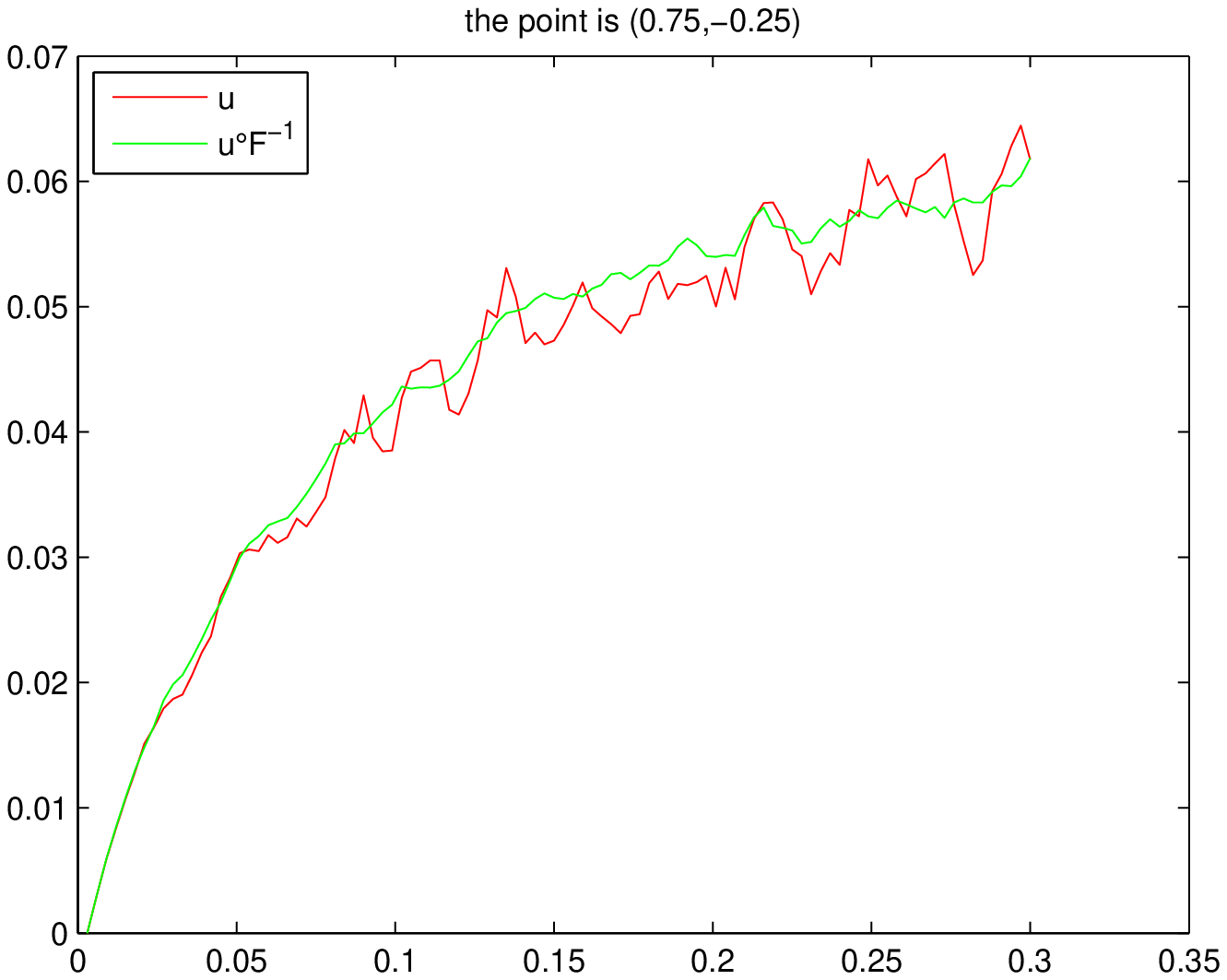}}
\caption{$t\rightarrow u(x_0,t), u\circ F^{-1}(x_0,t)$  from $t=0$
to $t=0.3$ with $x_0:=(0.75,-0.25)$} \label{fixpointu00}
\end{center}
\end{figure}

\begin{figure}[httb]
\begin{center}
    \subfigure[$\partial_x u$ and $\partial_x (u\circ F^{-1})$.]
    {\includegraphics[width=0.45\textwidth,height= 0.4\textwidth]{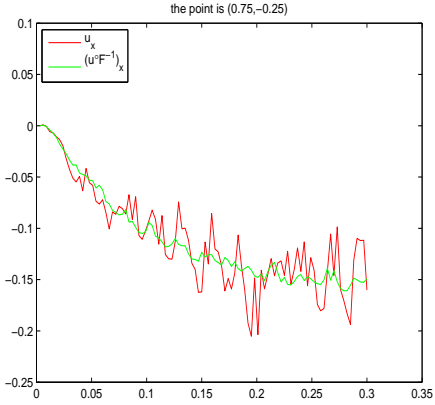}}
     \goodgap
    \subfigure[$\partial_y u$ and $\partial_y (u\circ F^{-1})$.]
    {\includegraphics[width=0.45\textwidth,height= 0.4\textwidth]{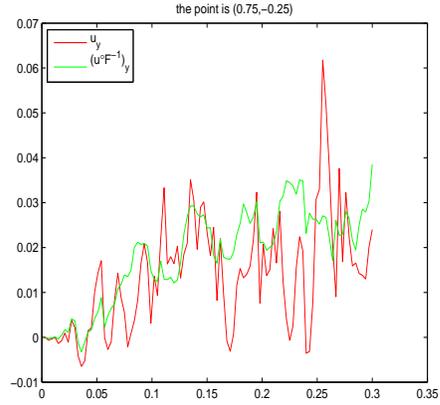}}\\
\caption{$t\rightarrow \nabla u(x_0,t),\nabla u\circ F^{-1}(x_0,t)$
from $t=0$ to $t=0.3$ with $x_0:=(0.75,-0.25)$} \label{fixpointu}
\end{center}
\end{figure}

In figure \ref{fixpointu00} and \ref{fixpointu}, the value of $x_0$
is set to $x_0:=(0.75,-0.25)$ and the curves $t\rightarrow u(x_0,t),
u\circ F^{-1}(x_0,t), \nabla u(x_0,t),\nabla u\circ F^{-1}(x_0,t)$
are plotted from $t=0$ to $t=0.3$

\clearpage

\subsection{Homogenization in space.}\label{sub2}
Let $X_h$ be a finite dimensional subspace of $H^1_0(\Omega)\cap
W^{1,\infty}(\Omega)$\footnote{$W^{1,\infty}$ is the usual space of
uniformly Lipschitz continuous functions.} with the following
approximation property: there exists a constant $C_X$ such that for
all $f\in W^{2,2}_D(\Omega)$

\begin{equation}\label{approp}
\inf_{v\in X_h} \|f-v\|_{H^1_0(\Omega)}\leq C_X h
\|f\|_{W^{2,2}_D(\Omega)}.
\end{equation}

It is known and easy to check that the set of piecewise linear
functions on a triangulation of $\Omega$ satisfies condition
\eref{approp}  provided that the length of the edges of the
triangles are bounded by $h$ ($C_X$ in \eref{approp} being given by
the aspect ratio of the triangles).

For media characterized by a continuum of time scales we will
consider twice differentiable elements
 satisfying the following usual
inverse inequalities (see section 1.7 of \cite{ErGu04}): for $v\in
X_h$,
\begin{equation}\label{appsrop3}
 \|v\|_{W^{2,2}_D(\Omega)}\leq
C_X h^{-1} \|v\|_{H^1_0(\Omega)}.
\end{equation}
and
\begin{equation}\label{appssswerop3}
 \| v\|_{H^1_0(\Omega)}\leq
C_X h^{-1} \|v\|_{L^2(\Omega)}.
\end{equation}
In this paper we will use splines to ensure that condition
\eref{appsrop3} is satisfied (observe that it requires the
quasi-uniformity of the (coarse) mesh, i.e. a bound on the aspect
ratio of the (coarse) triangles).

For $t\in (0,T)$ let us define
\begin{equation}
V_h(t):=\big\{\varphi \circ F(x,t)\, :\, \varphi \in X_h\big\}.
\end{equation}
Write $L^2\big(0,T;H^1_0(\Omega)\big)$ the usual Sobolev space
associated to the norm
\begin{equation}
\|v\|_{L^2(0,T;H^1_0(\Omega))}^2:=\int_0^T
\big\|v(.,t)\big\|_{H^1_0(\Omega)}^2\,dt.
\end{equation}
Write $Y_T$ the subspace of $L^2\big(0,T;H^1_0(\Omega)\big)$ such
that for each $v \in Y_T$ and $t\in [0,T]$, $x\rightarrow v(x,t)$
belongs to $V_h(t)$. \\Write $u_h$ the solution in $Y_T$ of the
following system of ordinary differential equations:

\begin{equation}\label{ghjbbfh52}
\begin{cases}
(\psi,\partial_t u_h)_{L^2(\Omega)}+a[\psi,u_h]=(\psi,g)_{L^2(\Omega)}\quad \text{for all $t\in (0,T)$ and $\psi \in V_h(t)$} \\
u_h(x,0)=0.
\end{cases}
\end{equation}
Write
\begin{equation}\label{gaz52}
a[v,w]:=\int_\Omega {^t\nabla v(x,t)}a(x,t)\nabla w(x,t)\,dx.
\end{equation}

\subsubsection{Time independent domain}
 We have the following theorem
\begin{Theorem}\label{slsakdeswj8kksjd23}
Assume that $\partial_t a \equiv 0$, $\Omega$ is convex,
$\beta_{\sigma}<1$ and $(\Tr[\sigma])^{-1}\in L^\infty(\Omega_T)$
then
\begin{equation}\label{sjswhsawdskjdkjxwsdlswe}
\big\|(u-u_h)(.,T)\big\|_{L^2(\Omega)}+\big\|u-u_h\big\|_{L^2(0,T;H^1_0(\Omega))}
\leq C h \|g\|_{L^2(\Omega_T)}.
\end{equation}
\end{Theorem}
\begin{Remark}
The constant $C$ depends on $C_X$, $n$, $\Omega$,
$\lambda_{\min}(a)$ and\\
$\big\|(\Tr[\sigma])^{-1}\big\|_{L^\infty(\Omega_T)}$. If $n\geq 5$
it also depends on $\big\|\Tr[\sigma]\big\|_{L^\infty(\Omega_T)}$
and if $n=1$ it also depends on $\lambda_{\max}(a)$.
\end{Remark}

\subsubsection{Medium with a continuum of time scales.}

\begin{Theorem}
Assume that $\Omega$ is convex, and condition \ref{slkssjk88271} is
satisfied then
\begin{equation}\label{ghhjaawssjncssbbhfhs52az}
\begin{split}
\big\|(u-u_h)(T)\big\|_{L^2(\Omega)}+\big\|u-u_h\big\|_{L^2(0,T,H^1_0(\Omega))}\leq
C h \|g\|_{L^2(\Omega_T)}.
\end{split}
\end{equation}
\end{Theorem}
\begin{Remark}
The constant $C$ depends on $C_X$, $n$, $\Omega$, $\delta$ and
$\epsilon$, $\lambda_{\min}(a)$ and $\lambda_{\max}(a)$.
\end{Remark}

The system of ordinary differential equations \eref{ghjbbfh52} is
still characterized by a continuum of time scales in situations
where the entries of $a$ merely belong to $L^\infty(\Omega_T)$. They
need to be discretized  (homogenized) in time in order to be solved
numerically. This will be the object of the next subsection. Loosely
speaking, although \eref{ghjh52} is associated to a fine
tessellation and fine time steps, it is possible to approximate its
operator on a coarse tessellation with coarse time steps.

\subsection{Homogenization in space and time.}\label{jksjhs89b}
Let $M\in \N^*$. Let $(t_n=n \frac{T}{M})_{0\leq n\leq M}$ be a
discretization of $[0,T]$. Let $(\varphi_i)$ be a basis of $X_h$.
Write $Z_T$ the subspace of $Y_T$ such that $w\in Z_T$ if and only
if $w$ can be written
\begin{equation}
w(x,t)=\sum_{i} c_i(t) \varphi_i(F(x,t)).
\end{equation}
and the functions $t\rightarrow c_i(t)$ are constants on each
intervals $(t_n,t_{n+1}]$. Write $V$ the subspace of $Y_T$ such that
its elements $\psi$ can be written
\begin{equation}
\psi(x,t)=\sum_{i} d_i \varphi_i(F(x,t)).
\end{equation}
where the parameters $d_i$ are constants (on $[0,T]$). For $w\in
Y_T$, define $w_n \in V$ by
\begin{equation}
w_n(x,t):=\sum_{i} c_i(t_n) \varphi_i(F(x,t)).
\end{equation}
Write $v$ the solution in $Z_T$ of the following system of implicit
ordinary differential equations (such that $v(x,0)\equiv 0$): for
$n\in \{0,\ldots,M-1\}$ and $\psi\in V$,
\begin{equation}\label{ghjddwsdcszdbbsfh52}
\begin{split}
 \big(\psi(t_{n+1}),
v_{n+1}(t_{n+1})\big)_{L^2(\Omega)}=&\big(\psi(t_{n}),
v_{n}(t_{n})\big)_{L^2(\Omega)}\\&+\int_{t_n}^{t_{n+1}}\Big(\big(\partial_t\psi(t),
v_{n+1}(t)\big)_{L^2(\Omega)}\\&-a\big[\psi(t),v_{n+1}(t)]+\big(\psi(t),g(t)\big)_{L^2(\Omega)}
\Big)\,dt.
\end{split}
\end{equation}

The following theorem shows the stability of the implicit scheme
\eref{ghjddwsdcszdbbsfh52}.
\begin{Theorem}\label{sshjhsjjs823}
Let $v\in Z_T$ be the solution of \eref{ghjddwsdcszdbbsfh52}. We
have
\begin{equation}\label{ghnnjddsdmmsswsdcszddsbbsfh52}
\begin{split}
  \big\|v(T)\big\|_{L^2(\Omega)}+\|v\|_{L^2(0,T,H^1_0(\Omega))}
  \leq C \|g\|_{L^2(\Omega_T)}.
  \end{split}
\end{equation}
\end{Theorem}
\begin{Remark}
The constant $C$ depends on $n$, $\Omega$ and $\lambda_{\min}(a)$.
\end{Remark}

The following theorem gives an error bound on the accuracy of time
discretization scheme \eref{ghjddwsdcszdbbsfh52} when $a$ does not
depend on time.
\begin{Theorem}\label{jskjdhdcxjdh723}
Let $v\in Z_T$ be the solution of \eref{ghjddwsdcszdbbsfh52} and
$u_h$ be the solution of \eref{ghjbbfh52}. Assume that $\partial_t
a\equiv 0$. We have
\begin{equation}\label{gazoszwscxbbfmbh52}
\begin{split}
 \big\| (u_h-v)(T)\big\|_{L^2(\Omega)}+
&\|u_h-v\|_{L^2(0,T,H^1_0(\Omega))} \leq C |\Delta t|
\\&\Big(\|\partial_t
g\|_{L^2(0,T,H^{-1}(\Omega))}+\big\|g(.,0)\big\|_{L^2(\Omega)}\Big).
\end{split}
\end{equation}
\end{Theorem}
\begin{Remark}
The constant $C$ depends on $n$, $\Omega$ and $\lambda_{\min}(a)$.
\end{Remark}

The following theorem gives an error bound on the accuracy of the
time discretization scheme \eref{ghjddwsdcszdbbsfh52} when $a$ has
no bounded time derivatives.
\begin{Theorem}\label{skjhseeeddfnmnddsf23}
Assume that $\Omega$ is convex, and condition \ref{slkssjk88271} is
satisfied. Let $v\in Z_T$ be the solution of
\eref{ghjddwsdcszdbbsfh52} and $u_h$ be the solution of
\eref{ghjbbfh52}, we have
\begin{equation}\label{ghjdzmmssasdssabfh52}
\begin{split}
\big\| (u_h-v)(T)\big\|_{L^2(\Omega)}+
\|u_h-v\|_{L^2(0,T,H^1_0(\Omega))} \leq C \frac{|\Delta t|}{h}
 \|g\|_{L^2(\Omega_T)}
\end{split}
\end{equation}
where $C$ depends on $\Omega$, $n$, $\delta$, $\epsilon$,
$\lambda_{\min}(a)$ and $\lambda_{\max}(a)$.
\end{Theorem}
\begin{Remark}
Observe that the accuracy of the time discretization scheme
\eref{ghjddwsdcszdbbsfh52} requires that $|\Delta t|<<h$ when $a$
has no bounded time derivatives.
\end{Remark}

We refer to section \ref{ksjsskhs821} for numerical experiments.

\subsection{Literature and further remarks.}
For early works on homogenization with random mixing coefficients we
refer to \cite{PaVa83}, \cite{KiVa86}, \cite{PaVa79},
\cite{MR0383530}, \cite{ZhKoOl79}, \cite{Ko85}, \cite{Ko87},
\cite{MR1429379}. Papanicolaou and Varadhan \cite{MR0461684} have
considered a two-component Markov process $(x(t), y(t))$ where
$y(t)$ is rapidly varying (and is not assumed to be ergodic in
dimension one) and enters in the coefficients of the stochastic
process driving $x(t)$. They have studied the convergence properties
of $x(t)$ as the fluctuations of $y(t)$ becomes more rapid using the
martingale approach to diffusion, developed by Stroock and Varadhan
\cite{MR0359025}, \cite{MR0359024}, \cite{MR532498},
\cite{MR0410912}.

The numerical homogenization method implemented in this paper is a
finite element method. The idea of using oscillating tests functions
can be  back tracked to the work of Murat and Tartar on
homogenization and H-convergence, we refer in particular to
\cite{MR557520} and \cite{MR1493039}. Those papers also contain
convergence proofs for the finite element method in an abstract
setting for a sequence of $H$-converging elliptic operators (recall
that the framework of H-convergence is independent from ergodicity
or scale separation assumptions and are based on the compactness of
any sequence of solutions of $-\diiv a_\epsilon \nabla u_\epsilon=g$
with uniformly bounded and elliptic conductivities $a_\epsilon$, we
also refer to the initial work of Spagnolo \cite{MR0240443} for
G-convergence).

The numerical implementation and practical application of
oscillating test functions in numerical finite element
homogenization have been called multi-scale finite element methods
and have been studied by several authors \cite{MR1286212},
\cite{MR1740386}, \cite{MR1613757}, \cite{MR1455261},
\cite{MR2123115}, \cite{MR1232956}, \cite{MR1194543}, \cite{AlBr05}.
The work of Hou and Wu \cite{MR1455261} has been a large source of
inspiration in numerical applications (particularly for reservoir
modeling in geophysics, we refer to \cite{MR1898136},
\cite{MR1956022}, \cite{MR2111701} and \cite{HouEf05} for recent
developments) since it was leading to a coarse scale operator while
keeping the fine scale structures of the solutions. With the
 method introduced by Hou and Wu, the
construction of the base functions is decoupled from element to
element leading to a scheme adapted to parallel computers. A proof
of the convergence of the method is given in periodic settings when
the size of the heterogeneities is smaller than the grid size and an
'oversampling technique' is proposed to remove the so called cell
resonance error \cite{MR2119937} when the size of the
heterogeneities is comparable to the grid size.

Allaire and Brizzi \cite{AlBr04} have observed that multiscale
finite element method with splines would have a higher accuracy and
have introduced the composition rule (we also refer to
\cite{MR1286212}). In \cite{OwZh05}, it has been observed that if
$u$ is the solution of the divergence form elliptic equation
\begin{equation}\label{ghjh5nsx2}
\begin{cases}
-\diiv\big(a(x)\nabla u(x)\big)=g\quad \text{in}\quad \Omega \\
u=0 \quad \text{in}\quad \partial \Omega.
\end{cases}
\end{equation}
and $F$ are harmonic coordinates defined by
\begin{equation}\label{dgdgfsghsf6bg2}
\begin{cases}
\diiv a \nabla F=0 \quad \text{in}\quad \Omega\\
F(x)=x \quad \text{on}\quad \partial \Omega.
\end{cases}
\end{equation}
then under the Cordes type condition $\beta_{\sigma}<1$ on $\sigma$
given by \eref{jhscczxhjd}, one has for some $p>2$.
\begin{equation}
\|u\circ F^{-1}\|_{W^{2,p}(\Omega)}\leq C \|g\|_{L^p(\Omega)}.
\end{equation}
It has been deduced from this compensation phenomenon that numerical
homogenization methods based on oscillating finite elements can
converge in the presence of a continuum of scales if one uses global
harmonic coordinates  to obtain the test functions instead of
solutions of a local cell problem \cite{OwZh05}. In dimension three
and higher it has been known since the work of Fenchenko and
Khruslov \cite{MR615994}, \cite{MR1145750} that the homogenization
of divergence form elliptic operators $-\diiv a_\epsilon \nabla
u_\epsilon=g$ can lead to a non local homogenized operator if the
sequence of matrices $a_\epsilon$ is uniformly elliptic but with
entries uniformly bounded only in $L^1(\Omega)$. From a numerical
point of view this non-local effects imply that a nonlocal numerical
homogenization method cannot be avoided to obtain accuracy. Hence in
 \cite{OwZh05}, it is shown that the accuracy of local methods
 depend on the aspect ratio of the triangles of the tessellation
 with respect to harmonic coordinates (which is not the case if one uses
 non local finite elements, we refer to \cite{OwZh05}
 for further discussions on the apparition of non local effects
 in numerical homogenization). Recently Briane has shown \cite{Bri05}
that this non-local effect is absent in dimension two in the
H-convergence setting.

The phenomenon is similar here, however observe that if one has
solved the initial parabolic equation at least $n$ times and those
solutions are (locally) linearly independent it is also possible to
use them as new coordinates for numerical homogenization. Observe
that in dimension higher than three the harmonic coordinates are not
always invertible, an idea to bypass this difficulty could be either
to choose the change of coordinates locally and adaptively  or to
enrich the coordinates by writing down the initial equations as
degenerate equations in a space of higher dimension \cite{Var05},
these points have not been explored. For divergence form elliptic
equations, recall that fast methods based on hierarchical
matrices\footnote{As for the fast multipole method and the
hierarchical multipole method designed by L. Greengard and V.
Rokhlin \cite{MR918448}, these methods are based on the singular
value decomposition of operators Green's function.} are available
\cite{MR1993936, BeCh05, Beb05, Bebe05, Beben05} for solving
\eref{ghjh5nsx2} and \eref{dgdgfsghsf6bg2} in $O\big(N (\ln
N)^{n+3}\big)$ operations ($N$ being the number of interior nodes of
the fine mesh).

The issue of numerical homogenization partial differential equations
with heterogeneous coefficients has received a great deal of
attention and many  methods have been proposed. A few of them are
cited below.
\begin{itemize}
\item Multi-scale finite volume methods \cite{JLP03}.
\item Heterogeneous Multi-scale Methods \cite{EV04}, \cite{MR2164241}.
\item Wavelet based homogenization \cite{MR1492791}, \cite{MR1618846}, \cite{LC04}, \cite{MR1614457}, \cite{MR1614980}, \cite{MR1354913}.
\item Residual free bubbles methods \cite{MR2006324}.
\item Discontinuous enrichment methods \cite{MR2007030}, \cite{MR1870426}.
\item Partition of Unity Methods \cite{FY05}.
\item Energy Minimizing Multi-grid Methods \cite{MR1756048}.
\end{itemize}

Following the methods of \cite{OwZh05}, it is possible to implement
a finite-volume method based on the compensation theorems given in
this paper. The elements given in this paper contain the fine scale
structure of $F$, as it has been done in \cite{OwZh05}, it is
possible to approximate the initial parabolic operator by a
homogenized parabolic operator associated to the coarse mesh (the
test functions in this case would be piecewise linear on the coarse
mesh and the approximation error associated to the homogenized
operator would depend on the aspect ratio of the triangles of the
coarse mesh in the metric induced by $F$).

Finally, in this paper $a$  has been assumed to be bounded and
uniformly elliptic. Without these assumptions the diffusion
associated to homogenized operator can be anomalously slow
\cite{BeOw00b}, \cite{Ow00a} or fast (super-diffusive) \cite{Ow04}.
If $a$ has an unbounded skew symmetric component, the homogenization
of \eref{ghjh52} can give rise to a degenerate operator
 \cite{Ow04}.

\section{Proofs}

\subsection{Compensation.}

\subsubsection{Time independent medium.}
We will need the following lemmas. Let $\A_T$ be the bilinear form
on $L^2\big(0,T;H^1_0(\Omega)\big)$ defined by
\begin{equation}
\A_T[v,w]:=\int_0^T a[v,w](t) \,dt
\end{equation}
where
\begin{equation}
a[v,u](t):=\int_\Omega {^t\nabla v}(x,t)a(x,t)\nabla u(x,t)\,dx.
\end{equation}
We write $\A_T[u]:=\A_T[u,u]$.
\begin{Lemma}\label{sjhsgh8}
We have
\begin{equation}\label{eqra1}
\big\|u(.,T)\big\|_{L^2(\Omega)}^2+ \A_T[u]\leq
\frac{C_{n,\Omega}}{\lambda_{\min}(a)} \|g\|_{L^2(\Omega_T)}^2.
\end{equation}
\end{Lemma}
\begin{proof}
Multiplying \eref{ghjh52}by $u$ and integrating with respect to time
we obtain that
\begin{equation}\label{gsesfdssw2}
\frac{1}{2}\big\|u(.,T)\big\|_{L^2(\Omega)}^2+
\A_T[u]=(u,g)_{L^2(\Omega_T)}.
\end{equation}
Using Poincar\'{e} and Minkowksi inequalities leads us to
\eref{eqra1}.
\end{proof}

\begin{Lemma}\label{jshsszskjsh61RED}
Assume $\partial_t a\equiv 0$. We have
\begin{equation}\label{sxdsgxssssxsesw2RED}
\begin{split}
\big\|\partial_t u\big\|_{L^2(\Omega_T)}^2+a\big[u(.,T)\big]\leq
 \big\|g\big\|_{L^2(\Omega_T)}^2.
 \end{split}
\end{equation}
\end{Lemma}
\begin{proof}

Multiplying \eref{ghjh52} by $\partial_t u$ and integrating by parts
we obtain that
\begin{equation}\label{ghjbsanjxzassxxsxch52RED}
\begin{split}
\big\|\partial_t u(.,t)\big\|_{L^2(\Omega)}^2+a[\partial_t
u,u]=(\partial_t u, g)_{L^2(\Omega)}.
\end{split}
\end{equation}
Observing that
\begin{equation}\label{ghjbsanjxzasswsxxsxch52RED}
\begin{split}
a[\partial_t u,u]=\frac{1}{2}\partial_t
\big(a[u]\big)-\frac{1}{2}\int_{\Omega}{^t\nabla u \partial_t a
\nabla u}.
\end{split}
\end{equation}
we conclude by  integration with respect to time and using Minkowski
inequality.
\end{proof}

\begin{Lemma}\label{jshsszskjsh61}
Assume $\partial_t a \equiv 0$. We have
\begin{equation}\label{sxdsgxssssxsesw2}
\begin{split}
\big\|\partial_t u(.,T)\big\|_{L^2(\Omega)}^2+\A_T[\partial_t u]\leq
& \frac{C_{n,\Omega}}{\lambda_{\min}(a)} \|\partial_t
g\|_{L^2(0,T,H^{-1}(\Omega))}^2+\big\|g(.,0)\big\|_{L^2(\Omega)}^2.
\end{split}
\end{equation}
\end{Lemma}
\begin{proof}
We obtain from \eref{ghjh52} that
\begin{equation}\label{gchudyhjh52}
\begin{split}
\partial_t^2 u=\diiv\big(a(x,t)\nabla \partial_t u(x,t)\big)+\diiv\big(\partial_t a(x,t)\nabla  u(x,t)\big) +\partial_t
g.
\end{split}
\end{equation}
Multiplying \eref{gchudyhjh52} by $\partial_t u$ and integrating
with respect to time we obtain that
\begin{equation}\label{ghjbsanjxzassxxsxch52p}
\begin{split}
\frac{1}{2}\big\|\partial_t
u(.,T)\big\|_{L^2(\Omega)}^2+\A_T[\partial_t u]=&\int_0^T
(\partial_t u,\partial_t g)_{L^2(\Omega)}dt-\int_0^T (\partial_t
a)\big[\partial_t u, u\big]\,dt\\&+\frac{1}{2}\big\|\partial_t
u(.,0)\big\|_{L^2(\Omega)}^2.
\end{split}
\end{equation}
We conclude by the $H^{-1}$-duality inequality and Minkowski
inequality.
\end{proof}

 We now need a variation of Campanato's result \cite{CM5}
on non-divergence form elliptic operators. Let us write for a
symmetric matrix $M$,
\begin{equation}
\nu_M:=\frac{\Tr(M)}{\Tr({^tM M})}.
\end{equation}
We consider the following Dirichlet problem:
\begin{equation}\label{dcaqssaslkwaq21}
L_M v=f
\end{equation}
with $L_M:=\sum_{i,j=1}^n M_{ij}(x) \partial_i \partial_j$. The
following  theorems \ref{hdgjhsdswasgd7} and \ref{hdgjhdgd7} are
straightforward adaptations of theorem 1.2.1 of
 \cite{MPG00}. They are proven in  \cite{MPG00} under the assumption that $M$ is bounded and
 elliptic. It is easy to check that the conditions $\beta_{M}<1$ and
 $\nu_M\in L^\infty(\Omega)$ are sufficient for the validity of
 those theorems. We refer to \cite{OwZh05} for that adaptation.

\begin{Theorem}\label{hdgjhsdswasgd7}
Assume that $\beta_M<1$,  $\nu_M\in L^\infty(\Omega)$ and $\Omega$
is convex. Then if $f\in L^2(\Omega)$ the Dirichlet problem
\eref{dcaqssaslkwaq21} has a unique solution satisfying
\begin{equation}\label{hdhazdgc7}
\|v\|_{W^{2,2}_D(\Omega)}\leq \frac{C }{1-\beta_M^\frac{1}{2}}
\|\nu_M f\|_{L^2(\Omega)}.
\end{equation}
\end{Theorem}
\begin{Remark}
$\beta_M$ is the Cordes parameter  associated to $M$.
\end{Remark}
\begin{Theorem}\label{hdgjhdgd7}
Assume that $\beta_M<1$,  $\nu_M\in L^\infty(\Omega)$ and $\Omega$
is convex. Then, there exists a real number $p_0>2$ depending only
on $n,\Omega$ and $\beta_M$ such that for each $f\in L^p(\Omega)$,
$2\leq p<p_0$ the Dirichlet problem \eref{dcaqssaslkwaq21} has a
unique solution satisfying
\begin{equation}\label{hdhdgc7}
\|v\|_{W^{2,p}_D(\Omega)}\leq \frac{C_{n,\Omega,p}
}{1-\beta_M^\frac{1}{2}} \|\nu_M f\|_{L^p(\Omega)}.
\end{equation}
\end{Theorem}

Let us now prove the compensation theorems. Choose
\begin{equation}\label{ksjshsj7622}
M:=\frac{\sigma}{|\det(\nabla F)|^\frac{1}{2}}\circ F^{-1}.
\end{equation}
It is easy to check that $\beta_\sigma<1$ implies that $F$ is an
homeomorphism from $\Omega$ onto $\Omega$, thus \eref{ksjshsj7622}
is well defined. Moreover observe that $\beta_M=\beta_\sigma$ and
\begin{equation}\label{hdhxcszswc7}
\|\nu_M\|_{L^\infty(\Omega_T)}^2 \leq
\frac{C_n}{(\lambda_{\min}(a))^\frac{n}{2}}
\big\|(\Tr[\sigma])^{\frac{n}{4}-1}\big\|^2_{L^\infty(\Omega_T)}.
\end{equation}
Fix $t\in [0,T]$. Choose
\begin{equation}\label{ghjhssalzz52}
\begin{split}
f:=\frac{(\partial_t u - g)}{|\det (\nabla F)|^\frac{1}{2}}\circ
F^{-1}.
\end{split}
\end{equation}
Observe that by the change of variable $y=F(x)$ one obtains that if
$\partial_t a\equiv 0$ (which implies that $F$ is time independent),
$\partial_t u \in L^2(\Omega)$ and $g(.,t)\in L^2(\Omega)$ that
$f\in L^2(\Omega)$ and
\begin{equation}\label{ghjhssadsfalzz52}
\begin{split}
\|f\|_{L^2(\Omega)}\leq \|\partial_t
u\|_{L^2(\Omega)}+\|g\|_{L^2(\Omega)}.
\end{split}
\end{equation}
It follows from theorem \ref{hdgjhsdswasgd7} that there exists a
unique $v\in W^{2,2}_D(\Omega)$ satisfying
\begin{equation}\label{hdhsszdgc7}
\|v\|_{W^{2,2}_D(\Omega)}^2\leq
\frac{C\|\nu_M\|_{L^\infty(\Omega_T)}^2
}{(1-\beta_\sigma^\frac{1}{2})^2} \big(\|\partial_t
u\|_{L^2(\Omega)}^2+\|g\|_{L^2(\Omega)}^2\big).
\end{equation}
 and the following equation
\begin{equation}\label{ghjhsxfazz52}
\begin{split}
\partial_t \hat{u}(y,t)=\sum_{i,j}\big(\sigma(F^{-1}(y,t),t)\big)_{i,j}
\partial_i\partial_j v(y,t)+\hat{g}(y,t).
\end{split}
\end{equation}
We use the notation $\hat{g}:=g\circ F^{-1}$ and $\hat{u}:=u\circ
F^{-1}$. Using the change of variable $y=F(x)$ and using the
property $\diiv a \nabla F=0$ when $\partial_t a\equiv 0$
 we obtain that
\eref{ghjhsxfazz52} can be written
\begin{equation}\label{ghjhsxfsaazsdz52j}
\begin{split}
\partial_t u=\diiv \big(a \nabla (v\circ F)\big)+g.
\end{split}
\end{equation}
If $\partial_t u \in L^2(\Omega)$ and $g(.,t)\in L^2(\Omega)$ we can
use the uniqueness property of the solution of the divergence form
elliptic Dirichlet  problem
\begin{equation}\label{ghjhsxfsaazsdz52}
\begin{split}
\diiv \big(a \nabla w\big)=\partial_t u -g.
\end{split}
\end{equation}
to obtain that $v\circ F=u$. Thus using lemma \ref{jshsszskjsh61} we
have proven theorem \ref{ksjshszjscxsxsdded8721}. Moreover assume
that $g\in L^2(\Omega_T)$ and $\partial_t u\in L^2(\Omega_T)$. It
follows that for $t\in [0,T]-B$, $g(.,t)\in L^2(\Omega)$ and
$\partial_t u(.,t)\in L^2(\Omega)$ where $B$ is a subset of $[0,T]$
of $0$-Lebesgue measure. It follows from the previous arguments that
on $[0,T]-B$, $u\circ F^{-1}(.,t) \in W^{2,2}_D(\Omega)$ and
satisfies
\begin{equation}\label{hdhsszwweeddgc7}
\|u\circ F^{-1}(.,t)\|_{W^{2,2}_D(\Omega)}^2\leq
\frac{C\|\nu_M\|_{L^\infty(\Omega_T)}^2
}{(1-\beta_\sigma^\frac{1}{2})^2} \big(\|\partial_t
u(.,t)\|_{L^2(\Omega)}^2+\|g(.,t)\|_{L^2(\Omega)}^2\big).
\end{equation}
Integrating \eref{hdhsszwweeddgc7} with respect to time we obtain
that \\$u\circ F^{-1}\in L^2(0,T, W^{2,2}_D(\Omega))$ and
\begin{equation}\label{hdhsszwxweeddgc7}
\|u\circ F^{-1}\|_{L^2(0,T, W^{2,2}_D(\Omega))}^2\leq
\frac{C\|\nu_M\|_{L^\infty(\Omega_T)}^2
}{(1-\beta_\sigma^\frac{1}{2})^2} \big(\|\partial_t
u\|_{L^2(\Omega_T)}^2+\|g\|_{L^2(\Omega_T)}^2\big).
\end{equation}
Thus using lemma \ref{jshsszskjsh61RED} we have obtained theorem
\ref{ksjshjsxxsd8721}.

Let us now prove theorem \ref{ksjshjsxgghxsd8721}. Assume that there
exists $q_0>2$ such that for $2\leq p<q_0$, $\partial_t u \in
L^p(\Omega_T)$ and $g\in L^p(\Omega_T)$. Let us now apply theorem
\ref{hdgjhdgd7} with $p<\min(p_0,q_0)$, $M$ given by
\eref{ksjshsj7622} and $f$ given by \eref{ghjhssalzz52}. It follows
that for $t\in [0,T]-B$ (where $B$ is a subset of $[0,T]$ of
$0$-Lebesgue measure), $g(.,t)\in L^p(\Omega)$ and $\partial_t
u(.,t)\in L^p(\Omega)$. We deduce from theorem \ref{hdgjhdgd7} and
the argumentation related to equation \eref{ghjhsxfsaazsdz52} that
on  $[0,T]-B$, $u\circ F^{-1}(.,t)\in W^{2,p}_D(\Omega)$ and
\begin{equation}\label{hdhssdwedeesewddgc7}
\|u\circ F^{-1}(.,t)\|_{W^{2,p}_D(\Omega)}^p\leq
\frac{C_{n,p,\Omega}\|\nu_M\|_{L^\infty(\Omega_T)}^p
}{(1-\beta_\sigma^\frac{1}{2})^p} \big(\|\partial_t
u(.,t)\|_{L^p(\Omega)}^p+\|g(.,t)\|_{L^p(\Omega)}^p\big).
\end{equation}
Integrating \eref{hdhssdwedeesewddgc7} with respect to time we
obtain that\\ $u\circ F^{-1}\in L^p(0,T,W^{2,p}_D(\Omega))$ and

\begin{equation}\label{hdhssssesewddgc7}
\|u\circ F^{-1}\|_{L^p(0,T,W^{2,p}_D(\Omega))}\leq
\frac{C_{n,p,\Omega}\|\nu_M\|_{L^\infty(\Omega_T)}
}{1-\beta_\sigma^\frac{1}{2}} \big(\|\partial_t
u\|_{L^p(\Omega_T)}+\|g\|_{L^p(\Omega_T)}\big).
\end{equation}

It remains to show that under the assumptions of theorem
\ref{ksjshjsxgghxsd8721}, $\partial_t u \in L^p(\Omega_T)$.

In order to bound $\big\|\partial_t u(.,t)\big\|_{L^p(\Omega)}$ we
use general Sobolev inequalities (chapter 5.6 of \cite{Evans97}).
\begin{itemize}
\item If $n\geq 3$, write $p^*=2n/(n-2)$. We have for $2<p\leq p^*$,
\begin{equation}\label{jskjsdh83}
\big(\int_{\Omega}(\partial_t u )^{p} \,dx\big)^\frac{2}{p} \leq
C_{n,\Omega} \big(\int_{\Omega}(\partial_t u )^{p^*}
\,dx\big)^\frac{2}{p^*}
\end{equation}
thus, using Gagliardo-Nirenberg-Sobolev inequality
\begin{equation}\label{hjsshsghs61}
\big(\int_{\Omega}(\partial_t u )^{p}  \,dx\big)^\frac{2}{p} \leq
C_{n,p,\Omega} \frac{1}{\lambda_{\min}(a)} a[\partial_t u].
\end{equation}
\item If $n=2$, we write for $(x_1,x_2,x_3)\in \Omega\times (0,1)$, $v(x_1,x_2,x_3):=\partial_t
u(x_1,x_2)$. Using Gagliardo-Nirenberg-Sobolev inequality in
dimension three we obtain that for $2<p\leq 6$
\begin{equation}
\big(\int_{\Omega}(\partial_t u )^{p} \,dx\big)^\frac{2}{p} \leq
C_{n,p,\Omega} \int_{\Omega}(\nabla \partial_t  u )^{2} \,dx.
\end{equation}
Which leads us to \eref{hjsshsghs61}.

\item If $n=1$ then using Morrey's inequality we obtain that
with $\gamma:=1/2$,
\begin{equation}
\|\partial_t u\|_{C^{0,\gamma}(\Omega)}^2 \leq C_{\Omega}
\frac{1}{\lambda_{\min}(a)} a[\partial_t u].
\end{equation}
\end{itemize}
We conclude the proof of theorem \ref{ksjshjsxgghxsd8721} by using
lemma \ref{jshsszskjsh61}.

We deduce theorem \ref{hdgwsssawwsd7} from Morrey's inequality and
theorem \ref{ksjshjsxgghxsd8721}.

\paragraph{H\"{o}lder continuity for $n\geq 3$ or non-convexity of
$\Omega$.} In this paragraph we will not assume $\Omega$ to be
convex. Let $N^{p,\lambda}(\Omega)$ $(1<p<\infty,\, 0<\lambda<n$) be
the weighted Morrey space formed by the functions
$v:\Omega\rightarrow \R$ such that
$\|v\|_{N^{p,\lambda}(\Omega)}<\infty$ with
\begin{equation}
\|v\|_{N^{p,\lambda}(\Omega)}:=\sup_{x_0\in
\Omega}\Big(\int_{\Omega}|x-x_0|^{-\lambda} |v(x)|^p
\Big)^\frac{1}{p}.
\end{equation}
To obtain the H\"{o}lder continuity of $u\circ F^{-1}$ in dimension
$n\geq 3$ we use corollary 4.1 of \cite{MR1903306}. We give the
result of S. Leonardi below in a form adapted to our context.
Consider the Dirichlet problem \eref{dcaqssaslkwaq21}. We do not
assume $\Omega$ to be bounded. We write $W^{2,p,\lambda}(\Omega)$
the functions in $W^{2,p}_D(\Omega)$ such that their second order
derivatives belong to $N^{p,\lambda}(\Omega)$.
\begin{Theorem}\label{ksjhs721}
There exist a constant $C^*=C^*(n,p,\lambda,\partial \Omega)>0$ such
that if $\beta_M <C^*$ and $f\in N^{p,\lambda}(\Omega)$ then the
Dirichlet problem \eref{dcaqssaslkwaq21} has a unique solution in
$W^{2,p,\lambda}\cap W^{1,p}_0(\Omega)$. Moreover, if $0<\lambda
<n<p$ then $\nabla v \in C^{\alpha}(\Omega)$ with $\alpha=1-n/p$ and
\begin{equation}
\|\nabla v\|_{C^\alpha(\Omega)}\leq
\frac{C}{\lambda_{\min}(M)}\|f\|_{N^{p,\lambda}(\Omega)}
\end{equation}
where $C=C(n,p,\lambda,\partial \Omega)$.
\end{Theorem}

The proof of theorem \ref{hdgwsswdswsd7}  is an application of
theorem \ref{ksjhs721}. We just need to observe that from H\"{o}lder
inequality we have for $0<\epsilon<0.5$
\begin{equation}
\|f\|_{N^{p,\epsilon}(\Omega)} \leq C_{n,p,\Omega,\epsilon}
\|f\|_{L^{p (1+\epsilon)}(\Omega)}.
\end{equation}
From this point the proof is similar to the proof of theorem
\ref{ksjshjsxgghxsd8721}.

\subsubsection{Medium with a continuum of time scales.}
We will need theorems 1.6.2 and 1.6.3 of \cite{MPG00}. For the sake
of completeness we will remind those theorems below in version
adapted to our framework. Consider the following parabolic problem:
\begin{equation}\label{dcaqssssswaslkwaq21}
\partial_t v=\sum_{i,j=1}^n M_{ij}(x) \partial_i \partial_j v+ f.
\end{equation}
We assume $M$ to be symmetric bounded and elliptic and $v=0$ at
$t=0$ and on the boundary $\partial \Omega$. Write
\begin{equation}
\eta_M:=\sup_{x\in \Omega_T}
\frac{\Tr[{^tMM}]+1}{\big(\Tr[M]+1\big)^2}.
\end{equation}
and
\begin{equation}
\alpha_M:=\sup_{x\in \Omega_T} \frac{\Tr[M]+1}{\Tr[{^tMM}]+1}.
\end{equation}
Write for $p\geq 2$
\begin{equation}
S_p(\Omega_T):=\Big\{v\in L^p\big(0,T,W^{2,p}_D(\Omega)\big);
\partial_t v \in L^p(\Omega_T); v(.,0)\equiv 0\Big\}
\end{equation}
and
\begin{equation}
\|v\|_{S_p(\Omega_T)}^p:=\int_{\Omega_T}\big(\sum_{i,j}(\partial_i
\partial_j v)^2+(\partial_t v)^2\big)^\frac{p}{2}\,dy\,dt.
\end{equation}

\begin{Theorem}\label{hgft5411}
Assume $\Omega$ to be convex and that there exists $\epsilon
\in(0,1)$ such that $\eta_M \leq 1/(n+\epsilon)$, then for each
$f\in L^2(\Omega_T)$ the Cauchy-Dirichlet problem
\eref{dcaqssssswaslkwaq21} admits a unique  solution in
$S_2(\Omega_T)$ which satisfies the bound
\begin{equation}
\|v\|_{S_2(\Omega_T)}\leq
\frac{\alpha_M}{1-\sqrt{1-\epsilon}}\|f\|_{L^2(\Omega_T)}.
\end{equation}
\end{Theorem}

\begin{Theorem}\label{hghjg6665542212}
Assume $\Omega$ to be convex and that there exists $\epsilon
\in(0,1)$ such that $\eta_M \leq 1/(n+\epsilon)$, then there exists
a number $p_0>2$ depending on $\Omega,n,\epsilon$ such that for each
$f\in L^p(\Omega_T)$ the Cauchy-Dirichlet problem
\eref{dcaqssssswaslkwaq21} admits a unique  solution in
$S_p(\Omega_T)$ which satisfies the bound
\begin{equation}\label{jshshsg651}
\|v\|_{S_p(\Omega_T)}\leq C_p
\frac{\alpha_M}{1-\sqrt{1-\epsilon}}\|f\|_{L^p(\Omega_T)}.
\end{equation}
\end{Theorem}
\begin{Remark}
In fact theorem 1.6.3 of \cite{MPG00} is written with
$1-C(p)\sqrt{1-\epsilon}$ in the denominator of \eref{jshshsg651}
but it is easy to modify it to obtain \eref{jshshsg651} by lowering
the value of $p_0$ changing the value of $C_p$.
\end{Remark}
Let $\delta >0$. Let us now apply theorem \ref{hgft5411} on
$[0,T/\delta ]$ with
\begin{equation}
M:=\delta \sigma \circ F^{-1}(y,\delta t)
\end{equation}
and
\begin{equation}
f:=\delta (g\circ F^{-1})(y,\delta t).
\end{equation}
Observe that if condition \ref{slkssjk88271} is satisfied then $F$
is an homeomorphism and $M$ is well defined, bounded and elliptic.
Moreover $\eta_M<\infty$ and $\alpha_M<\infty$ since
\begin{equation}
\esssup_{\Omega_{\frac{T}{\delta}}}\frac{\Tr[{^tMM}]+1}{\big(\Tr[M]+1\big)^2}=\esssup_{\Omega_T}
\frac{\delta^2\Tr[{^t\sigma\sigma}]+1}{\big(\delta\Tr[\sigma]+1\big)^2}\,.
\end{equation}
It follows  that the following equation admits a unique solution in
$S_2(\Omega_{\frac{T}{\delta}})$.
\begin{equation}
\begin{split}
\partial_t w(y,t)=\sum_{i,j}M_{i,j}(y,t)
\partial_i\partial_j w(y,t)+k(y,t)
\end{split}
\end{equation}
with $k(y,t)=\delta \hat{g}(y,\delta t)$. And we have

\begin{equation}
\int_0^\frac{T}{\delta}\int_{\Omega} \big((\partial_t
w)^2+\sum_{i,j}(\partial_i
\partial_j w)^2\big)\,dy\,dt\leq
\frac{C}{(1-\sqrt{1-\epsilon})^2}\|f\|_{L^2(\Omega_{\frac{T}{\delta}})}.
\end{equation}
Using the change of variables $t\rightarrow \delta t$ and writing
\begin{equation}
w(y,t):=v(y,\delta t).
\end{equation}
we obtain that $v$ satisfies the following equation on $\Omega_T$
\begin{equation}\label{ghjhdrsazz52}
\begin{split}
\partial_t v(y,t)=\sum_{i,j}\big(\sigma(F^{-1}(y,t),t)\big)_{i,j}
\partial_i\partial_j v(y,t)+\hat{g}(y,t).
\end{split}
\end{equation}
Using the change of variable $y=F(x)$ and observing that $\partial_t
F=\diiv a\nabla F$ we obtain that $v\circ F$ satisfies
\begin{equation}\label{ghjhdrsdwazdz52}
\begin{split}
\partial_t (v\circ F)=\diiv \big(a\nabla (v\circ F)\big)+g.
\end{split}
\end{equation}
It follows from the uniqueness of the solution of
\eref{ghjhdrsdwazdz52} that $u=v\circ F$. In resume we have obtained
theorem \ref{skjhsj823} (we use lemma \ref{dkdjjddh} to control the
constants). The proof of \ref{skjhsjd823} is similar and based on
theorem \ref{hghjg6665542212}. The proof of \ref{skjhsssej823}
follows from \ref{skjhsjd823} and Morrey's inequality.

Let us now prove proposition \ref{kwjdjwj7}. Write $x=\Tr[\sigma]$
and $z=n \frac{\Tr[{^t\sigma \sigma}]}{(\Tr[\sigma])^2}$ (observe
that $1\leq z\leq n$). It is easy to check that condition
\ref{slkssjk88271} can be written
\begin{equation}\label{kjshs71}
-\delta^2 x^2 (\frac{\epsilon +n}{n}z-1)+2x\delta-(n+\epsilon-1)\geq
0.
\end{equation}
Choose $\delta=n
\big\|(\Tr[\sigma])^{-1}\big\|_{L^\infty(\Omega_T)}$. Observing that
$\delta x \geq n$ and $\delta x \leq n y_\sigma$
 it is easy to conclude the proof of proposition \ref{kwjdjwj7}.
Similarly obtains the following lemma by straightforward computation
from equation \eref{kjshs71}.
\begin{Lemma}\label{dkdjjddh}
Assume that condition \ref{slkssjk88271} is satisfied then
$\mu_\sigma<C(n,\epsilon,\delta)$
\begin{equation}
\big\|(\Tr[\sigma])^{-1}\big\|_{L^\infty(\Omega_T)}\leq
C(n,\epsilon,\delta)
\end{equation}
and
\begin{equation}
\big\|\Tr[\sigma]\big\|_{L^\infty(\Omega_T)}\leq
C(n,\epsilon,\delta).
\end{equation}
\end{Lemma}

\subsection{Convergence of the finite element method.}\label{PreLem}
 Write $\Rh$ the projection operator
mapping $L^2\big(0,T;H^1_0(\Omega)\big)$ onto $Y_T$ defined by: for
all $v\in Y_T$
\begin{equation}
\A_T[v,u-\Rh u]=0.
\end{equation}
Write $\rho:=u-\Rh u$ and $\theta:=\Rh u-u_h$.

\begin{Lemma}\label{jkshssde732}
\begin{equation}\label{ghhjaawjszsbbhfh52az}
\frac{1}{2}\big\|(u-u_h)(T)\big\|^2_{L^2(\Omega)}+\A_T[u-u_h]=\int_{\Omega_T}\rho\partial_t(u-u_h)+\A_T[\rho,u-u_h].
\end{equation}
\end{Lemma}
\begin{proof}
Subtracting
 \eref{ghjh52} (integrated against $\psi$) and \eref{ghjbbfh52} we obtain that
\begin{equation}\label{ghjjjbaabbfh52}
\big(\psi,\partial_t (u-u_h)\big)+a[\psi,u-u_h]=0\quad \text{for all
$\psi \in V_h(t)$}.
\end{equation}
Integrating by parts with respect to time we deduce that
\begin{equation}\label{ghjjjbssabfh52}
\big(\psi, (u-u_h)(.,t)\big)+a[\psi,u-u_h]=\int_{\Omega_t}
\partial_t \psi (u-u_h).
\end{equation}
Taking $\psi=\theta$ in \eref{ghjjjbssabfh52} we deduce that
\begin{equation}\label{ghjjjbssabsasafh52}
\begin{split}
\big\|(u-u_h)(.,t)\big\|^2_{L^2(\Omega)}+\A_t[u-u_h]=&\int_{\Omega_t}
\partial_t \theta (u-u_h)+\big(\rho, (u-u_h)(.,t)\big)\\&+\A_t[\rho,u-u_h].
\end{split}
\end{equation}
Observing that
\begin{equation}\label{ghjjjsbssabsasdafh52}
\int_0^t \big(\partial_t \theta, u-u_h\big)+\big(\rho,
(u-u_h)(.,t)\big)=\frac{1}{2}\big\|(u-u_h)(.,t)\big\|^2_{L^2(\Omega)}+\int_0^t(
\rho, \partial_t (u-u_h)).
\end{equation}
we deduce the lemma.
\end{proof}

\subsubsection{Time independent medium.}

\begin{Lemma}\label{jkshssssde7a32}
\begin{equation}\label{ghhjaawjbbhfhs52az}
\begin{split}
\big\|(u-u_h)(T)\big\|^2_{L^2(\Omega)}+\A_T[u-u_h]\leq  2\Big(&
\|\rho\|_{L^2(\Omega_T)}\|\partial_t u-\partial_ t
u_h\|_{L^2(\Omega_T)}\\&+ \A_T[\rho]\Big).
\end{split}
\end{equation}
\end{Lemma}
\begin{proof}
Lemma \ref{jkshssssde7a32} is a straightforward consequence of lemma
\ref{jkshssde732} and Cauchy-Schwartz and Minkowski inequalities.
\end{proof}

\begin{Lemma}\label{sjhsgh7}
We have
\begin{equation}\label{gsessnbx6fdssw2}
\big\|u_h(.,T)\big\|_{L^2(\Omega)}^2+ \A_T[u_h]\leq
\frac{C_{n,\Omega}}{\lambda_{\min}(a)} \|g\|_{L^2(\Omega_T)}^2.
\end{equation}
\end{Lemma}
\begin{proof}
Taking $\psi=u_h$ in \eref{ghjbbfh52} and integrating with respect
to time we obtain that
\begin{equation}\label{gsesfdssmmw2}
\frac{1}{2}\big\|u_h(.,T)\big\|_{L^2(\Omega)}^2+
\A_T[u_h]=(u_h,g)_{L^2(\Omega_T)}.
\end{equation}
Using Poincar\'{e} and Minkowksi inequalities leads us to
\eref{gsessnbx6fdssw2}.
\end{proof}

\begin{Lemma}\label{jshsszskssjsbh61RED}
Assume $\partial_t a \equiv 0$. We have
\begin{equation}\label{sxdsgxssssbxsesw2RED}
\begin{split}
\big\|\partial_t u_h\big\|_{L^2(\Omega_T)}^2+a\big[u_h(.,T)\big]\leq
 \big\|g\big\|_{L^2(\Omega_T)}^2.
 \end{split}
\end{equation}
\end{Lemma}
\begin{proof}
The proof is similar to lemma \ref{jshsszskjsh61RED}. We need to
take $\psi=\partial_t u_h$ in \eref{ghjbbfh52}.
\end{proof}

Let $t\in [0,T]$ and $v\in H^1_0(\Omega)$,   we will write $\Rht
v(.,t)$ the solution of:
\begin{equation}
\int_\Omega {^t\nabla \psi}a(x,t)(\psi,v-\Rht v)\,dx=0 \quad
\text{for all $\psi\in V_h(t)$}.
\end{equation}

We will need the following lemma,
\begin{Lemma}\label{lemrho}
 Assume the mapping $x\rightarrow F(x,t)$ to be
invertible, then for $v\in H^1_0(\Omega)$ we have
\begin{itemize}
\item For $n=1$,
\begin{equation}\label{ghhjsdaswwwfh52}
\big(a[v-\Rht v]\big)^\frac{1}{2} \leq   C_{X} h
 \|v\circ F^{-1}(.,t)\|_{W^{2,2}_D} \|a\nabla
 F\|_{L^\infty(\Omega_T)}^\frac{1}{2}.
\end{equation}
\item For $n\geq 2$,
\begin{equation}\label{ghhjsdaswwszwfh52}
\begin{split}
\big(a[v-\Rht v]\big)^\frac{1}{2}  \leq  & C_{X} h
 \|v\circ F^{-1}(.,t)\|_{W^{2,2}_D} \\& \times C_n \mu_\sigma^{\frac{n-1}{4}}
\big\|(\Tr[\sigma])^{-1}\big\|_{L^\infty(\Omega_T)}^{\frac{n-2}{4}}.
\end{split}
\end{equation}
\end{itemize}
\end{Lemma}
\begin{Remark}
Recall that $\mu_{\sigma}$ is given by equation \eref{sjjdhd27}and
it is easy to check that $\mu_{\sigma}$ is bounded by an increasing
function of $(1-\beta_{\sigma})^{-1}$.
\end{Remark}
\begin{proof}
Using  the change of coordinates $y=F(x,t)$ we obtain that (we write
$\hat{v}:=v \circ F^{-1}$)
\begin{equation}\label{ghhjasabbhfh52}
a[v]=Q[\hat{v}]
\end{equation}
with
\begin{equation}\label{ghhssabbhfh52}
\Q[w]:=\int_\Omega {^t\nabla w(y,t)}Q(y,t)\nabla w(y,t)\,dy
\end{equation}
and
\begin{equation}\label{ghhsssh52}
Q(y,t):=\frac{\sigma}{\det(\nabla F)}\circ F^{-1}.
\end{equation}
Using the definition of $\Rht v$ we obtain that
\begin{equation}\label{ghhjasszbbhfh52}
\Q[\hat{v}-\widehat{\Rht v}]=\inf_{\varphi \in X_h
}Q[\hat{v}-\varphi].
\end{equation}
Using property \eref{approp} we obtain that
\begin{equation}\label{ghhjawjbbhfh52}
\Q[\hat{v}-\widehat{\Rht v}]\leq \lambda_{\max}(Q) C_{X}^2 h^2
\|\hat{v}\|_{W^{2,2}_D (T)}^2.
\end{equation}
It is easy to obtain that
\begin{itemize}
\item $n=1$.
\begin{equation}
\lambda_{\max}(Q)\leq \|a\nabla F\|_{L^\infty(\Omega_T)}.
\end{equation}
\item $n\geq 2$.
\begin{equation}
\lambda_{\max}(Q)\leq C_n \mu_\sigma^{\frac{n-1}{2}}
\big\|(\Tr[\sigma])^{-1}\big\|_{L^\infty(\Omega_T)}^{\frac{n}{2}-1}.
\end{equation}
\end{itemize}
\end{proof}

\begin{Lemma}\label{slskjdddd823}
Assume that $\partial_t a \equiv 0$, $\Omega$ is convex,
$\beta_{\sigma}<1$ and $(\Tr[\sigma])^{-1}\in L^\infty(\Omega_T)$
then
\begin{equation}\label{sjswhdawdsdlswe}
\A_T[\rho] \leq C h^2 \|g\|_{L^2(\Omega_T)}^2.
\end{equation}
\end{Lemma}
\begin{Remark}
The constant $C$ depends on $C_X$, $n$, $\Omega$,
$\lambda_{\min}(a)$ and\\
$\big\|(\Tr[\sigma])^{-1}\big\|_{L^\infty(\Omega_T)}$. If $n\geq 5$
it also depends on $\big\|\Tr[\sigma]\big\|_{L^\infty(\Omega_T)}$
and if $n=1$ it also depends on $\lambda_{\max}(a)$.
\end{Remark}
\begin{proof}
The proof is a straightforward application of lemma \ref{lemrho} and
theorem \ref{ksjshjsxxsd8721}. Observe that in dimension one
$a\nabla F=(\int_\Omega a^{-1})^{-1}$
\end{proof}

\begin{Lemma}\label{slskj823}
Assume that $\partial_t a \equiv 0$, $\Omega$ is convex,
$\beta_{\sigma}<1$ and $(\Tr[\sigma])^{-1}\in L^\infty(\Omega_T)$
then
\begin{equation}\label{sjswhawdsdlswe}
\|\rho\|_{L^2(\Omega_T)} \leq C h^2 \|g\|_{L^2(\Omega_T)}.
\end{equation}
\end{Lemma}
\begin{Remark}
The constant $C$ depends on $C_X$, $n$, $\Omega$,
$\lambda_{\min}(a)$ and\\
$\big\|(\Tr[\sigma])^{-1}\big\|_{L^\infty(\Omega_T)}$. If $n\geq 5$
it also depends on $\big\|\Tr[\sigma]\big\|_{L^\infty(\Omega_T)}$
and if $n=1$ it also depends on $\lambda_{\max}(a)$.
\end{Remark}

\begin{proof}
The proof follows from standard duality techniques (see for instance
theorem 5.7.6 of \cite{BreSco02}). We choose $v\in
L^2(0,T,H^1_0(\Omega))$ to be the solution of the following linear
problem: for all $w\in L^2(0,T,H^1_0(\Omega))$
\begin{equation}\label{dhjdssh61}
A_T[w,v]=(w,\rho)_{L^2(\Omega_T)}.
\end{equation}
Choosing $w=\rho$ in equation \eref{dhjdssh61}  we deduce that
\begin{equation}
\|\rho\|_{L^2(\Omega_T))}^2 =\A_T[\rho,v-\Rht v].
\end{equation}
Using Cauchy Schwartz inequality we deduce that
\begin{equation}\label{sjhdlssswe}
\|\rho\|_{L^2(\Omega_T)}^2 \leq \big(\A_T[\rho]\big)^\frac{1}{2}
\big(\A_T[v-\Rht v]\big)^\frac{1}{2}.
\end{equation}
Using theorem \ref{ksjshjsxxsd8721} we obtain that
\begin{equation}\label{ksjsjh61}
\|\hat{v}\|_{L^2(0,T,W^{2,2}_D(\Omega))}\leq C
\|\rho\|_{L^2(\Omega_T)}.
\end{equation}
Using lemma \ref{lemrho} we obtain that
\begin{equation}\label{sjhddsdlswe}
\big(\A_T[v-\Rht v]\big)^\frac{1}{2} \leq C h
\|\rho\|_{L^2(\Omega_T)}.
\end{equation}
It follows that
\begin{equation}\label{sjswhawdsdlsmxwe}
\|\rho\|_{L^2(\Omega_T)} \leq C h \big(\A_T[\rho]\big)^\frac{1}{2}.
\end{equation}
 We deduce the lemma by applying  lemma \ref{slskjdddd823}  to
bound $A_T[\rho]$.
\end{proof}

\begin{Theorem}\label{slsakdeswj823}
Assume that $\partial_t a \equiv 0$, $\Omega$ is convex,
$\beta_{\sigma}<1$ and $(\Tr[\sigma])^{-1}\in L^\infty(\Omega_T)$
then
\begin{equation}\label{sjswhsawdswsdlswe}
\big\|(u-u_h)(.,T)\big\|_{L^2(\Omega)}^2+\big\|u-u_h\big\|_{L^2(0,T;H^1_0(\Omega))}^2
\leq C h^2 \|g\|_{L^2(\Omega_T)}^2.
\end{equation}
\end{Theorem}
\begin{Remark}
The constant $C$ depends on $C_X$, $n$, $\Omega$,
$\lambda_{\min}(a)$ and\\
$\big\|(\Tr[\sigma])^{-1}\big\|_{L^\infty(\Omega_T)}$. If $n\geq 5$
it also depends on $\big\|\Tr[\sigma]\big\|_{L^\infty(\Omega_T)}$
and if $n=1$ it also depends on $\lambda_{\max}(a)$.
\end{Remark}
\begin{proof}
The proof is a straightforward application of lemmas \ref{slskj823},
\ref{slskjdddd823},  \ref{jshsszskssjsbh61RED}, \ref{jkshssssde7a32}
and \ref{jshsszskjsh61RED}.
\end{proof}

\subsubsection{Medium with a continuum of time scales.}
In this subsection we will assume that the finite elements are in
$H^2(\Omega)\cap H^1_0(\Omega)$ and satisfy inverse inequality
\eref{appsrop3}.
\begin{Lemma}\label{jkshssssdess7a32}
\begin{equation}\label{ghhjaawjssbbhfhs52az}
\begin{split}
\frac{1}{2}\big\|(u-u_h)(t)\big\|^2_{L^2(\Omega)}+\A_t[u-u_h]=&\int_{\Omega_t}\frac{\hat{\rho}}{|\det
\nabla F|\circ F^{-1}}\\& \big(\hat{g}+\sum_{i,j=1}^n
\sigma_{i,j}\circ F^{-1}\partial_i \partial_j \hat{u}_h-\partial_t
\hat{u}_h\big).
\end{split}
\end{equation}
\end{Lemma}
\begin{proof}
Consider equation \eref{ghhjaawjszsbbhfh52az}. We have
\begin{equation}
\int_{\Omega_t}\rho\partial_t(u-u_h)=\int_{\Omega_t}\frac{\hat{\rho}}{|\det
\nabla F|\circ F^{-1}}\partial_t (\hat{u}-\hat{u}_h)+
\int_{\Omega_t}\rho \partial_t F (\nabla F)^{-1}\nabla (u-u_h).
\end{equation}
Using equation \eref{ghjagsash52} we obtain that
\begin{equation}
\int_{\Omega_t}\rho \partial_t F (\nabla F)^{-1}\nabla
(u-u_h)=-\A_t[\rho,u-u_h]-\sum_{i,j=1}^n \int_{\Omega_t}\hat{\rho}
Q_{i,j}\partial_i \partial_j (\hat{u}-\hat{u}_h).
\end{equation}
\end{proof}

\begin{Lemma}\label{jskjshsj73}
\begin{equation}
 \big\|\frac{\partial_t \hat{u}_h}{|\det(\nabla
F)|^\frac{1}{2}\circ F^{-1} }\big\|_{L^2(\Omega_T)}\leq 2
\|g\|_{L^2(\Omega_T)}+ C \|\hat{u}_h\|_{L^2(0,T,W^{2,2}_D(\Omega))}.
\end{equation}
where the constant $C$ depends on $n$, $\lambda_{\max}(a)$,
$\|\Tr[\sigma]\|_{L^\infty(\Omega_T)}$, $\mu_\sigma$.
\end{Lemma}
\begin{proof}
Using the change of variable $y=F(x,t)$ in \eref{ghjbbfh52} we
obtain that for all $\varphi\in X_h$
\begin{equation}\label{ghjbbfh533e2}
\begin{cases}
(\varphi,\frac{\partial_t \hat{u}_h}{|\det(\nabla F)|\circ F^{-1} }
)_{L^2(\Omega)}=&\sum_{i,j=1}^n \int_{\Omega}
(\varphi,Q_{i,j}\partial_i\partial_j
\hat{u}_h)_{L^2(\Omega)}\\&+(\varphi,\frac{\hat{g}}{|\det(\nabla
F)|\circ F^{-1} })_{L^2(\Omega)}\\
 \hat{u}_h(x,0)=0.
\end{cases}
\end{equation}
Recall that $Q$ is given by \eref{ghhsssh52}. We choose
$\varphi=\partial_t \hat{u}$ and observe that
\begin{equation}
\frac{\sigma}{|\det \nabla F|^\frac{1}{2}}=\frac{\sigma}{|\det
\sigma|^\frac{1}{4}}|\det a|^\frac{1}{4}.
\end{equation}
Thus
\begin{equation}
 \big\|\frac{\sigma}{|\det \nabla
F|^\frac{1}{2}}\big\|\leq
C\big(n,\lambda_{\max}(a),\|\Tr[\sigma]\|_{L^\infty(\Omega_T)},\mu_\sigma\big).
\end{equation}
We deduce the lemma by Minkowski inequality.
\end{proof}
Combining lemma \ref{jkshssssdess7a32} and lemma \ref{jskjshsj73} we
obtain the following lemma

\begin{Lemma}\label{jkshsssssdess7a32}
\begin{equation}\label{ghhjaawssjssbcbhfhs52az}
\begin{split}
\frac{1}{2}\big\|(u-u_h)(T)\big\|^2_{L^2(\Omega)}+\A_T[u-u_h]\leq
&\|\rho\|_{L^2(\Omega_T)}\Big(\|g\|_{L^2(\Omega_T)}\\&+ C
\|\hat{u}_h\|_{L^2(0,T,W^{2,2}_D(\Omega))} \Big).
\end{split}
\end{equation}
where the constant $C$ depends on $n$, $\lambda_{\max}(a)$,
$\|\Tr[\sigma]\|_{L^\infty(\Omega_T)}$, $\mu_\sigma$.
\end{Lemma}

\begin{Lemma}\label{slskjsseee823l}
Assume that $\Omega$ is convex, and condition \ref{slkssjk88271} is
satisfied then
\begin{equation}\label{sjswhfawssdsdlswe}
\|\rho\|_{L^2(\Omega_T)} \leq C h^2 \|g\|_{L^2(\Omega_T)}.
\end{equation}
\end{Lemma}
\begin{Remark}
The constant $C$ depends on $C_X$, $n$, $\Omega$, $\delta$ and
$\epsilon$, $\lambda_{\min}(a)$ and $\lambda_{\max}(a)$.
\end{Remark}
\begin{proof}
The proof is similar to the proof of lemma \ref{slskj823}. As in
\ref{slskj823} we choose $v\in L^2(0,T,H^1_0(\Omega))$ to be the
solution of the following linear problem: for all $w\in
L^2(0,T,H^1_0(\Omega))$
\begin{equation}\label{dhjdssh61m}
A_T[w,v]=(w,\rho)_{L^2(\Omega_T)}.
\end{equation}
Choosing $w=\rho$ in equation \eref{dhjdssh61m}  we deduce that
\begin{equation}
\|\rho\|_{L^2(\Omega_T))}^2 =\A_T[\rho,v-\Rht v].
\end{equation}
Using Cauchy Schwartz inequality we deduce that
\begin{equation}\label{sjhdlssswem}
\|\rho\|_{L^2(\Omega_T)}^2 \leq \big(\A_T[\rho]\big)^\frac{1}{2}
\big(\A_T[v-\Rht v]\big)^\frac{1}{2}.
\end{equation}
Using theorem \ref{skjhsj823} we obtain that
\begin{equation}\label{ksjsjh61m}
\|\hat{v}\|_{L^2(0,T,W^{2,2}_D(\Omega))}\leq C
\|\rho\|_{L^2(\Omega_T)}.
\end{equation}
Using lemma \ref{lemrho} we obtain that
\begin{equation}\label{sjhddsdlswem}
\big(\A_T[v-\Rht v]\big)^\frac{1}{2} \leq C h
\|\rho\|_{L^2(\Omega_T)}.
\end{equation}
It follows that
\begin{equation}\label{sjswhawdsdlsmxwem}
\|\rho\|_{L^2(\Omega)} \leq C h \big(\A_T[\rho]\big)^\frac{1}{2}.
\end{equation}
 We deduce the lemma by applying  lemma \ref{lemrho} and theorem \ref{skjhsj823}  to
bound $A_T[\rho]$.
\end{proof}

\begin{Lemma}\label{slskjsseee823}
Assume that $\Omega$ is convex, and that $\Tr[\sigma]\in
L^\infty(\Omega_T)$.
\begin{equation}\label{sjswhawssdsdlswe}
\|\hat{u}_h\|_{L^2(0,T,W^{2,2}_D(\Omega))}\leq \frac{C}{h}
\|g\|_{L^2(\Omega_T)}.
\end{equation}
\end{Lemma}
\begin{Remark}
The constant $C$ depends on $C_X$, $n$, $\Omega$,
$\lambda_{\min}(a)$ and $\lambda_{\max}(a)$ and
$\|\Tr[\sigma]\|_{L^\infty(\Omega_T)}$.
\end{Remark}
\begin{proof}
Using the inverse inequality \eref{appsrop3} of the finite elements
we obtain that
\begin{equation}\label{sjswhawssdsdflswe}
\|\hat{u}_h\|_{L^2(0,T,W^{2,2}_D(\Omega))}\leq \frac{C_X}{h}
\|\nabla \hat{u}_h\|_{L^2(0,T,W^{2,2}_D(\Omega))}.
\end{equation}
Using the change of variables $y=F(x)$ we obtain that
\begin{equation}\label{sjswhawssdsdfflswe}
 \|\nabla
\hat{u}_h\|_{L^2(0,T,W^{2,2}_D(\Omega))}^2 \leq C \A_T[u_h]
\end{equation}
where $C$ depends on $n$, $\lambda_{\min}(a)$ and
$\|\Tr[\sigma]\|_{L^\infty(\Omega_T)}$. We deduce the lemma by using
lemma \ref{sjhsgh7}.
\end{proof}

\begin{Theorem}
Assume that $\Omega$ is convex, and condition \ref{slkssjk88271} is
satisfied then
\begin{equation}\label{ghhjaawssjssbbhfhsnb52az}
\begin{split}
\frac{1}{2}\big\|(u-u_h)(T)\big\|^2_{L^2(\Omega)}+\A_T[u-u_h]\leq C
h \|g\|_{L^2(\Omega_T)}^2.
\end{split}
\end{equation}
\end{Theorem}
\begin{Remark}
The constant $C$ depends on $C_X$, $n$, $\Omega$, $\delta$ and
$\epsilon$, $\lambda_{\min}(a)$ and $\lambda_{\max}(a)$.
\end{Remark}
\begin{proof}
The proof is a straightforward application of lemma
\ref{jkshsssssdess7a32}, lemma \ref{slskjsseee823l} and lemma
\ref{slskjsseee823}.
\end{proof}

\subsubsection{Homogenization and discretization in time.}
We use the notation of subsection \ref{jksjhs89b}. First let us
observe that the numerical scheme associated to
\eref{ghjddwsdcszdbbsfh52} is stable. Indeed choosing $\psi=v_{n+1}$
one gets
\begin{equation}\label{ghjddsdsswsdcszdsbbsfh52}
\begin{split}
 \big|v_{n+1}(t_{n+1})\big|^2_{L^2(\Omega)}=&\big(v_{n+1}(t_{n}),
v_{n}(t_{n})\big)_{L^2(\Omega)}\\+\frac{1}{2}\big( \big|
v_{n+1}(t_{n+1})&\big|^2_{L^2(\Omega)}- \big|
v_{n+1}(t_{n})\big|^2_{L^2(\Omega)}\big)
\\-\int_{t_n}^{t_{n+1}}\Big(a\big[v_{n+1}(t)]&+\big(v_{n+1}(t),g(t)\big)_{L^2(\Omega)}\Big)\,dt.
\end{split}
\end{equation}
It follows by  Cauchy-Schwartz inequality that
\begin{equation}\label{ghjddsdscswsdcszddsbbsfh52}
\begin{split}
 \frac{1}{2}\big|v_{n+1}(t_{n+1})\big|^2_{L^2(\Omega)}\leq & \frac{1}{2}\big|v_{n}(t_{n})\big|^2_{L^2(\Omega)}
\\-\int_{t_n}^{t_{n+1}}\Big(a\big[v_{n+1}(t)]&+\big(v_{n+1}(t),g(t)\big)_{L^2(\Omega)}\Big)\,dt.
\end{split}
\end{equation}
Hence  using Poincar\'{e} and Minkowski inequalities one obtains
that
\begin{equation}\label{ghjddsdsswsdcszddsbbsfh52}
\begin{split}
  \big|v_{n+1}(t_{n+1})\big|^2_{L^2(\Omega)}+\int_{t_n}^{t_{n+1}}a\big[v_{n+1}(t)]\,dt  \leq& \big|v_{n}(t_{n})\big|^2_{L^2(\Omega)}
\\+\frac{C_{n,\Omega}}{\lambda_{\min}(a)}&\int_{t_n}^{t_{n+1}}\big|g(t)\big|^2_{L^2(\Omega)}\,dt.
\end{split}
\end{equation}
which implies theorem \ref{sshjhsjjs823} and the stability of the
scheme. Integrating \eref{ghjbbfh52} with respect to time we obtain
that for $\psi\in V$,
\begin{equation}\label{ghjddwszdbbfh52}
\begin{split}
 \big(\psi(t_{n+1}),
u_h(t_{n+1})\big)_{L^2(\Omega)}=&\big(\psi(t_{n}),
u_h(t_{n})\big)_{L^2(\Omega)}+\int_{t_n}^{t_{n+1}}\Big(\big(\partial_t\psi(t),
u_h(t)\big)_{L^2(\Omega)}\\&-a\big[\psi(t),u_h(t)]+\big(\psi(t),g(t)\big)_{L^2(\Omega)}
\Big)\,dt.
\end{split}
\end{equation}
Let us write $(u^i)$ the coordinates of $u_h$ associated to the
basis $(\varphi_i \circ F)$, i.e.
\begin{equation}
u_h(x,t):=\sum_{i} u^i(t) \varphi_i(F(x,t)).
\end{equation}
Let us define
\begin{equation}
u_n(x,t):=\sum_{i} u^i(t_n) \varphi_i(F(x,t)).
\end{equation}
Subtracting \eref{ghjddwszdbbfh52} and \eref{ghjddwsdcszdbbsfh52} we
obtain that for $\psi \in Z_T$,
\begin{equation}\label{ghjdddxxwszdssdbbfh52}
\begin{split}
 \big(\psi(t_{n+1}),
(u_{n+1}-v_{n+1})(t_{n+1})\big)_{L^2(\Omega)}=&\big(\psi(t_{n}),
(u_n-v_{n})(t_{n})\big)_{L^2(\Omega)}\\&+\int_{t_n}^{t_{n+1}}\Big(\big(\partial_t\psi(t),
(u_h-v_{n+1})(t)\big)_{L^2(\Omega)}\\&-a\big[\psi(t),(u_h-v_{n+1})(t)]
\Big)\,dt.
\end{split}
\end{equation}
Choosing $\psi=u_{n+1}-v_{n+1}$ we deduce using Cauchy-Schwartz
inequality that
\begin{equation}\label{ghjdzddxxwszsssssdbbfh52}
\begin{split}
\frac{1}{2} & \big| (u_{n+1}-v_{n+1})(t_{n+1})\big|^2_{L^2(\Omega)}+
\int_{t_n}^{t_{n+1}}a\big[(u_{n+1}-v_{n+1})(t)] \,dt \leq\\&
\frac{1}{2}\big|
(u_n-v_{n})(t_{n})\big|^2_{L^2(\Omega)}+\int_{t_n}^{t_{n+1}}\Big(\big(\partial_t(u_{n+1}-v_{n+1})(t),
(u_h-u_{n+1})(t)\big)_{L^2(\Omega)}\\&-a\big[(u_{n+1}-v_{n+1})(t),(u_{h}-u_{n+1})(t)]
\Big)\,dt.
\end{split}
\end{equation}

\paragraph{Time independent medium.}
Observe that if the medium is time independent then
\eref{ghjdzddxxwszsssssdbbfh52} can be written
\begin{equation}\label{gazoszsssssdbbfh52}
\begin{split}
\frac{1}{2} & \big| (u_{n+1}-v_{n+1})(t_{n+1})\big|^2_{L^2(\Omega)}+
\int_{t_n}^{t_{n+1}}a\big[(u_{n+1}-v_{n+1})(t)] \,dt \leq\\&
\frac{1}{2}\big|
(u_n-v_{n})(t_{n})\big|^2_{L^2(\Omega)}-\int_{t_n}^{t_{n+1}}
a\big[(u_{n+1}-v_{n+1})(t),(u_{h}-u_{n+1})(t)] \,dt
\end{split}
\end{equation}
which leads us to
\begin{equation}\label{gazoszwsbbfh52}
\begin{split}
\frac{1}{2} & \big| (u_{n+1}-v_{n+1})(t_{n+1})\big|^2_{L^2(\Omega)}+
\int_{t_n}^{t_{n+1}}a\big[(u_{n+1}-v_{n+1})(t)] \,dt \leq\\&
\frac{1}{2}\big|
(u_n-v_{n})(t_{n})\big|^2_{L^2(\Omega)}\\&+\int_{t_n}^{t_{n+1}}
\int_{t_n}^{t_{n+1}} 1(t<s)\,a\big[(u_{n+1}-v_{n+1})(t),\partial_s
u_h(s)] \,ds\,dt.
\end{split}
\end{equation}
Write $\Delta t:=t_{n+1}-t_n$. Using Minkowski inequality we obtain
that
\begin{equation}\label{gazoszwsbbfh52ppp}
\begin{split}
a\big[(u_{n+1}-v_{n+1})(t),\partial_s u_h(s)]\leq &\frac{1}{2 \Delta
t}a\big[(u_{n+1}-v_{n+1})(t)]\\&+ \frac{1}{2}\Delta t
a\big[\partial_s u_h(s)].
\end{split}
\end{equation}

It follows  from \eref{gazoszwsbbfh52} that
\begin{equation}\label{gazoszwsbbfh52c}
\begin{split}
 & \big| (u_{n+1}-v_{n+1})(t_{n+1})\big|^2_{L^2(\Omega)}+
\int_{t_n}^{t_{n+1}}a\big[(u_{n+1}-v_{n+1})(t)] \,dt \leq\\& \big|
(u_n-v_{n})(t_{n})\big|^2_{L^2(\Omega)}+ |\Delta t|^2
\int_{t_n}^{t_{n+1}} \,a\big[\partial_s u_h(s)] \,ds.
\end{split}
\end{equation}
Observing that
\begin{equation}\label{gazoszwsbbfh52cc}
\begin{split}
\int_{t_n}^{t_{n+1}}a\big[(u_{n+1}-v_{n+1})(t)] \,dt\geq& 0.5
\int_{t_n}^{t_{n+1}}a\big[(u_{h}-v_{n+1})(t)]
\,dt\\&-\int_{t_n}^{t_{n+1}}a\big[(u_{h}-u_{n+1})(t)] \,dt
\end{split}
\end{equation}
and
\begin{equation}\label{gazoszwsbbfh52vcb}
\begin{split}
\int_{t_n}^{t_{n+1}}a\big[(u_{h}-u_{n+1})(t)] \,dt \leq |\Delta t|^2
\int_{t_n}^{t_{n+1}} \,a\big[\partial_s u_h(s)] \,ds
\end{split}
\end{equation}
we obtain that
\begin{equation}\label{gazoszwsbbfh52vcdf}
\begin{split}
 & \big| (u_{n+1}-v_{n+1})(t_{n+1})\big|^2_{L^2(\Omega)}+ 0.5
\int_{t_n}^{t_{n+1}}a\big[(u_{h}-v_{n+1})(t)] \,dt \leq\\& \big|
(u_n-v_{n})(t_{n})\big|^2_{L^2(\Omega)}+\frac{3}{2} |\Delta t|^2
\int_{t_n}^{t_{n+1}} \,a\big[\partial_s u_h(s)] \,ds.
\end{split}
\end{equation}
In conclusion we have obtained the following following lemma
\begin{Lemma}\label{jskjdhdjdh723}
Let $v\in Z_T$ be the solution of \eref{ghjddwsdcszdbbsfh52}. We
have
\begin{equation}\label{gazoszwsbbfhgf52}
\begin{split}
 \big\| (u_h-v)(T)\big\|^2_{L^2(\Omega)}+
\int_{0}^{T}a\big[(u_h-v)(t)] \,dt \leq 3 |\Delta t|^2 \int_{0}^{T}
\,a\big[\partial_s u_h(s)] \,ds.
\end{split}
\end{equation}
\end{Lemma}
Combining lemma \ref{jshsszskjsh61} with lemma \ref{jskjdhdjdh723}
we obtain the following theorem:
\begin{Theorem}\label{jskjdhdjdhxzs723}
Let $v\in Z_T$ be the solution of \eref{ghjddwsdcszdbbsfh52}. We
have
\begin{equation}\label{gazoszwsbdsbfh52}
\begin{split}
 \big\| (u_h-v)(T)\big\|^2_{L^2(\Omega)}+
&\int_{0}^{T}a\big[(u_h-v)(t)] \,dt \leq 3 |\Delta t|^2
\\&\Big(\frac{4}{\lambda_{\min}(a)} \|\partial_t
g\|_{L^2(0,T,H^{-1}(\Omega))}^2+\big\|g(.,0)\big\|_{L^2(\Omega)}^2\Big).
\end{split}
\end{equation}
\end{Theorem}

\paragraph{Time dependent medium.}
Observe that $$\partial_t(u_{n+1}-v_{n+1})=\partial_t F (\nabla
F)^{-1}\nabla (u_{n+1}-v_{n+1}).$$ It follows after writing
$\partial_t F =\diiv a \nabla F$, integration by parts and using the
change of variables $y=F(x,t)$ in \eref{ghjdzddxxwszsssssdbbfh52}
that
\begin{equation}\label{ghjdzddxxwszsskksssdddbbfh52}
\begin{split}
\frac{1}{2} & \big| (u_{n+1}-v_{n+1})(t_{n+1})\big|^2_{L^2(\Omega)}+
\int_{t_n}^{t_{n+1}}a\big[(u_{n+1}-v_{n+1})(t)] \,dt \leq\\&
\frac{1}{2}\big|
(u_n-v_{n})(t_{n})\big|^2_{L^2(\Omega)}-2\int_{t_n}^{t_{n+1}}
a\big[(u_{n+1}-v_{n+1})(t),(u_{h}-u_{n+1})(t)\big] \,dt\\&
-\sum_{i,j} \int_{t_n}^{t_{n+1}} \int_{\Omega}
(\hat{u}_{h}-\hat{u}_{n+1})Q_{i,j}\partial_i \partial_j
(\hat{u}_{n+1}-\hat{v}_{n+1}) \,dt\,dy.
\end{split}
\end{equation}
Hence using Minkowski inequality  we obtain that
\begin{equation}\label{ghjdzddxxwszsskksssdddbbfcih52}
\begin{split}
 & \big| (u_{n+1}-v_{n+1})(t_{n+1})\big|^2_{L^2(\Omega)}+
 \int_{t_n}^{t_{n+1}}a\big[(u_{n+1}-v_{n+1})(t)] \,dt \leq\\&
\big| (u_n-v_{n})(t_{n})\big|^2_{L^2(\Omega)}+4 \int_{t_n}^{t_{n+1}}
a\big[(u_{h}-u_{n+1})(t)\big] \,dt\\& -2\sum_{i,j}
\int_{t_n}^{t_{n+1}} \int_{\Omega}
(\hat{u}_{h}-\hat{u}_{n+1})Q_{i,j}\partial_i \partial_j
(\hat{u}_{n+1}-\hat{v}_{n+1}) \,dt\,dy.
\end{split}
\end{equation}
Using Minkowski inequality we obtain that
\begin{equation}\label{ghjdzddxxwszsskksssdddiibbfh52}
\begin{split}
\Big|\sum_{i,j} \int_{t_n}^{t_{n+1}}& \int_{\Omega}
(\hat{u}_{h}-\hat{u}_{n+1})Q_{i,j}\partial_i \partial_j
(\hat{u}_{n+1}-\hat{v}_{n+1}) \,dt\,dy\Big|\leq\\& C_A n^2
\int_{t_n}^{t_{n+1}} \int_{\Omega}
|\hat{u}_{h}-\hat{u}_{n+1}|^2\,dt\,dy
\\&+\frac{\lambda_{\max}(Q)}{C_A} \int_{t_n}^{t_{n+1}} \sum_{i,j}\int_{\Omega}
|\partial_i \partial_j (\hat{u}_{n+1}-\hat{v}_{n+1})|^2
\,dt\,dy.
\end{split}
\end{equation}
Using the inverse inequality \eref{appsrop3} and the change of
variable $y=F(x)$ we obtain that
\begin{equation}\label{ghjdzddxxwsiizsskksssdddbbfh52}
\begin{split}
 \int_{t_n}^{t_{n+1}}
\sum_{i,j}\int_{\Omega} |\partial_i \partial_j
(\hat{u}_{n+1}-\hat{v}_{n+1})&|^2 \,dt\,dy\leq \frac{C_X}{h^2
\lambda_{\min}(Q)}
\\&\int_{t_n}^{t_{n+1}}a\big[(u_{n+1}-v_{n+1})(t)] \,dt
.\end{split}
\end{equation}
In resume, choosing $C_A= \frac{4 C_X \lambda_{\max}(Q)}{h^2
\lambda_{\min}(Q)}$ we have obtained that
\begin{equation}\label{ghjdzddxxwsdbfh52}
\begin{split}
 & \big| (u_{n+1}-v_{n+1})(t_{n+1})\big|^2_{L^2(\Omega)}+0.5
 \int_{t_n}^{t_{n+1}}a\big[(u_{n+1}-v_{n+1})(t)] \,dt \leq\\&
\big| (u_n-v_{n})(t_{n})\big|^2_{L^2(\Omega)}+8 \int_{t_n}^{t_{n+1}}
a\big[(u_{h}-u_{n+1})(t)\big] \,dt\\& + \frac{8 C_X
\lambda_{\max}(Q)}{h^2 \lambda_{\min}(Q)} n^2 \int_{t_n}^{t_{n+1}}
\int_{\Omega} |\hat{u}_{h}-\hat{u}_{n+1}|^2\,dt\,dy .\end{split}
\end{equation}
And a computation similar to the one leading to
\eref{gazoszwsbbfh52ppp} gives us
\begin{equation}\label{ghjdzsxxwsdbfh52o}
\begin{split}
 & \big| (u_{n+1}-v_{n+1})(t_{n+1})\big|^2_{L^2(\Omega)}+\frac{1}{4}
 \int_{t_n}^{t_{n+1}}a\big[(u_h-v_{n+1})(t)] \,dt \leq\\&
\big| (u_n-v_{n})(t_{n})\big|^2_{L^2(\Omega)}+9 \int_{t_n}^{t_{n+1}}
a\big[(u_{h}-u_{n+1})(t)\big] \,dt\\& + \frac{8 C_X
\lambda_{\max}(Q)}{h^2 \lambda_{\min}(Q)} n^2 \int_{t_n}^{t_{n+1}}
\int_{\Omega} |\hat{u}_{h}-\hat{u}_{n+1}|^2\,dt\,dy .\end{split}
\end{equation}
Moreover using the change of variables $F(x)=y$ and the inverse
inequality \eref{appssswerop3} we obtain that
\begin{equation}\label{ghjdzsxxwsdbfh52}
\begin{split}
 \int_{t_n}^{t_{n+1}}
a\big[(u_{h}-u_{n+1})(t)\big] \,dt \leq \frac{C_X
\lambda_{\max}(Q)}{h^2} \int_{t_n}^{t_{n+1}} \int_{\Omega}
|\hat{u}_{h}-\hat{u}_{n+1}|^2\,dt\,dy.
\end{split}
\end{equation}
Let us also observe that
\begin{equation}\label{ghjdzsswaqwsdbfh52}
\begin{split}
 \int_{t_n}^{t_{n+1}} \int_{\Omega}
|\hat{u}_{h}-\hat{u}_{n+1}|^2\,dt\,dy\leq |\Delta t|^2
\int_{t_n}^{t_{n+1}} \int_{\Omega} |\partial_t
\hat{u}_{h}|^2\,dt\,dy.
\end{split}
\end{equation}
It follows that
\begin{equation}\label{ghjdzsxxwsdbfh52oo}
\begin{split}
 & \big| (u_{n+1}-v_{n+1})(t_{n+1})\big|^2_{L^2(\Omega)}+\frac{1}{4}
 \int_{t_n}^{t_{n+1}}a\big[(u_h-v_{n+1})(t)] \,dt \leq\\&
\big| (u_n-v_{n})(t_{n})\big|^2_{L^2(\Omega)}+ C_B \frac{|\Delta
t|^2}{h^2} \int_{t_n}^{t_{n+1}} \int_{\Omega} |\partial_t
\hat{u}_{h}|^2\,dt\,dy.
\end{split}
\end{equation}
with
\begin{equation}\label{ghssawsdbfh52}
\begin{split}
C_B=C_n C_X \lambda_{\max}(Q) (1+\frac{1}{\lambda_{\min}(Q)})
.\end{split}
\end{equation}
We deduce that
\begin{equation}\label{ghjdzsxxwsdssabfh52}
\begin{split}
 \big\| (u_h-v)(T)\big\|^2_{L^2(\Omega)}+\frac{1}{4}
 \int_{0}^{T}a\big[(u_h-v)(t)] \,dt \leq C_B \frac{|\Delta t|^2}{h^2} \int_{0}^{T} \int_{\Omega}
|\partial_t \hat{u}_{h}|^2\,dt\,dy.
\end{split}
\end{equation}
Using lemma \ref{dkdjjddh} to control $C_B$ and combining
\eref{ghjdzsxxwsdssabfh52} with theorem \ref{skjhsj823} we obtain
the following theorem.
\begin{Theorem}\label{skjhseeeddfddsf23}
Assume that $\Omega$ is convex, and condition \ref{slkssjk88271} is
satisfied. Let $v\in Z_T$ be the solution of
\eref{ghjddwsdcszdbbsfh52}, we have
\begin{equation}\label{ghjdzssasdssabfh52}
\begin{split}
 \big\| (u_h-v)(T)\big\|^2_{L^2(\Omega)}+\frac{1}{4}
 \int_{0}^{T}a\big[(u_h-v)(t)] \,dt \leq C \frac{|\Delta t|^2}{h^2}
 \|g\|_{L^2(\Omega_T)}^2
.\end{split}
\end{equation}
where $C$ depends on $\Omega$, $n$, $\delta$, $\epsilon$,
$\lambda_{\min}(a)$ and $\lambda_{\max}(a)$.
\end{Theorem}

\clearpage

\section{Numerical Experiments}\label{ksjsskhs821}
The purpose of this section is to give several illustrations of the
implementation of this method. The domain is the unit square in
dimension two. Equation \eref{ghjh52} is solved on a fine
tessellation characterized by $16129$ interior nodes (degree of
freedoms).

Three different coarse tessellations are considered, one with $9$
 degrees of freedoms (noted \emph{dof} in the tables), one with $49$ and the last one with $225$.

The parabolic operator associated to equation \eref{ghjh52} has been
homogenized onto these coarse meshes using the method the method
presented in this paper. We have chosen splines to span the space
$X_h$ introduced in subsection \ref{sub2}.

\subsection{Time independent examples.}
\begin{figure}[httb]
  \begin{center}
    \subfigure[$u$.]
    {\includegraphics[width=0.35\textwidth,height= 0.3\textwidth]{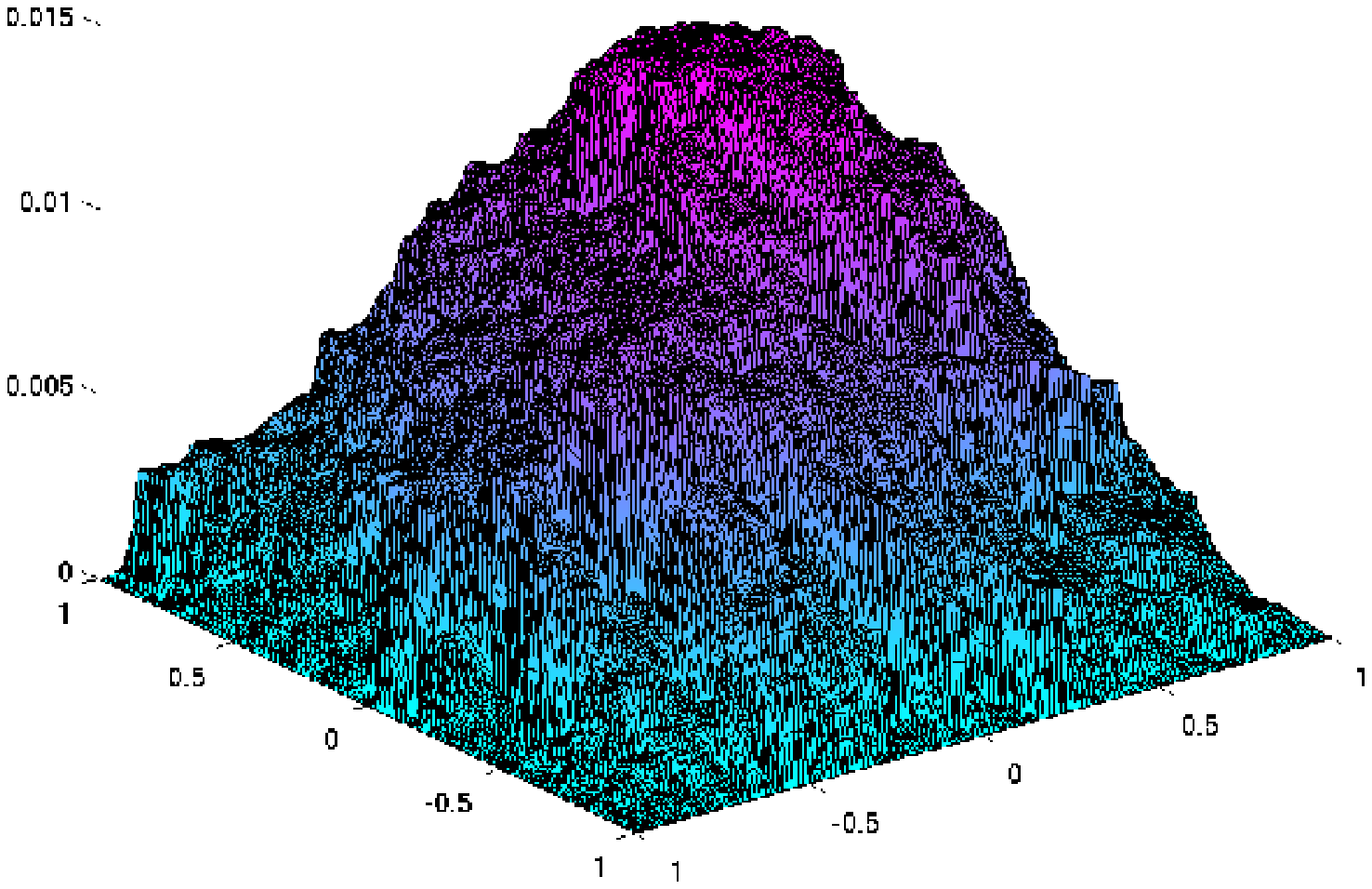}}
    \goodgap
    \subfigure[$u_h$.]
    {\includegraphics[width=0.35\textwidth,height= 0.3\textwidth]{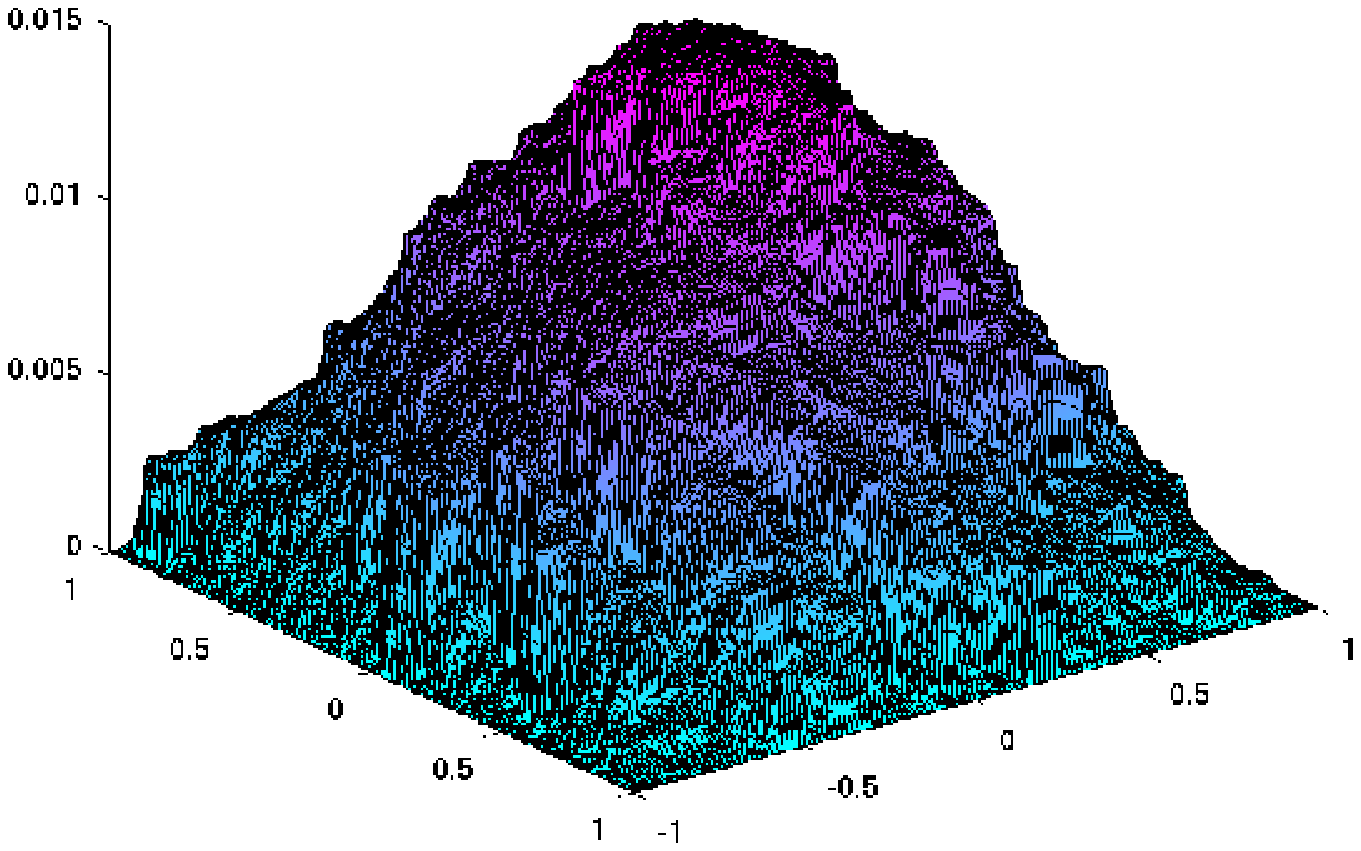}}\\
    \caption{$u$ computed on $16641$ interior nodes and $u_h$ computed on $9$ interior nodes.}
    \label{eruuhzaa4}
\end{center}
\end{figure}

\begin{example}
\label{exa:siteprco}Time independent site percolation.
\end{example}
In this example we consider the site percolating medium associated
 to figure \ref{ap7}. The fine mesh is characterized by
$16641$ nodes. \eref{ghjh52} has been homogenized to three different
coarse meshes with $9$, $49$ and $225$ interior nodes using the
method described here and splines for the space $X_h$. \eref{ghjh52}
has been solved with the fine mesh operator and the coarse mesh
operators with $g=1$ and $g=\sin(2.4 x-1.8 y+2 \pi t)$. The fine
mesh and coarse mesh errors are given in tables \ref{cerrtindepp7},
\ref{ferrtindepp7}, \ref{cerrtindeptdrhsp7},
\ref{ferrtindeptdrhsp7}. Figure \ref{eruuhzaa4} shows $u$ computed
on $16641$ interior nodes and $u_h$ computed on $9$ interior nodes
in the case $g=1$ at time $1$.

\begin{table}
\begin{center}
\caption{Coarse Mesh Error. Time Independent Site Percolation with
$g=1$.} \label{cerrtindepp7}
\begin{tabular}{|c|c|c|c|c|}
\hline dof& $L^{1}$& $L^{\infty}$& $L^{2}$& $H^{1}$\tabularnewline
\hline 9& 0.0142& 0.0389& 0.0168& 0.0366\tabularnewline \hline 49&
0.0077& 0.0450& 0.0101& 0.0482\tabularnewline \hline 225& 0.0035&
0.0228& 0.0060& 0.0293\tabularnewline \hline
\end{tabular}
\end{center}
\end{table}

\begin{table}
\begin{center}
\caption{Fine Mesh Error. Time Independent Site Percolation with
 $g=1$.} \label{ferrtindepp7}
\begin{tabular}{|c|c|c|c|c|}
\hline dof& $L^{1}$& $L^{\infty}$& $L^{2}$& $H^{1}$\tabularnewline
\hline 9& 0.0196& 0.0843& 0.0251& 0.1193\tabularnewline \hline 49&
0.0136& 0.0698& 0.0184& 0.1028\tabularnewline \hline 225& 0.0040&
0.0243& 0.0070& 0.0485\tabularnewline \hline
\end{tabular}
\end{center}
\end{table}

\begin{table}
\begin{center}
\caption{Coarse Mesh Error. Time Independent Percolation Case with
$g=\sin(2.4 x-1.8 y+2 \pi t)$.} \label{cerrtindeptdrhsp7}
\begin{tabular}{|c|c|c|c|c|}
\hline dof& $L^{1}$& $L^{\infty}$& $L^{2}$& $H^{1}$\tabularnewline
\hline 9& 0.0236& 0.0569& 0.0262& 0.0477\tabularnewline \hline 49&
0.0181& 0.0571& 0.0215& 0.0558\tabularnewline \hline 225& 0.0119&
0.0774& 0.0167& 0.0939\tabularnewline \hline
\end{tabular}
\end{center}
\end{table}

\begin{table}
\begin{center}
\caption{Fine Mesh Error. Time Independent Percolation  with
$g=\sin(2.4 x-1.8 y+2 \pi t)$.} \label{ferrtindeptdrhsp7}
\begin{tabular}{|c|c|c|c|c|}
\hline dof& $L^{1}$& $L^{\infty}$& $L^{2}$& $H^{1}$\tabularnewline
\hline 9& 0.0424& 0.1099& 0.0512& 0.1712\tabularnewline \hline 49&
0.0277& 0.0985& 0.0348& 0.1451\tabularnewline \hline 225& 0.0174&
0.0886& 0.0242& 0.1192\tabularnewline \hline
\end{tabular}
\end{center}
\end{table}

\begin{figure}[httb]
\begin{center}
\includegraphics[%
  scale=0.3]{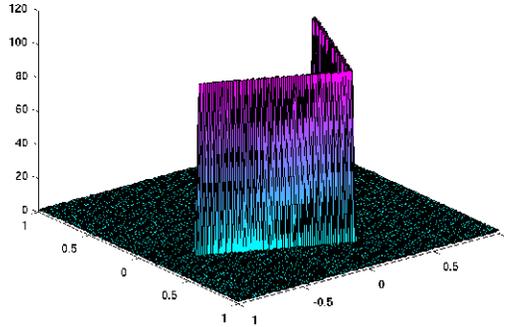}
\caption{High Conductivity Channel superposed on a random medium.}
\label{ap6}
\end{center}
\end{figure}

\clearpage
\begin{example}
\label{exa:chan}Time independent high conductivity channel.
\end{example}
 In this example $a$ is random and characterized by a fine and long
ranged high conductivity channel. We choose $a(x)=100$, if $x$ is in
the channel, and $a(x)=O(1)$ and random, if $x$ is not in the
channel. The media is illustrated in figure \ref{ap6}

Tables \ref{cerrp6tindepsp} and \ref{ferrp6tindepsp} give the coarse
and fine meshes errors.

\begin{table}
\begin{center}
\caption{Coarse Mesh Error, high conductivity channel.}
\label{cerrp6tindepsp}
\begin{tabular}{|c|c|c|c|c|}
\hline dof& $L^{1}$& $L^{\infty}$& $L^{2}$& $H^{1}$\tabularnewline
\hline 9& 0.0159& 0.0496& 0.0207& 0.0477\tabularnewline \hline 49&
0.0067& 0.0389& 0.0102& 0.0345\tabularnewline \hline 225& 0.0035&
0.0228& 0.0060& 0.0293\tabularnewline \hline
\end{tabular}
\end{center}
\end{table}

\begin{table}
\begin{center}
\caption{Fine Mesh Error, high conductivity channel.}
\label{ferrp6tindepsp}
\begin{tabular}{|c|c|c|c|c|}
\hline dof& $L^{1}$& $L^{\infty}$& $L^{2}$& $H^{1}$\tabularnewline
\hline 9& 0.0178& 0.0564& 0.0257& 0.0947\tabularnewline \hline 49&
0.0079& 0.0388& 0.0129& 0.0660\tabularnewline \hline 225& 0.0040&
0.0243& 0.0070& 0.0485\tabularnewline \hline
\end{tabular}
\end{center}
\end{table}

\clearpage

\subsection{Time dependent examples.}
In the following examples we consider media characterized by a
continuum of time scales. In the following examples the ODE obtained
on the coarse mesh from the homogenization of \eref{ghjh52} have
also been homogenized in time according to the method described in
subsection \ref{jksjhs89b}.

\begin{example}
\label{exa:tdep_trignometric}Time Dependent Multiscale
trigonometric.
\end{example}
In this example $a$ is given by equation \eref{ksjhskjhdhkdjh872y2}.
Although the number fine time steps to solve \eref{ghjh52} is
$2663$, only $134$ coarse time steps have been used to solve the
homogenized equation. Hence if  one also takes into account
homogenization in space, the compression factor is of the order of
$35000$ for the coarse mesh with $9$ interior nodes.

\begin{figure}[httb]
  \begin{center}
    \subfigure[$t\rightarrow a(x_0,t)$.]
    {\includegraphics[width=0.35\textwidth,height= 0.3\textwidth]{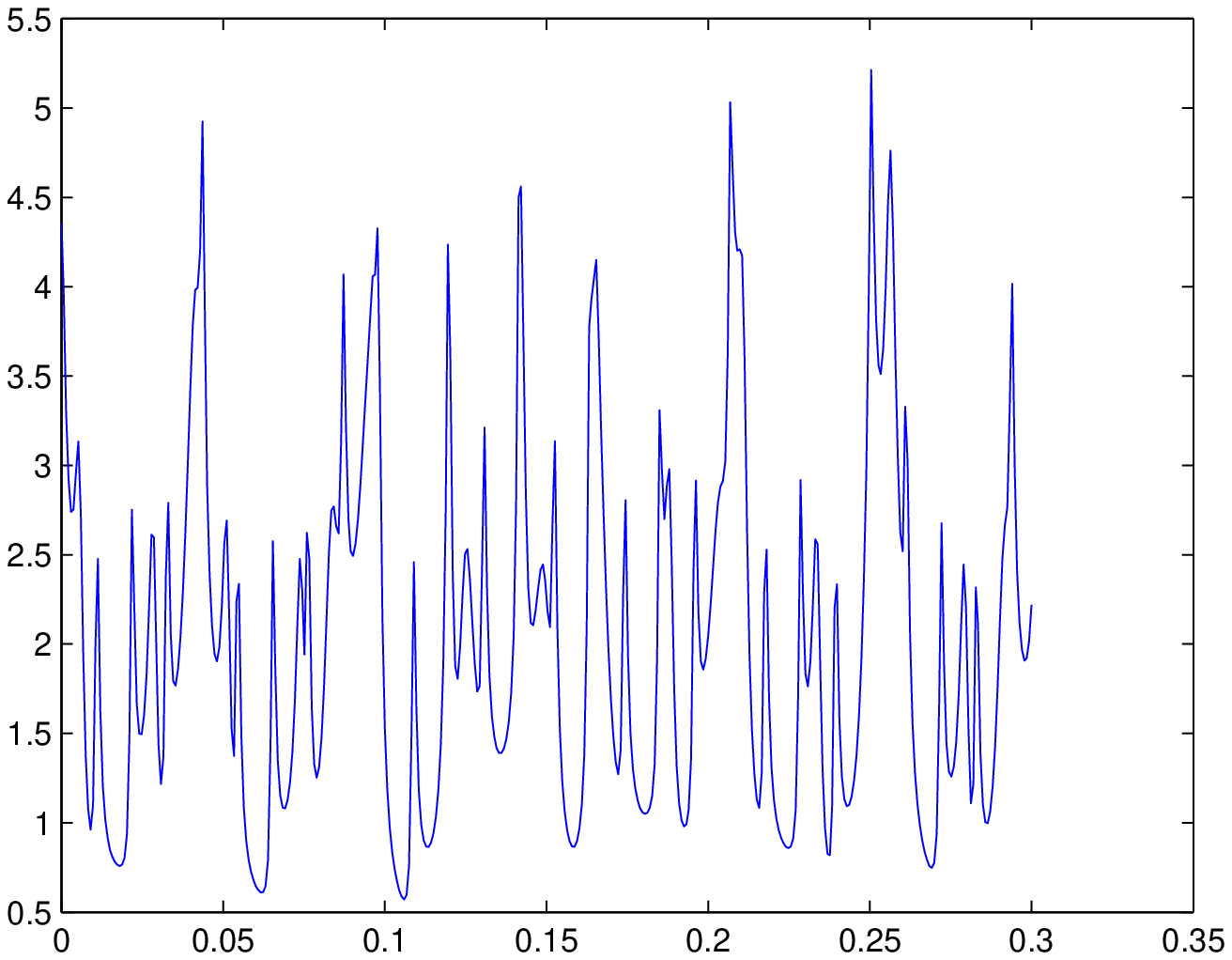}}
    \goodgap
    \subfigure[Top view of $t\rightarrow F(x_0,t)$.]
    {\includegraphics[width=0.35\textwidth,height= 0.3\textwidth]{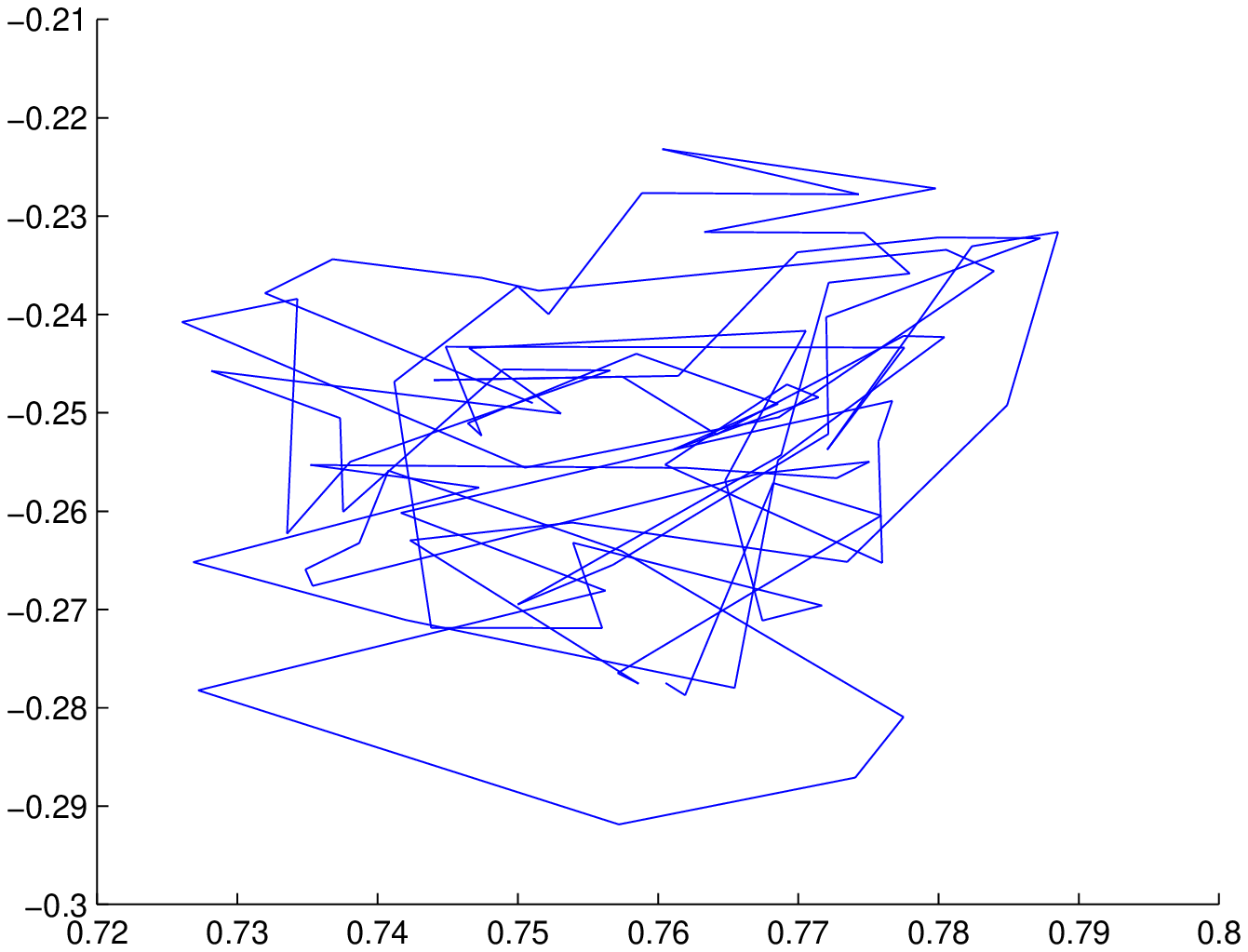}}\\
    \caption{Multiscale time dependent trigonometric medium.}
    \label{fixpointF}
\end{center}
\end{figure}

Figure \ref{fixpointF} shows the curves of $t\rightarrow a(x_0,t)$
and $t\rightarrow F(x_0,t)$ for a given $x_0\in \Omega$.

\begin{figure}[httb]
  \begin{center}
    \subfigure[Coarse mesh $L^1$ error.]
    {\includegraphics[width=0.35\textwidth,height= 0.3\textwidth]{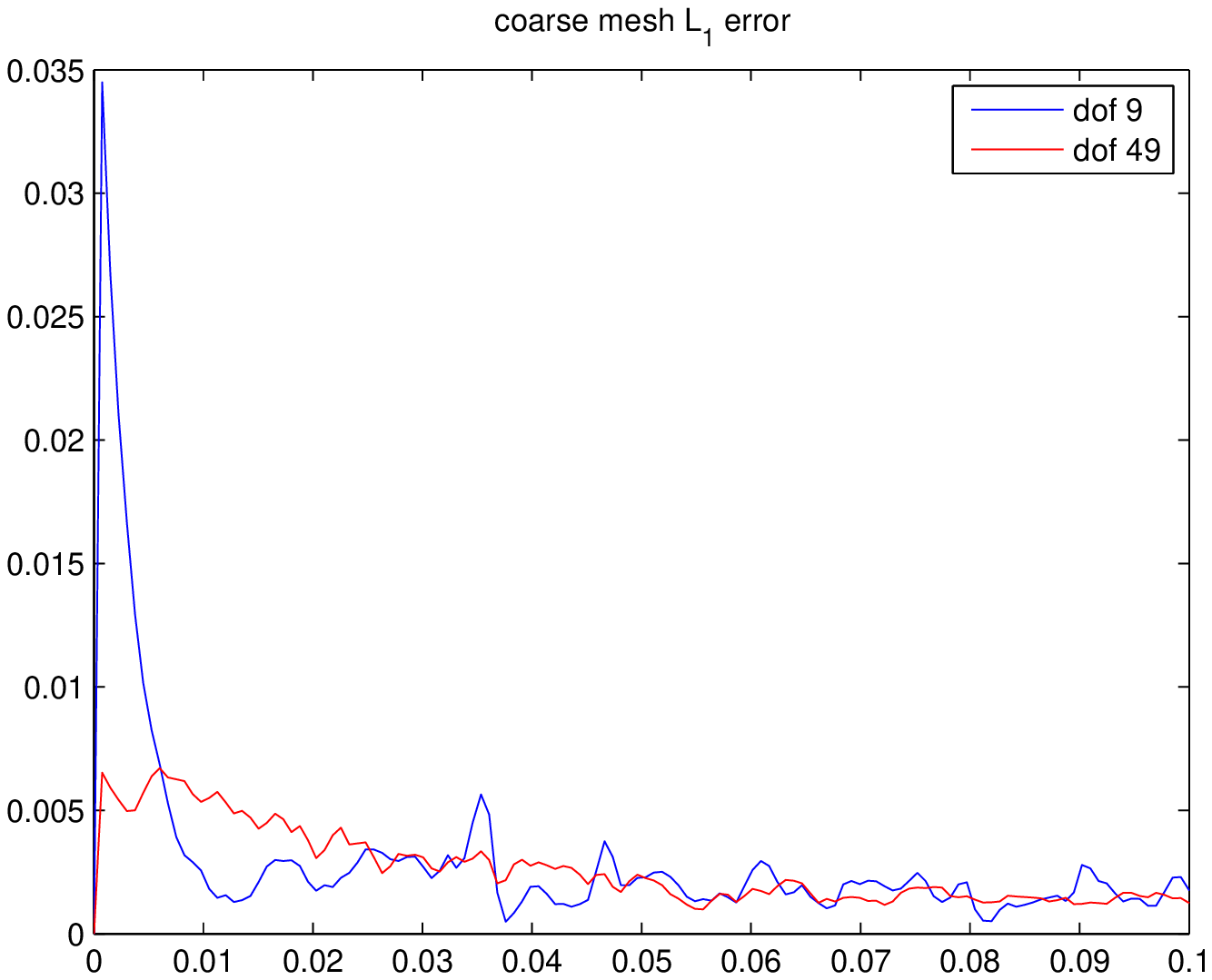}}
    \goodgap
    \subfigure[Fine mesh $L^1$ error.]
    {\includegraphics[width=0.35\textwidth,height= 0.3\textwidth]{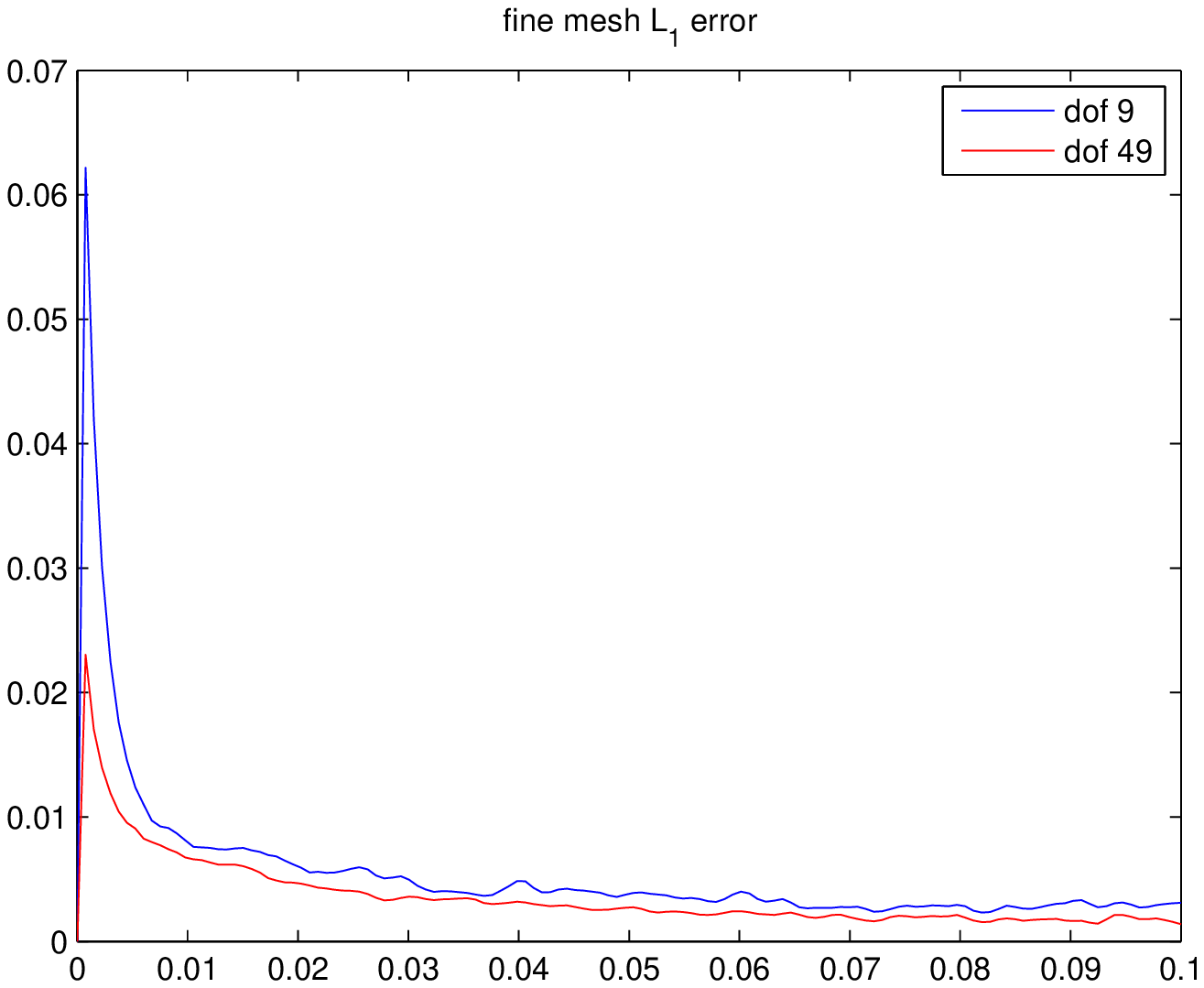}}\\
    \caption{$L^1$ error. Multiscale Trigonometric time dependent Medium.}
    \label{errL1tdepp4}
\end{center}
\end{figure}

\begin{figure}[httb]
  \begin{center}
    \subfigure[Coarse mesh $L^2$ error.]
    {\includegraphics[width=0.35\textwidth,height= 0.3\textwidth]{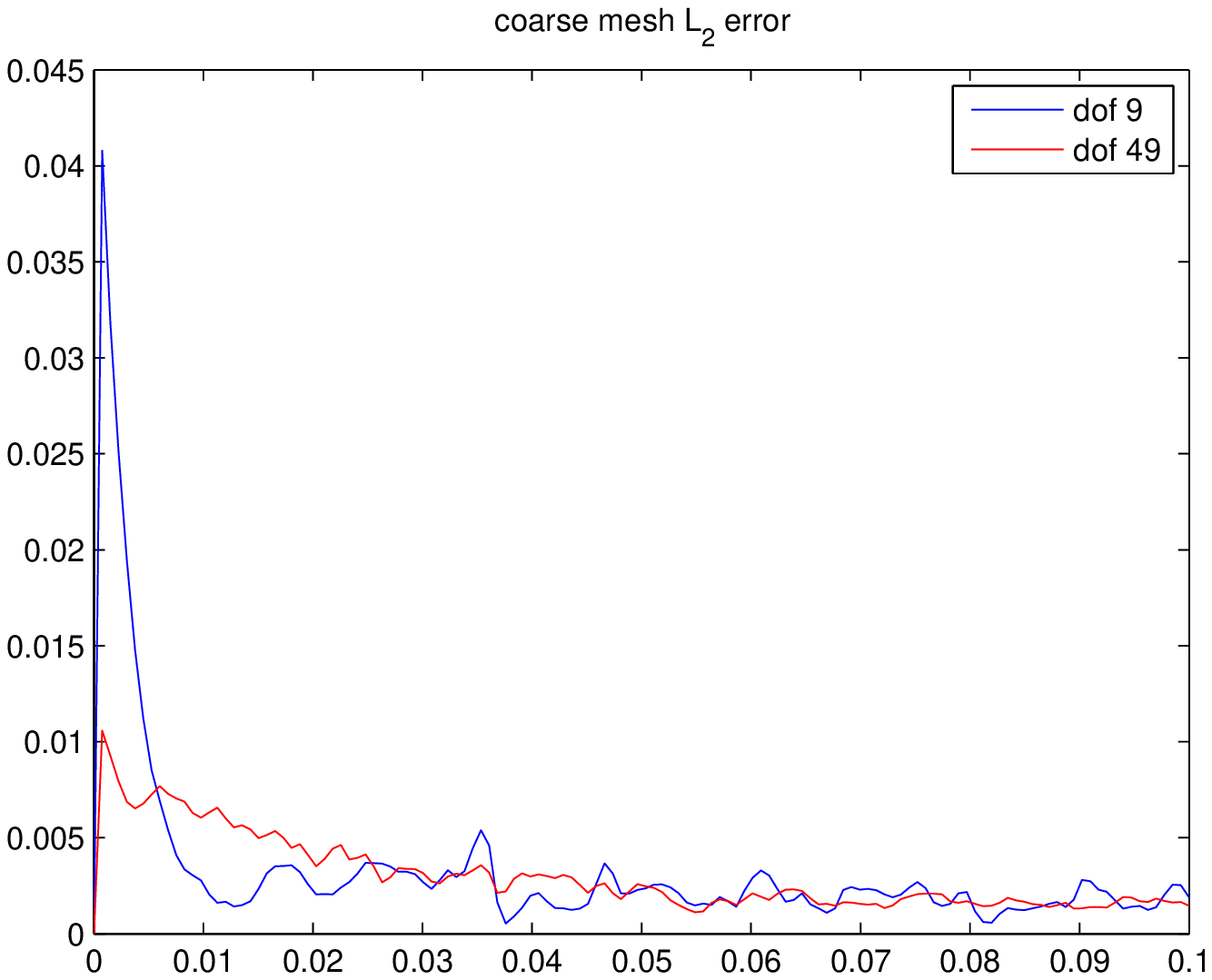}}
    \goodgap
    \subfigure[Fine Mesh $L^2$ error.]
    {\includegraphics[width=0.35\textwidth,height= 0.3\textwidth]{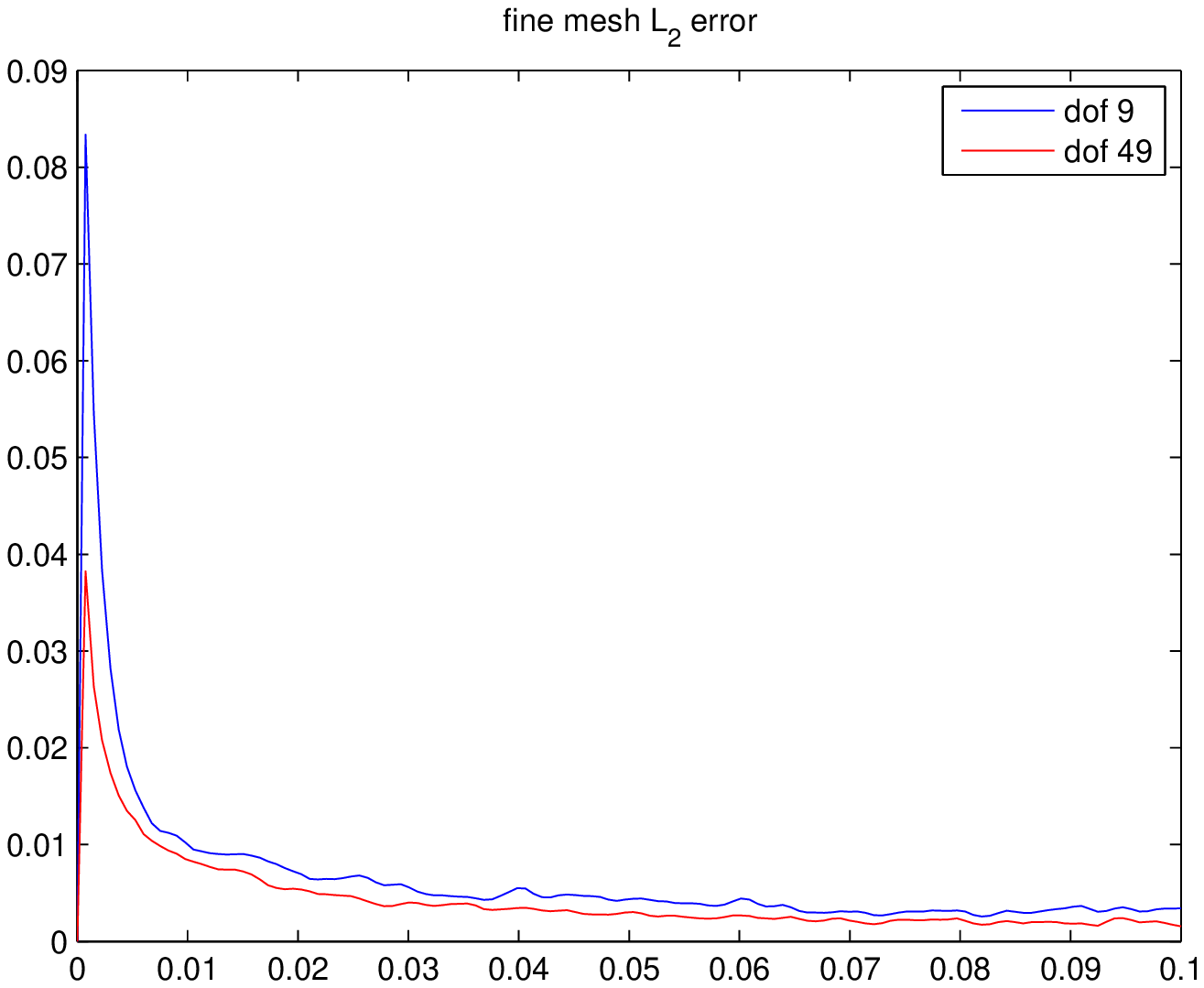}}\\
    \caption{$L^2$ error. Multiscale Trigonometric time dependent Medium.}
    \label{errL2tdepp4}
\end{center}
\end{figure}

\begin{figure}[httb]
  \begin{center}
    \subfigure[Coarse mesh $L_{\infty}$ error.]
    {\includegraphics[width=0.35\textwidth,height= 0.3\textwidth]{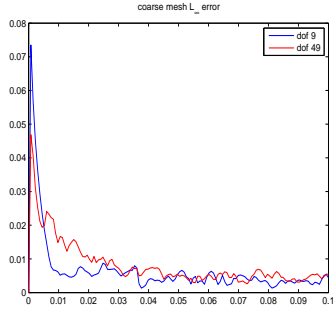}}
    \goodgap
    \subfigure[Fine mesh $L_{\infty}$ error.]
    {\includegraphics[width=0.35\textwidth,height= 0.3\textwidth]{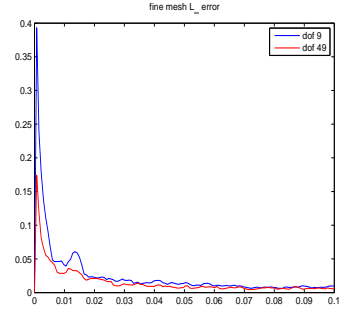}}\\
    \caption{$L^{\infty}$ error. Multiscale Trigonometric time dependent Medium.}
    \label{errLitdepp4}
\end{center}
\end{figure}

\begin{figure}[httb]
  \begin{center}
    \subfigure[Coarse Mesh $H^1$ error. \label{cH1p4}]
    {\includegraphics[width=0.35\textwidth,height= 0.3\textwidth]{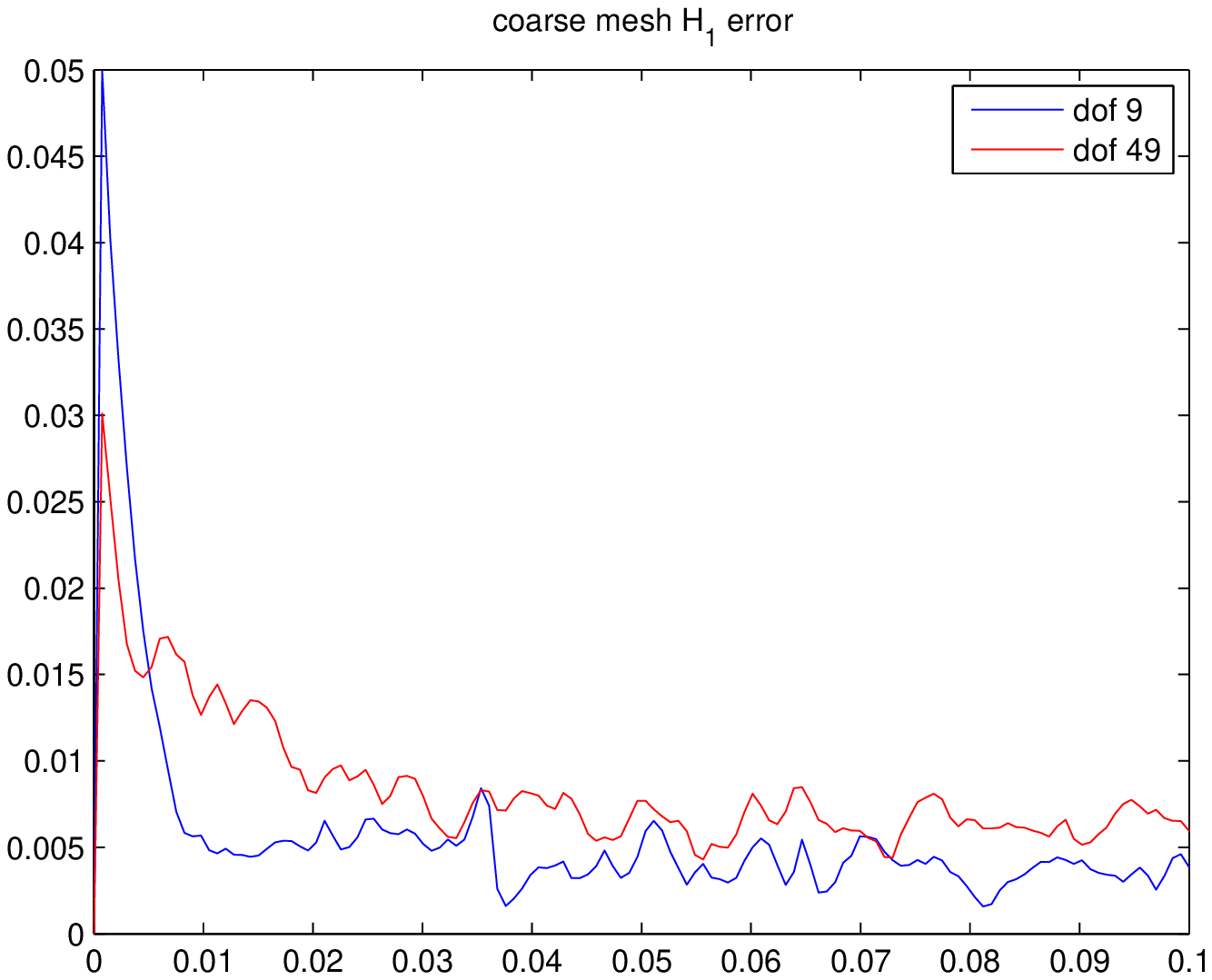}}
    \goodgap
    \subfigure[Fine Mesh $H^1$ error. \label{fH1p4}]
    {\includegraphics[width=0.35\textwidth,height= 0.3\textwidth]{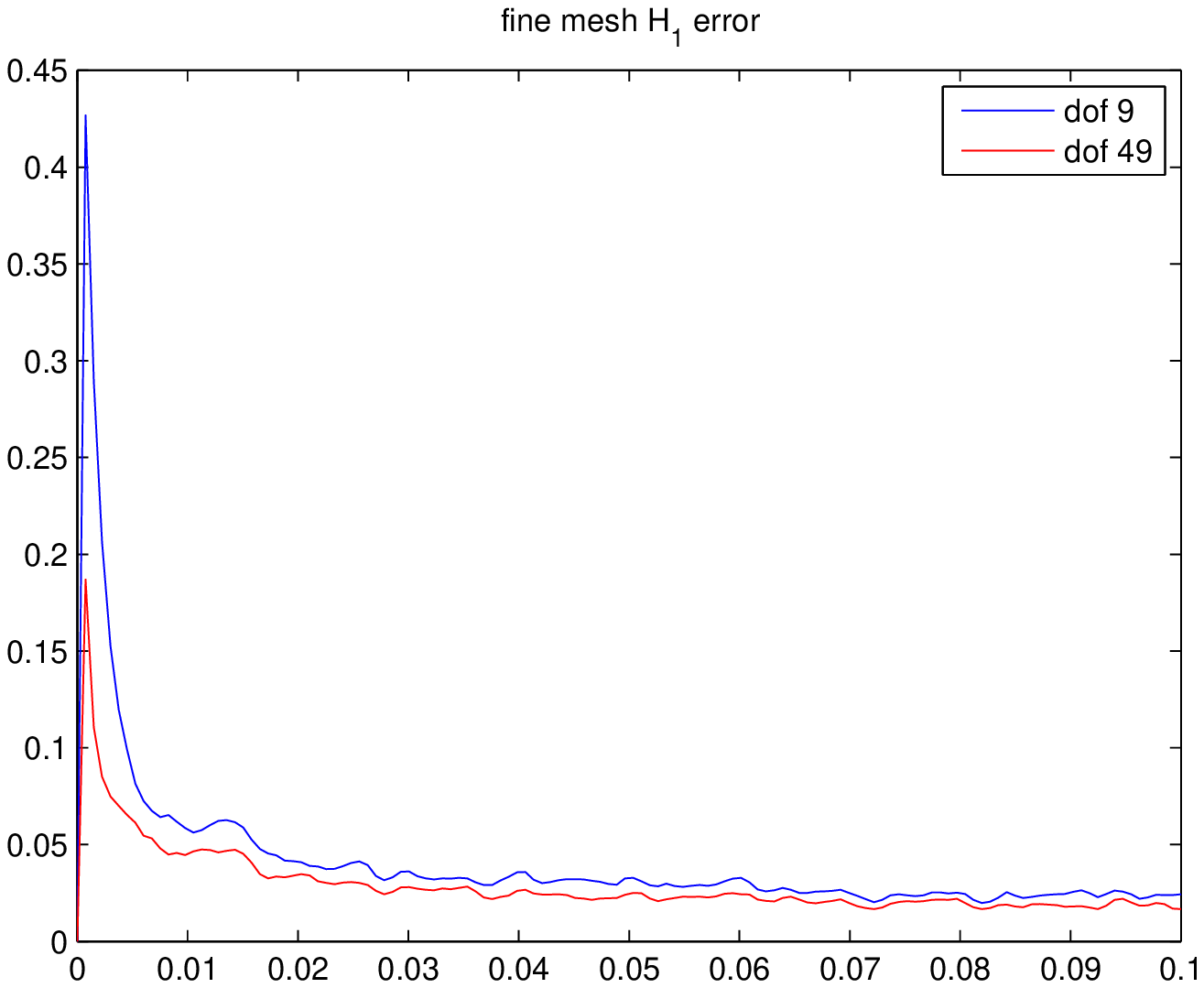}}\\
    \caption{$H^1$ error. Multiscale Trigonometric time dependent Medium.}
    \label{errH1tdepp4}
\end{center}
\end{figure}

The coarse  and fine mesh relative $L^1$, $L^2$, $L^\infty$, and
$H^1$ errors with respect to time have been plotted in figures
\ref{errL1tdepp4}, \ref{errL2tdepp4}, \ref{errLitdepp4} and
\ref{errH1tdepp4}. The initial increase of the relative error has
its origin in the initial value $u\equiv 0$ at time $0$.

The coarse and fine meshes errors are given in tables
\ref{cerrtdp4sp} and \ref{ferrtdp4sp} for $g=1$ at $t=0.1$, those
errors are given in tables \ref{cerrtdrhsp4} and \ref{ferrtdrhsp4}
for $g=\sin(2.4 x-1.8y+2 \pi t)$ at $t=0.1$

\begin{table}
\begin{center}
\caption{Coarse Mesh Error. Multiscale trigonometric time dependent
Medium. $g=1$.} \label{cerrtdp4sp}
\begin{tabular}{|c|c|c|c|c|}
\hline dof& $L^{1}$& $L^{\infty}$& $L^{2}$& $H^{1}$\tabularnewline
\hline 9& 0.0018& 0.0045& 0.0019& 0.0039\tabularnewline \hline 49&
0.0012& 0.0054& 0.0015& 0.0060\tabularnewline \hline
\end{tabular}
\end{center}
\end{table}

\begin{table}
\begin{center}
\caption{Fine Mesh Error.  Multiscale trigonometric time dependent
Medium. $g=1$.} \label{ferrtdp4sp}
\begin{tabular}{|c|c|c|c|c|}
\hline dof& $L^{1}$& $L^{\infty}$& $L^{2}$& $H^{1}$\tabularnewline
\hline 9& 0.0031& 0.0096& 0.0034& 0.0242\tabularnewline \hline 49&
0.0014& 0.0059& 0.0016& 0.0166\tabularnewline \hline
\end{tabular}
\end{center}
\end{table}

\begin{table}
\begin{center}
\caption{Coarse mesh error. Multiscale trigonometric time dependent
Medium. $g=\sin(2.4 x-1.8y+2 \pi t)$. } \label{cerrtdrhsp4}
\begin{tabular}{|c|c|c|c|c|}
\hline dof& $L^{1}$& $L^{\infty}$& $L^{2}$& $H^{1}$\tabularnewline
\hline 9& 0.0043& 0.0087& 0.0044& 0.0085\tabularnewline \hline 49&
0.0033& 0.0079& 0.0035& 0.0084\tabularnewline \hline
\end{tabular}
\end{center}
\end{table}

\begin{table}
\begin{center}
\caption{Fine mesh error. Multiscale trigonometric time dependent
medium. $g=\sin(2.4 x-1.8y+2 \pi t)$. } \label{ferrtdrhsp4}
\begin{tabular}{|c|c|c|c|c|}
\hline dof& $L^{1}$& $L^{\infty}$& $L^{2}$& $H^{1}$\tabularnewline
\hline 9& 0.0082& 0.0199& 0.0087& 0.0379\tabularnewline \hline 49&
0.0038& 0.0104& 0.0040& 0.0244\tabularnewline \hline
\end{tabular}
\end{center}
\end{table}

\clearpage

\begin{example}
\label{exa:tdep_fourier}Time Dependent Random Fourier Modes.
\end{example}

In this example  $a(x,y,t)=e^{h(x,y,t)}$ where $h$ is given by the
following equation
\[ h(x,y,t)=\sum_{|k|\leq R}(a_{k}\sin(2\pi
k.x')+b_{k}\cos(2\pi k.x'))\] where $R=4$, $x'=x+\sqrt{2}t$,
$y'=y-\sqrt{2}t$, $a_{k}$ and $b_{k}$ are independent identically
distributed random variables on $[-0.2,0.2]$. In this example, one
can compute that $\frac{\lambda_{\max}(a)}{\lambda_{\min}(a)}=95.7$.
\begin{figure}[httb]
\begin{center}
\includegraphics[%
  scale=0.3]{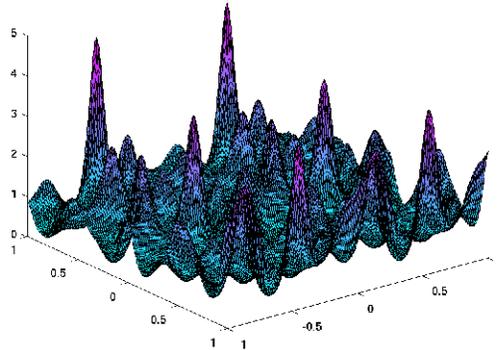}
\caption{Random Fourier Mode} \label{mediap5}
\end{center}
\end{figure}

Figure \ref{mediap5} is a plot of $a$ at time $0$. Although the
number fine time steps to solve \eref{ghjh52} is $1332$ only $67$
coarse time steps have been used to solve the homogenized equation.
Hence if  one also takes into account homogenization in space, the
compression factor is of the order of $35000$ for the coarse mesh
with $9$ interior nodes.

\begin{figure}[httb]
  \begin{center}
    \subfigure[coarse mesh $L^1$ error.]
    {\includegraphics[width=0.35\textwidth,height= 0.3\textwidth]{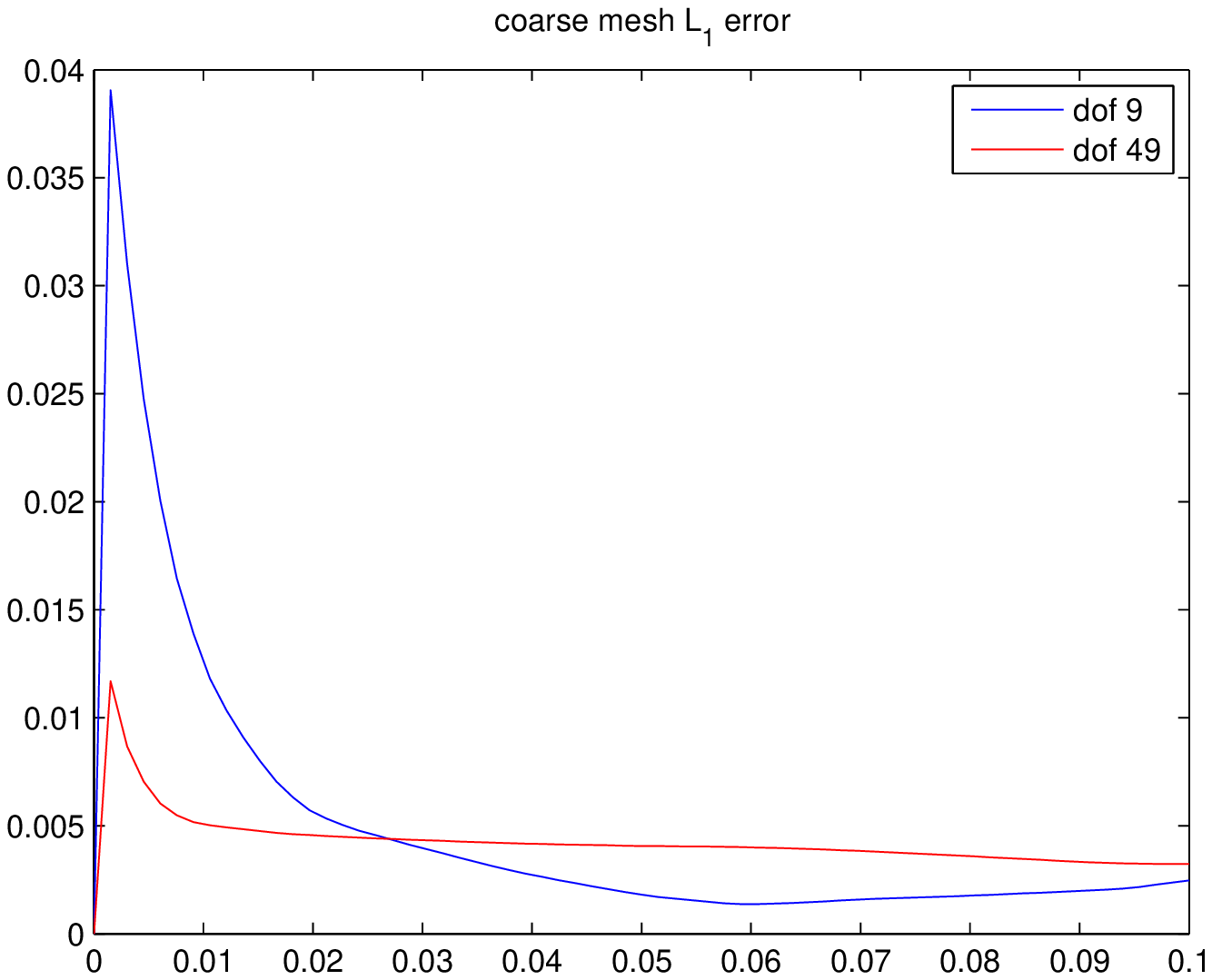}}
    \goodgap
    \subfigure[fine mesh $L^1$ error.]
    {\includegraphics[width=0.35\textwidth,height= 0.3\textwidth]{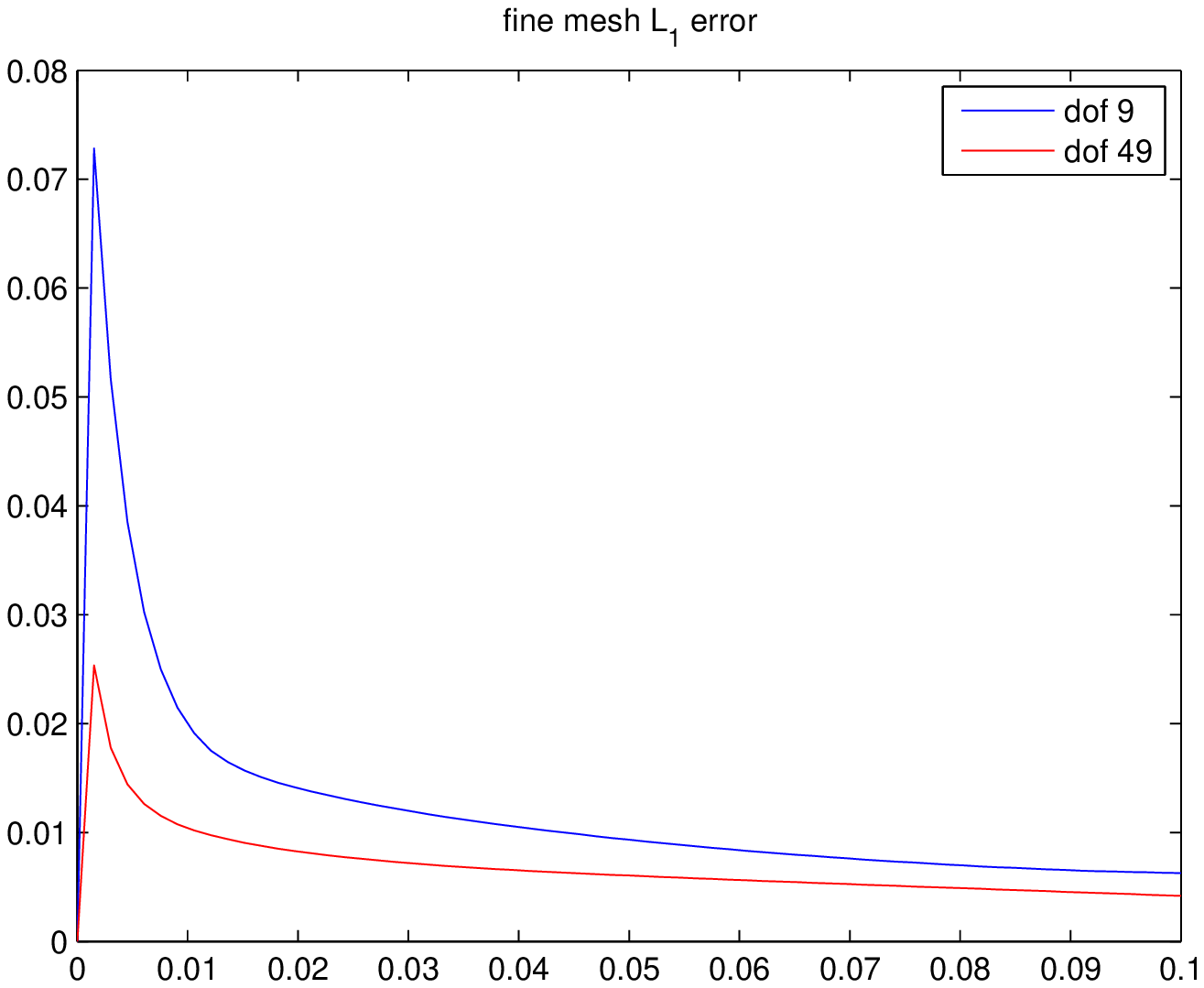}}\\
    \caption{$L^1$ error. Random Fourier modes up to  $t=.1$.}
    \label{errL1tdepp5}
\end{center}
\end{figure}

\begin{figure}[httb]
  \begin{center}
    \subfigure[coarse mesh $L^2$ error.]
    {\includegraphics[width=0.35\textwidth,height= 0.3\textwidth]{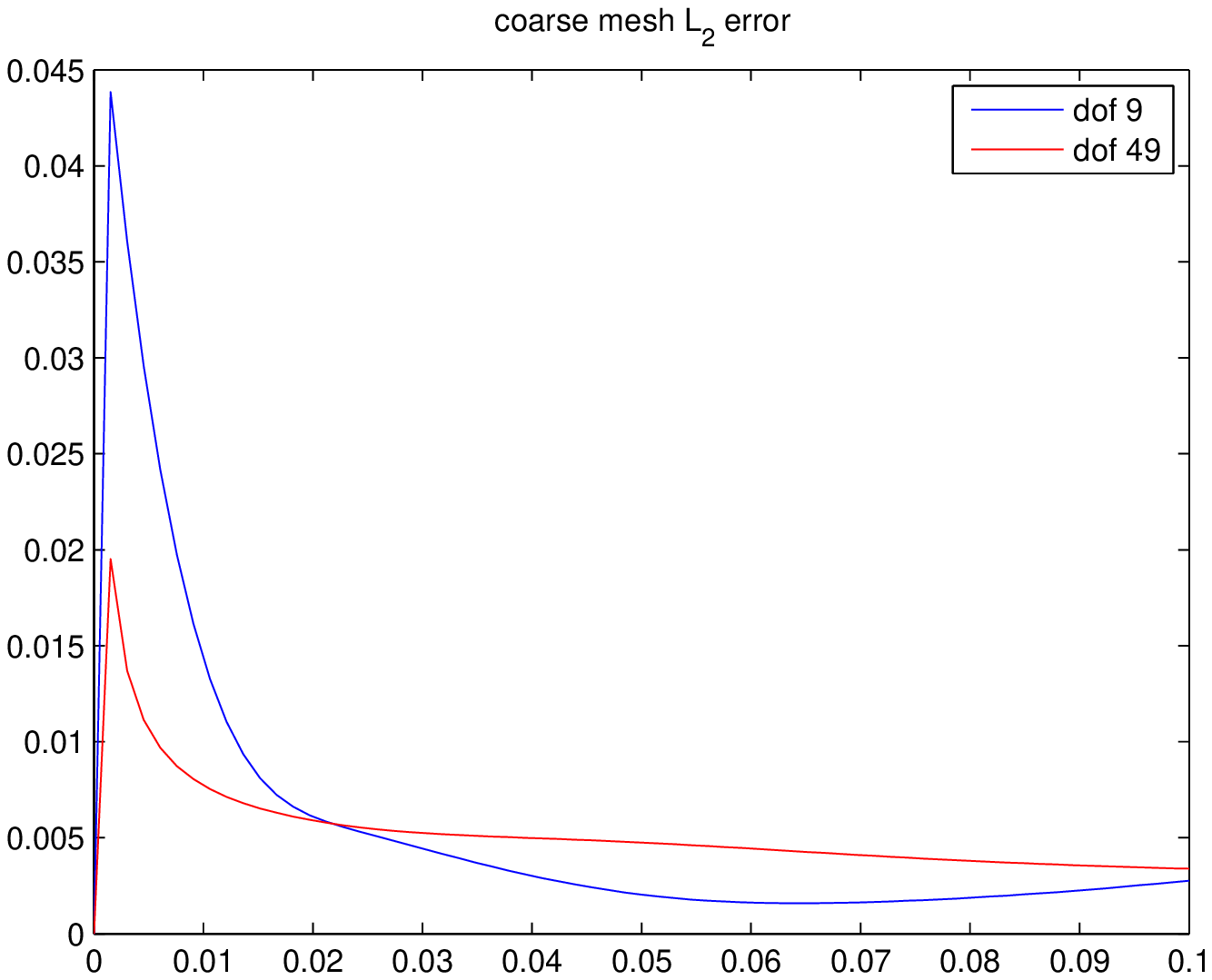}}
    \goodgap
    \subfigure[Fine Mesh $L^2$ error.]
    {\includegraphics[width=0.35\textwidth,height= 0.3\textwidth]{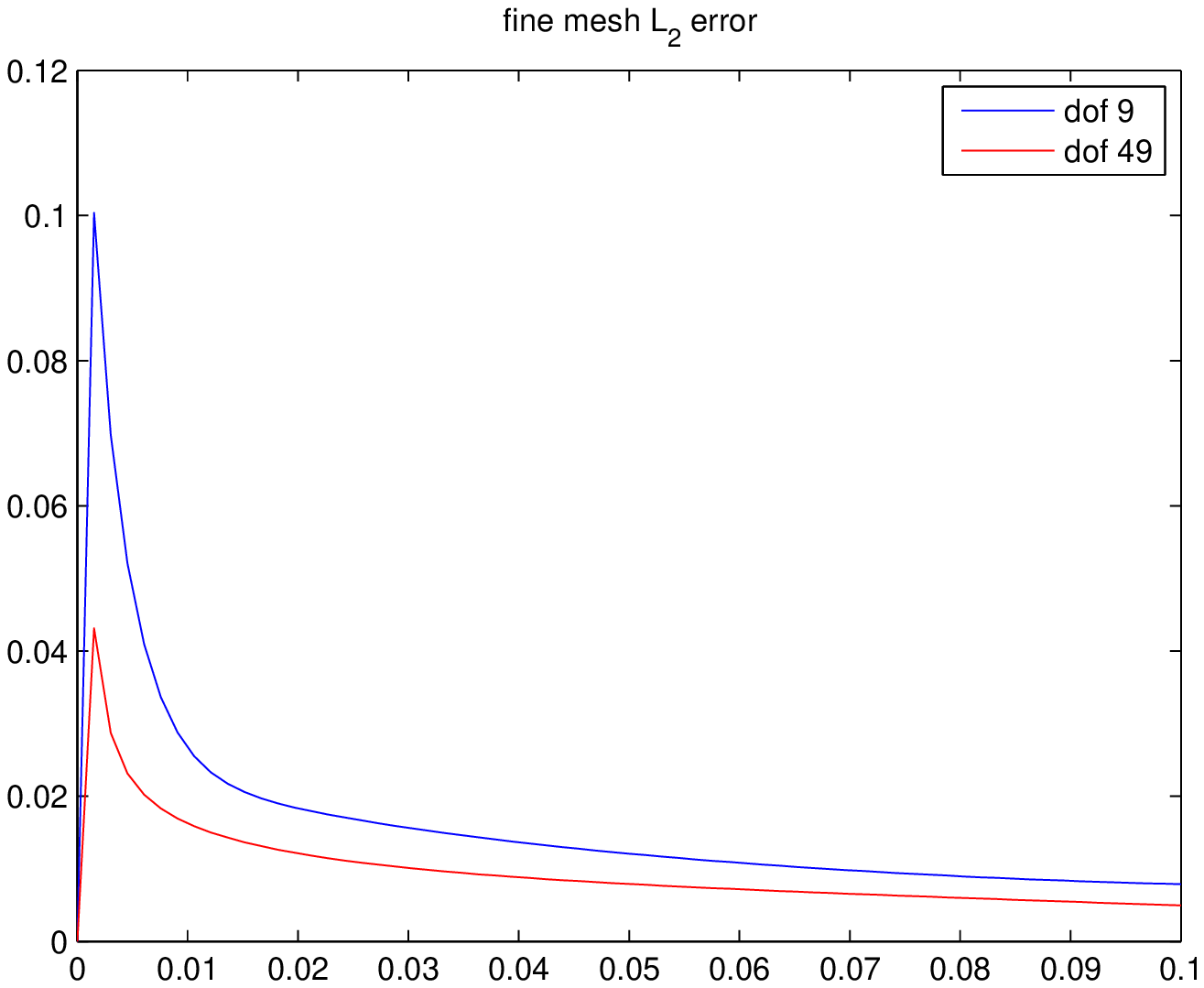}}\\
    \caption{$L^2$ error. Random Fourier modes up to $t=.1$..}
    \label{errL2tdepp5}
\end{center}
\end{figure}

\begin{figure}[httb]
  \begin{center}
    \subfigure[coarse mesh $L_{\infty}$ error.]
    {\includegraphics[width=0.35\textwidth,height= 0.3\textwidth]{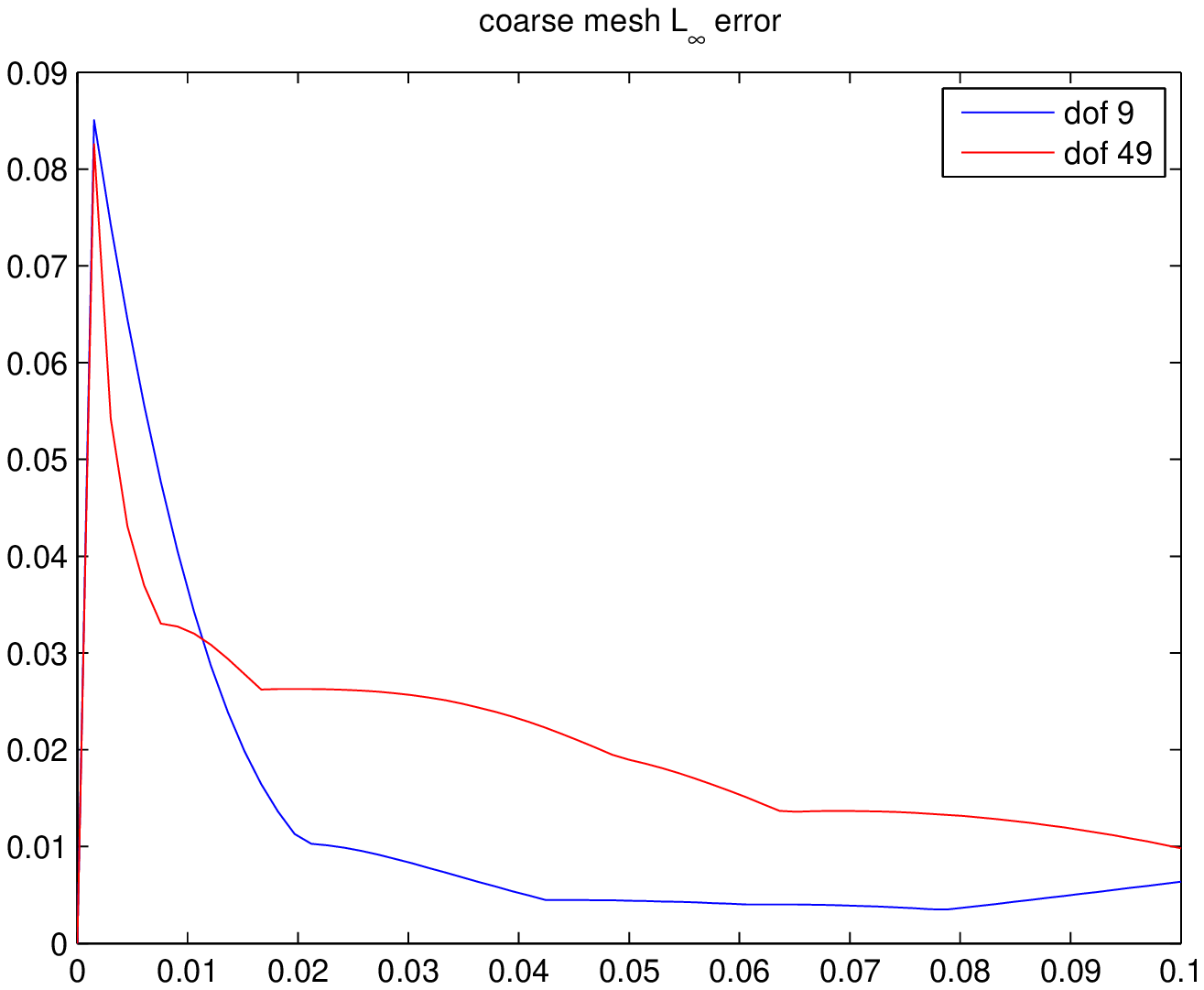}}
    \goodgap
    \subfigure[fine mesh $L_{\infty}$ error.]
    {\includegraphics[width=0.35\textwidth,height= 0.3\textwidth]{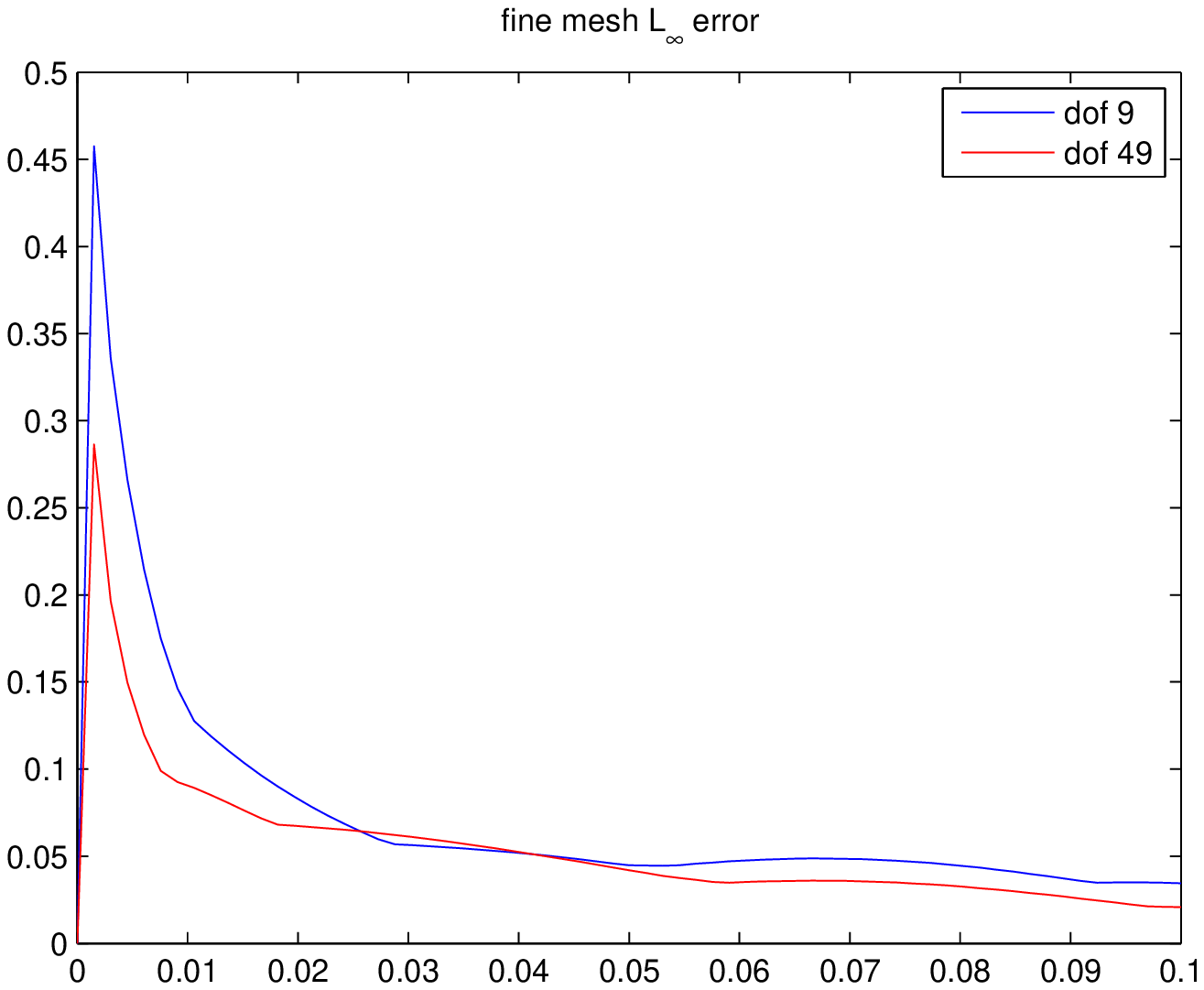}}\\
    \caption{$L^{\infty}$ error. Random Fourier modes up to $t=.1$..}
    \label{errLitdepp5}
\end{center}
\end{figure}

\begin{figure}[httb]
  \begin{center}
    \subfigure[coarse mesh $H^1$ error.]
    {\includegraphics[width=0.35\textwidth,height= 0.3\textwidth]{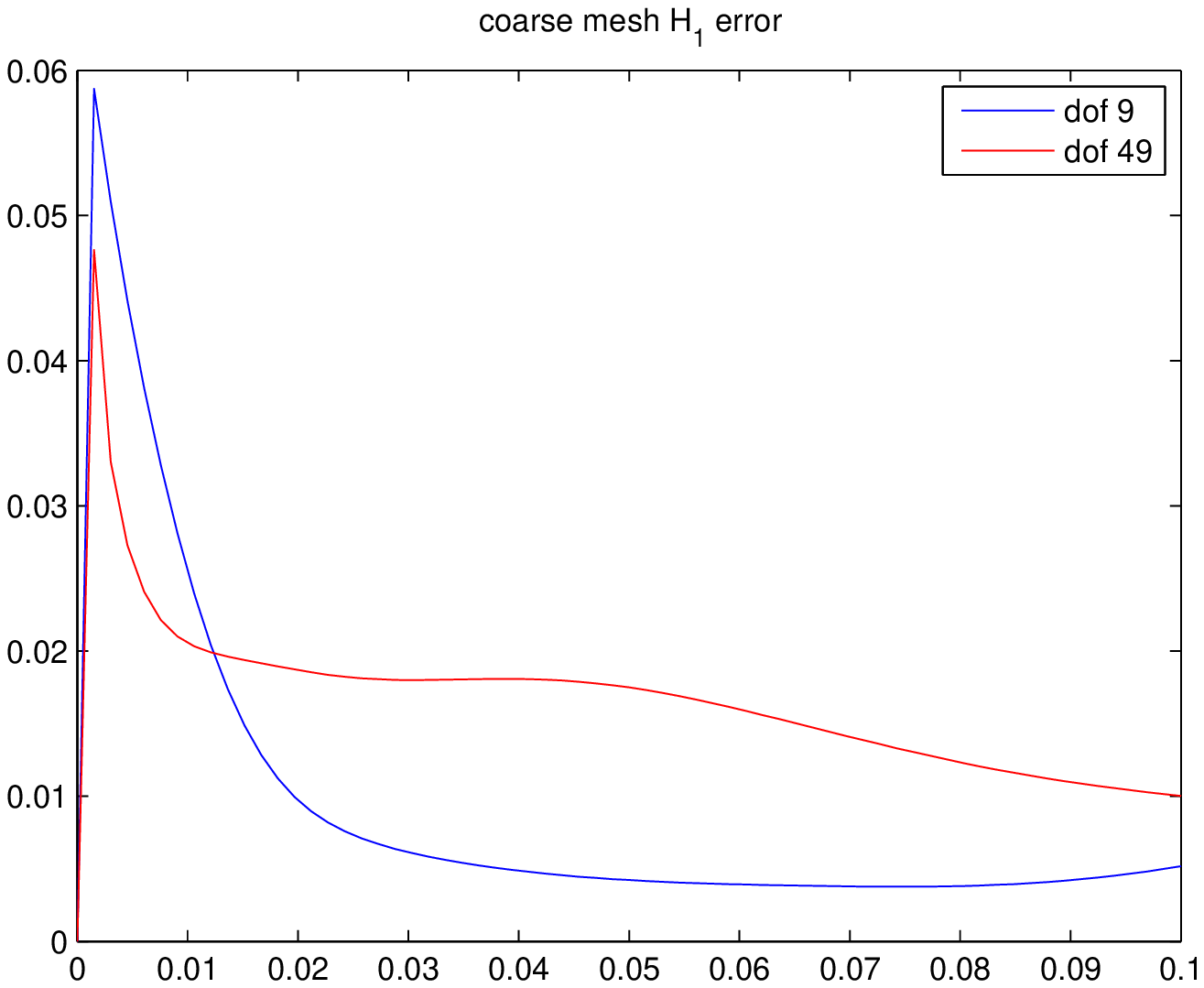}}
    \goodgap
    \subfigure[fine mesh $H^1$ error.]
    {\includegraphics[width=0.35\textwidth,height= 0.3\textwidth]{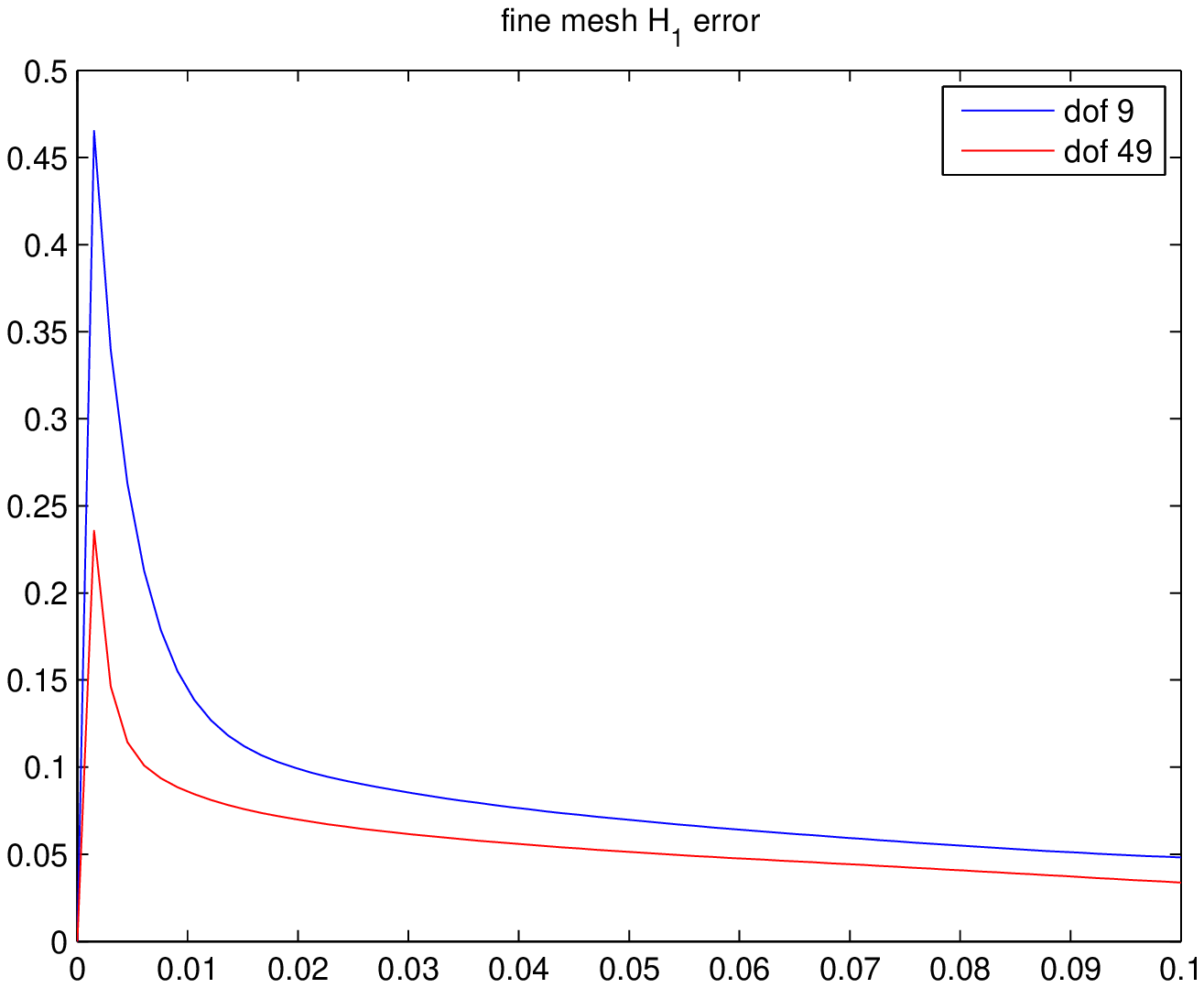}}\\
    \caption{$H^1$ error. Random Fourier modes up to $t=.1$..}
    \label{errH1tdepp5}
\end{center}
\end{figure}

\begin{figure}[httb]
  \begin{center}
    \subfigure[coarse mesh $L^1$ error.]
    {\includegraphics[width=0.35\textwidth,height= 0.3\textwidth]{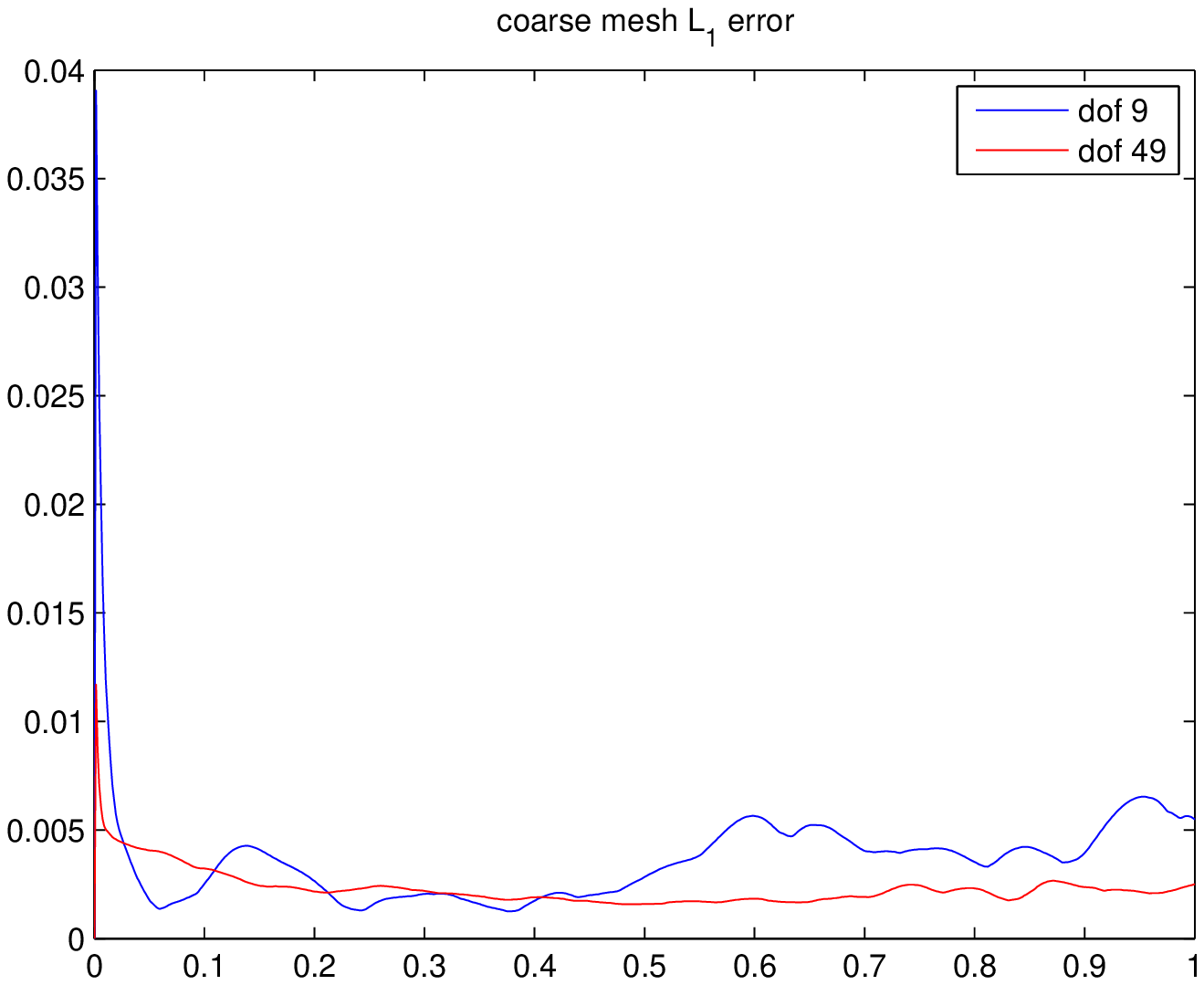}}
    \goodgap
    \subfigure[fine mesh $L^1$ error.]
    {\includegraphics[width=0.35\textwidth,height= 0.3\textwidth]{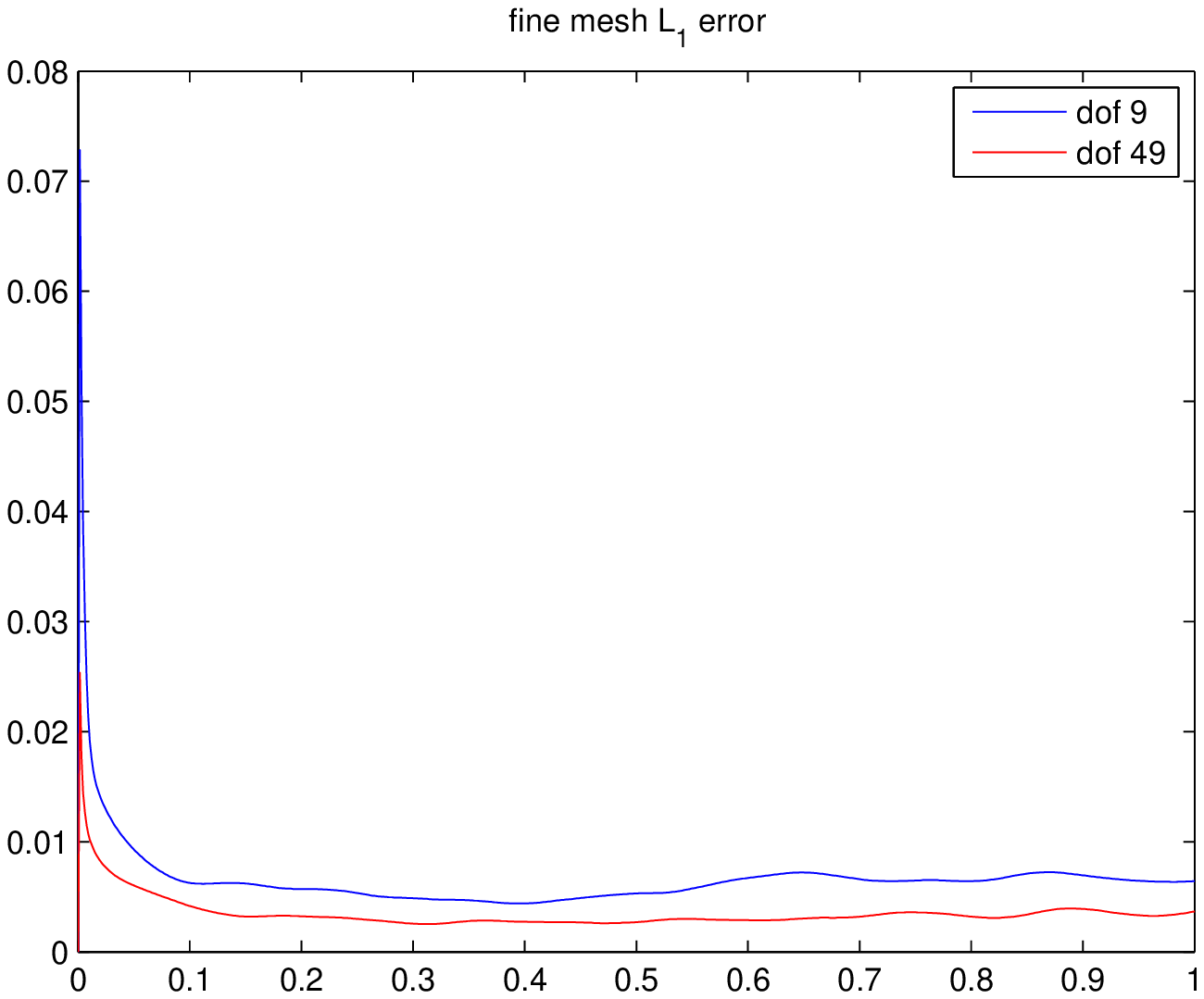}}\\
    \caption{$L^1$ error. Random Fourier modes up to $t=1$.}
    \label{errL1t100p5}
\end{center}
\end{figure}

\begin{figure}[httb]
  \begin{center}
    \subfigure[coarse mesh $L^2$ error.]
    {\includegraphics[width=0.35\textwidth,height= 0.3\textwidth]{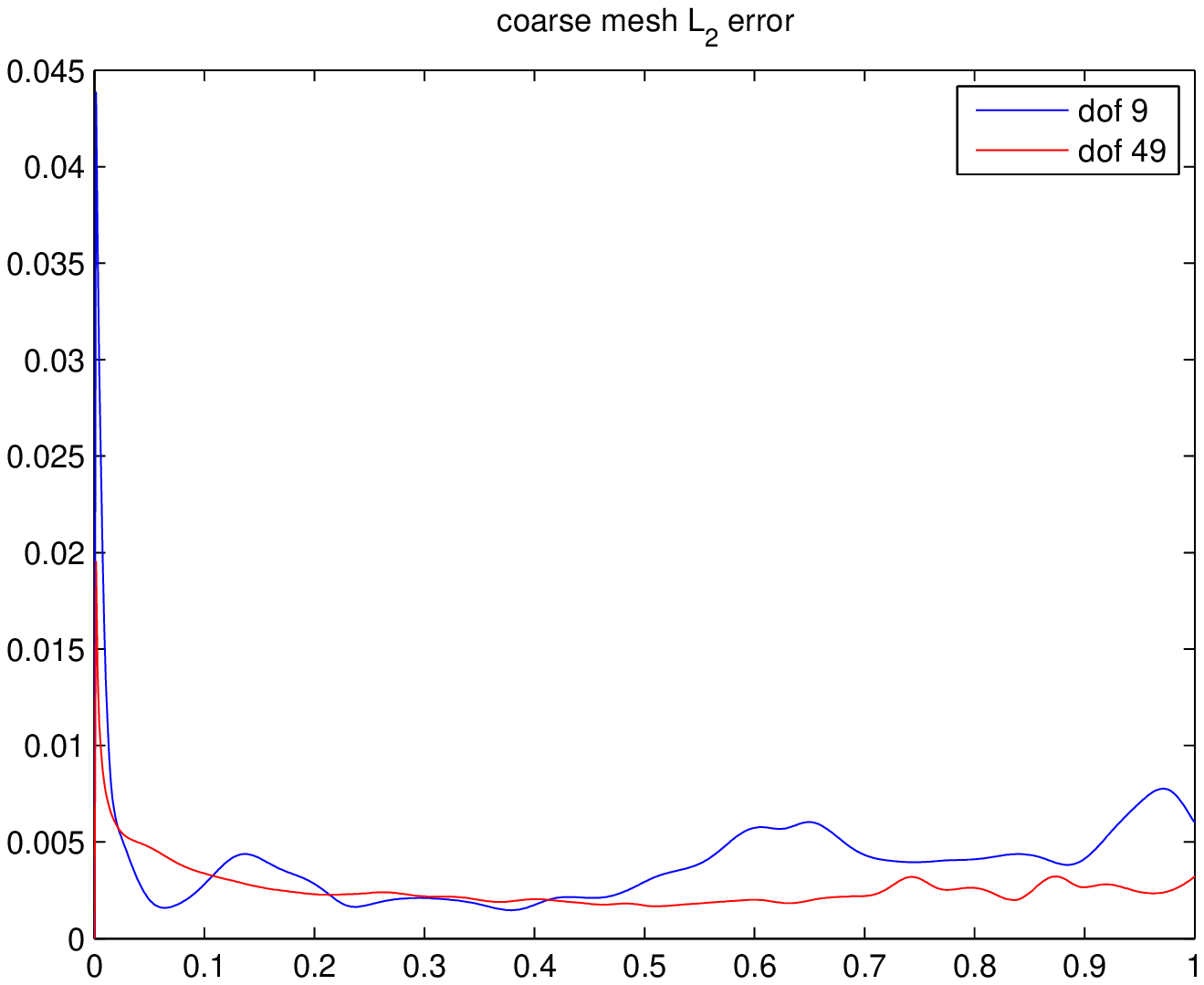}}
    \goodgap
    \subfigure[Fine Mesh $L^2$ error.]
    {\includegraphics[width=0.35\textwidth,height= 0.3\textwidth]{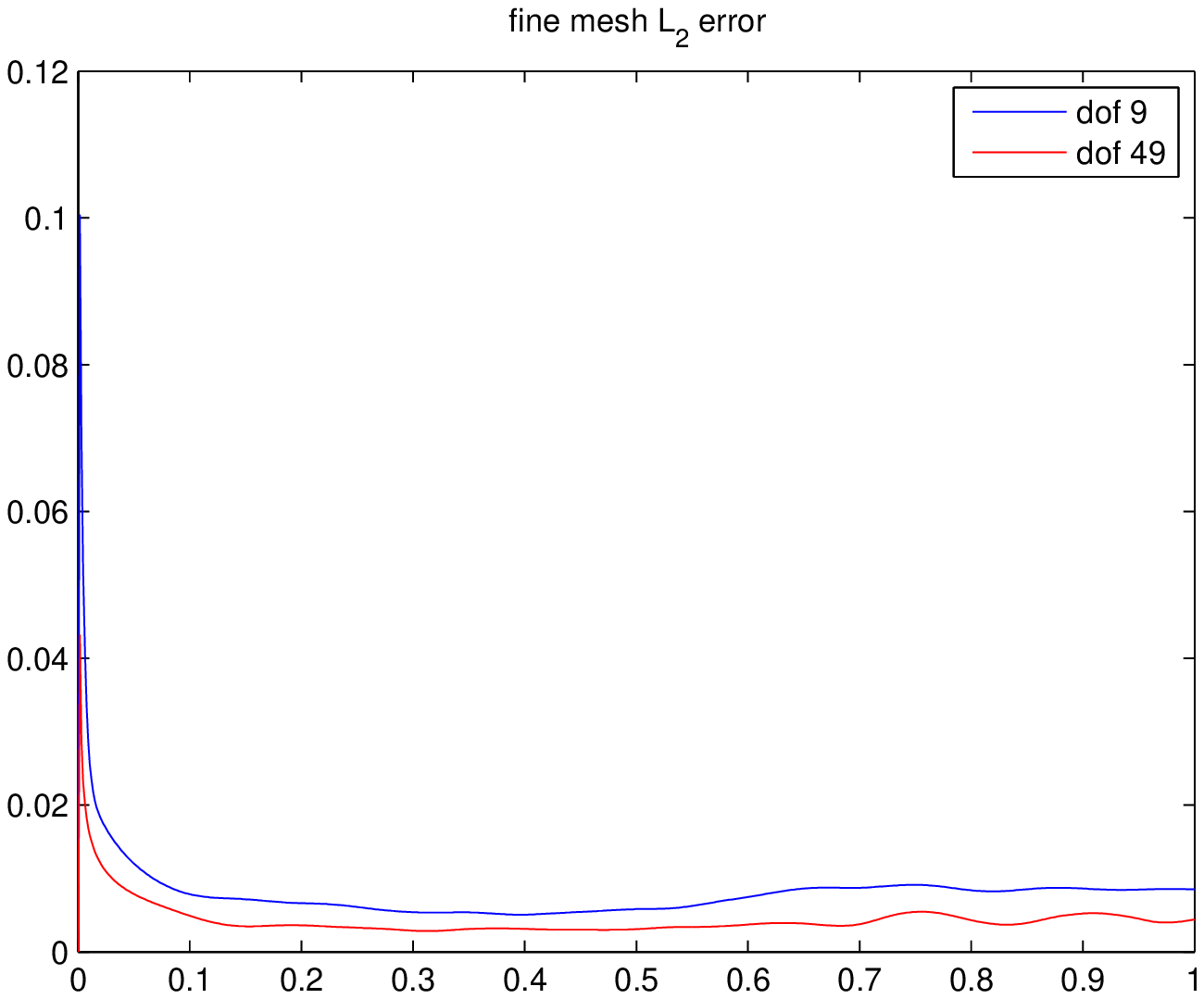}}\\
    \caption{$L^2$ error. Random Fourier modes up to $t=1$.}
    \label{errL2t100p5}
\end{center}
\end{figure}

\begin{figure}[httb]
  \begin{center}
    \subfigure[coarse mesh $L_{\infty}$ error.]
    {\includegraphics[width=0.35\textwidth,height= 0.3\textwidth]{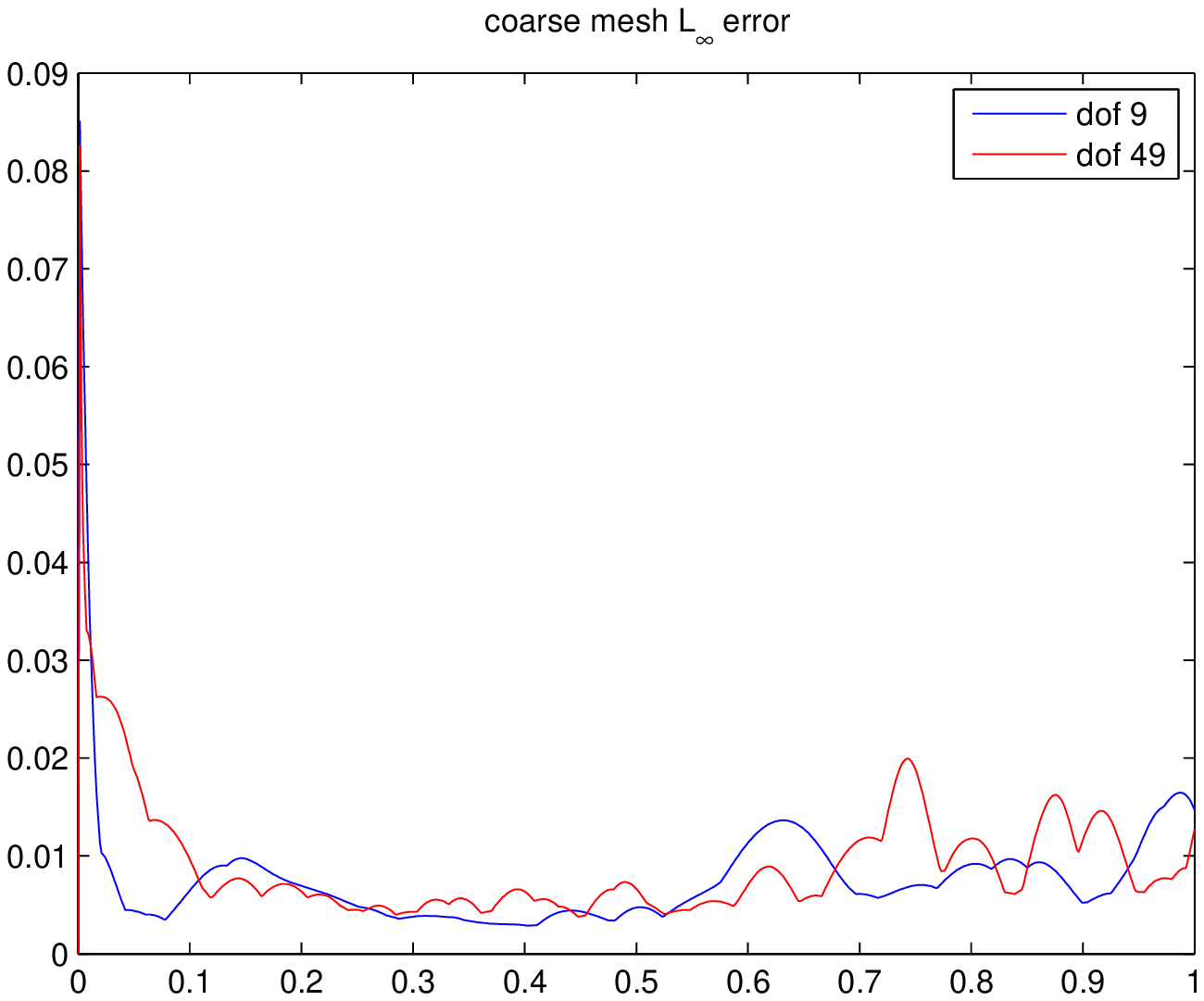}}
    \goodgap
    \subfigure[fine mesh $L_{\infty}$ error.]
    {\includegraphics[width=0.35\textwidth,height= 0.3\textwidth]{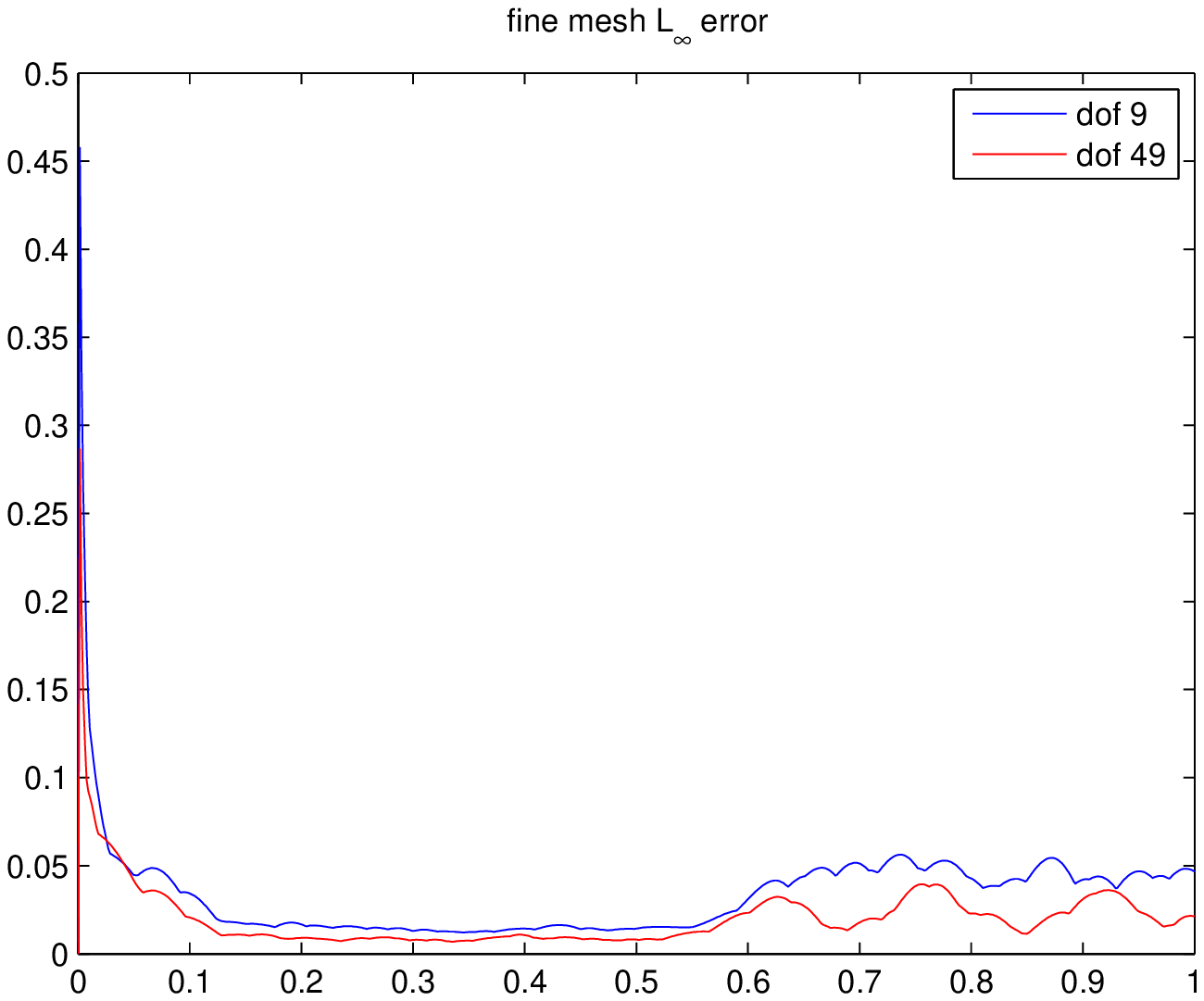}}\\
    \caption{$L^{\infty}$ error. Random Fourier modes up to $t=1$.}
    \label{errLit100p5}
\end{center}
\end{figure}

\begin{figure}[httb]
  \begin{center}
    \subfigure[coarse Mesh $H^1$ error. \label{cH1t1}]
    {\includegraphics[width=0.35\textwidth,height= 0.3\textwidth]{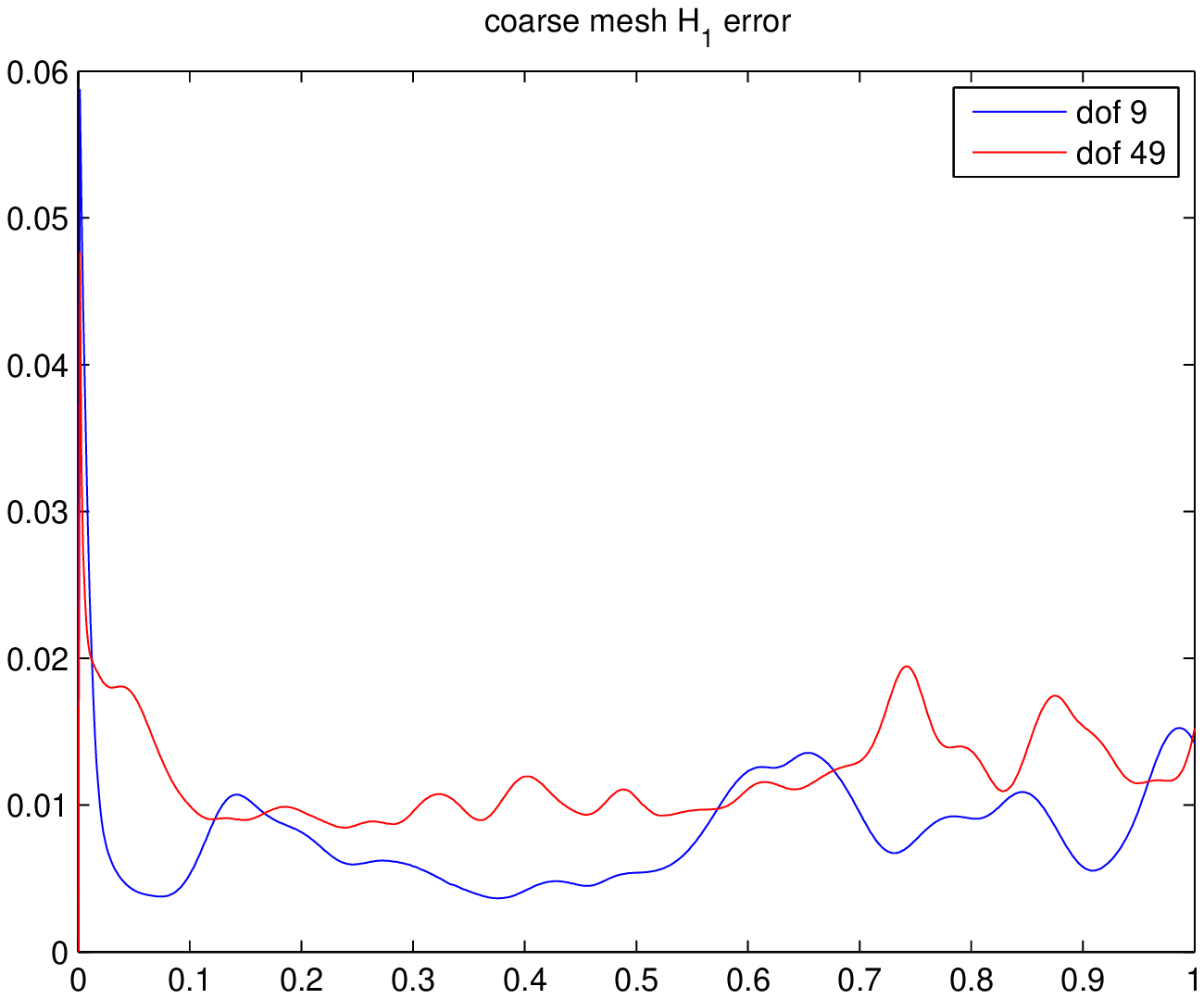}}
    \goodgap
    \subfigure[fine Mesh $H^1$ error. \label{fH1t1}]
    {\includegraphics[width=0.35\textwidth,height= 0.3\textwidth]{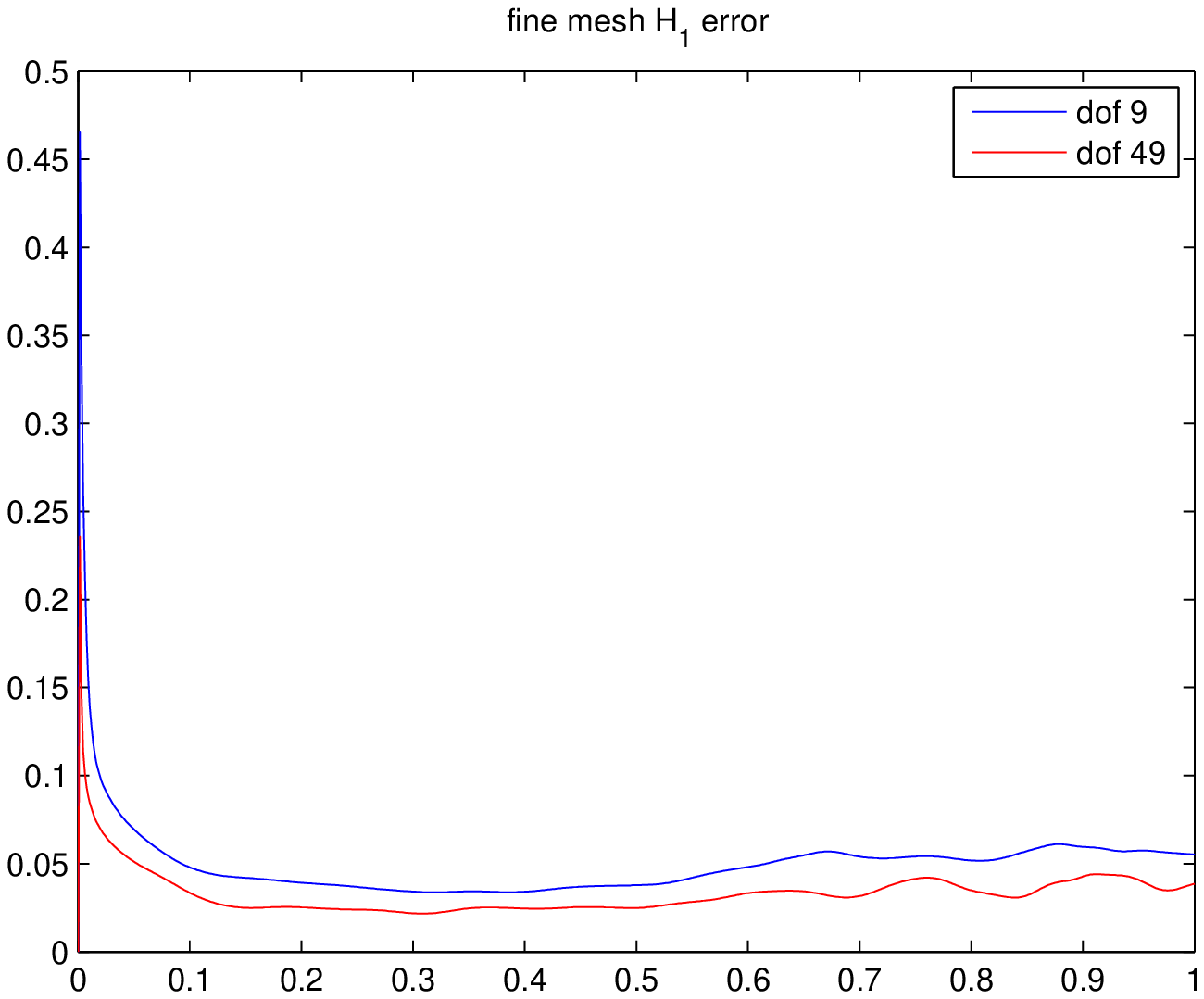}}\\
    \caption{$H^1$ error. Random Fourier modes up to $t=1$.}
    \label{errH1t100p5}
\end{center}
\end{figure}

The coarse  and fine mesh relative $L^1$, $L^2$, $L^\infty$, and
$H^1$ errors with respect to time (up to time $t=0.1$) have been
plotted in figures \ref{errL1tdepp5}, \ref{errL2tdepp5},
\ref{errLitdepp5} and \ref{errH1tdepp5}. Those errors are also given
up to $t=1$ in figures \ref{errL1t100p5}, \ref{errL2t100p5},
\ref{errLit100p5}, \ref{errH1t100p5}.

 The coarse and fine meshes errors are given in tables
\ref{cerrtdp5sp} and \ref{ferrtdp5sp} for $g=1$ at $t=0.1$.

\begin{table}
\begin{center}
\caption{Coarse mesh error at $t=0.1$. Time dependent random Fourier
modes.} \label{cerrtdp5sp}
\begin{tabular}{|c|c|c|c|c|}
\hline dof& $L^{1}$& $L^{\infty}$& $L^{2}$& $H^{1}$\tabularnewline
\hline 9& 0.0025& 0.0028& 0.0064& 0.0052\tabularnewline \hline 49&
0.0032& 0.0098& 0.0034& 0.0100\tabularnewline \hline
\end{tabular}
\end{center}
\end{table}

\begin{table}
\begin{center}
\caption{Fine mesh error at $t=0.1$. Time dependent random Fourier
modes.} \label{ferrtdp5sp}
\begin{tabular}{|c|c|c|c|c|}
\hline dof& $L^{1}$& $L^{\infty}$& $L^{2}$& $H^{1}$\tabularnewline
\hline 9& 0.0063& 0.0344& 0.0079& 0.0481\tabularnewline \hline 49&
0.0042& 0.0207& 0.0049& 0.0337\tabularnewline \hline
\end{tabular}
\end{center}
\end{table}

\clearpage

\begin{figure}[httb]
  \begin{center}
    \subfigure[$a$ at $t=0$.]
    {\includegraphics[width=0.35\textwidth,height= 0.3\textwidth]{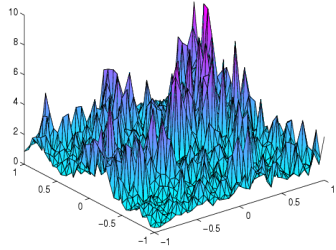}}
    \goodgap
    \subfigure[$(F_1,F_2)$ at $t=0$.]
    {\includegraphics[width=0.35\textwidth,height= 0.3\textwidth]{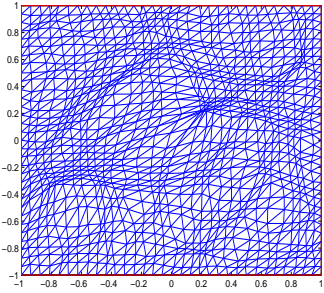}}\\
    \subfigure[$a$ at $t=0.1$.]
    {\includegraphics[width=0.35\textwidth,height= 0.3\textwidth]{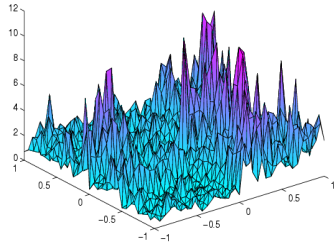}}
    \goodgap
    \subfigure[$(F_1,F_2)$ at $t=0.1$.]
    {\includegraphics[width=0.35\textwidth,height= 0.3\textwidth]{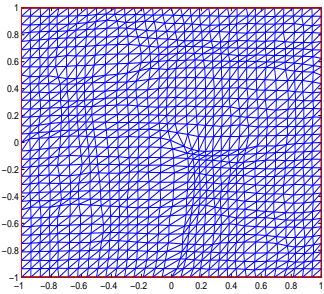}}\\
    \caption{$a$ and $(F_1,F_2)$ at time $t=0$, $t=0.1$ for the time dependent random fractal medium.}
    \label{mediatdp6}
\end{center}
\end{figure}

\begin{figure}[httb]
  \begin{center}
    \subfigure[coarse mesh $L^1$ error.]
    {\includegraphics[width=0.35\textwidth,height= 0.3\textwidth]{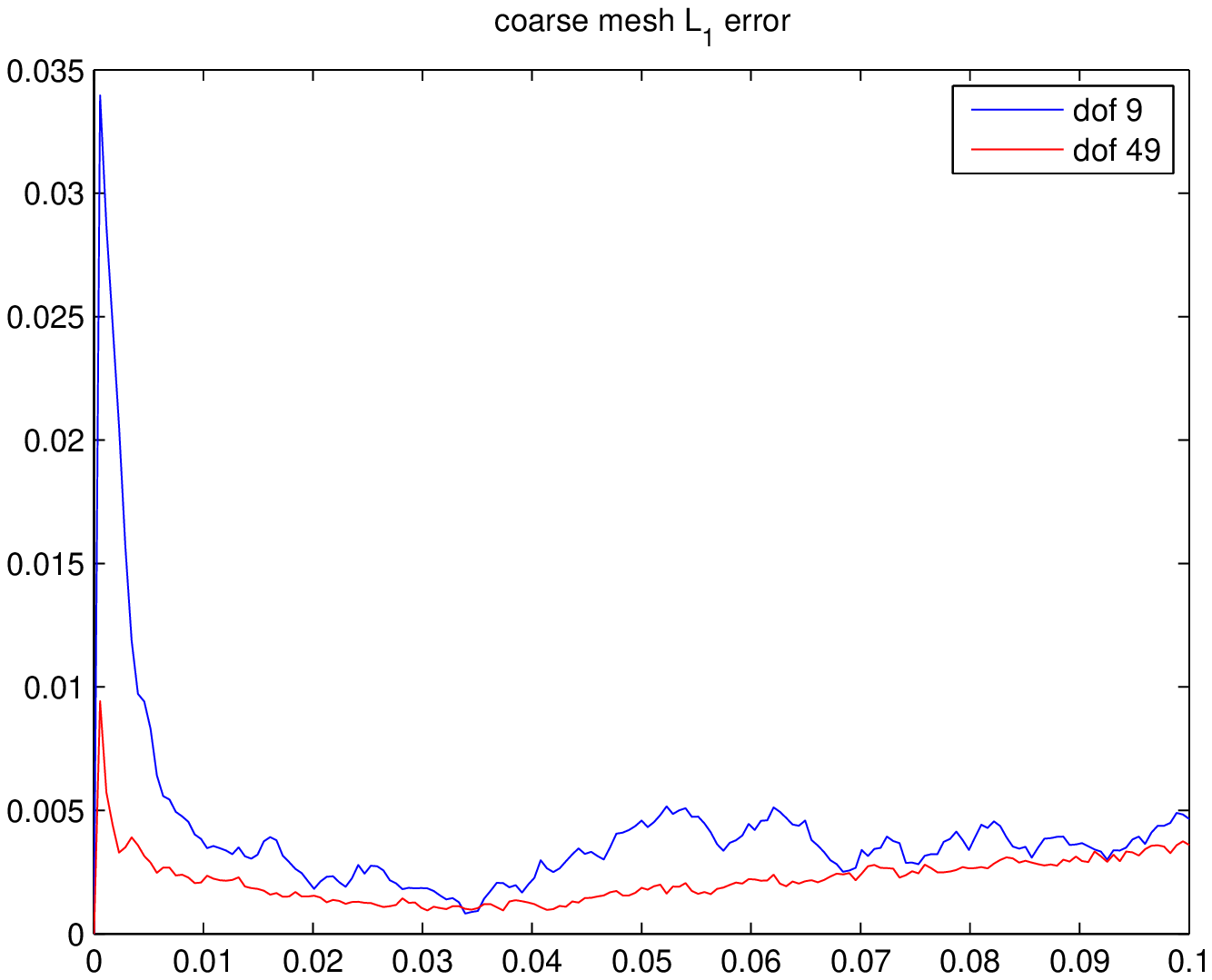}}
    \goodgap
    \subfigure[fine mesh $L^1$ error.]
    {\includegraphics[width=0.35\textwidth,height= 0.3\textwidth]{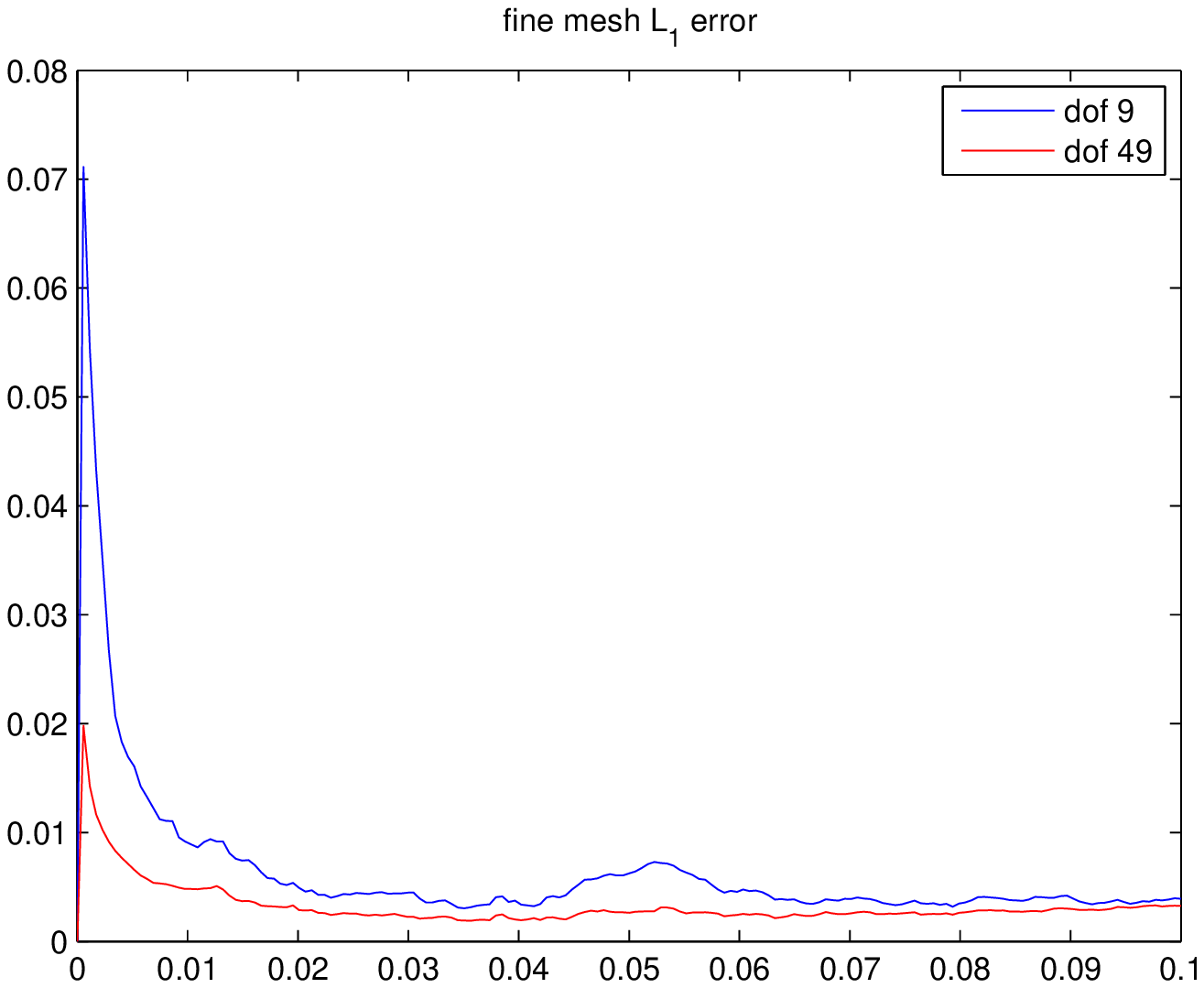}}\\
    \caption{$L^1$ error for the time dependent random fractal medium at $t=0.1$.}
    \label{errL1t10tdp6}
\end{center}
\end{figure}

\begin{figure}[httb]
  \begin{center}
    \subfigure[coarse mesh $L^2$ error.]
    {\includegraphics[width=0.35\textwidth,height= 0.3\textwidth]{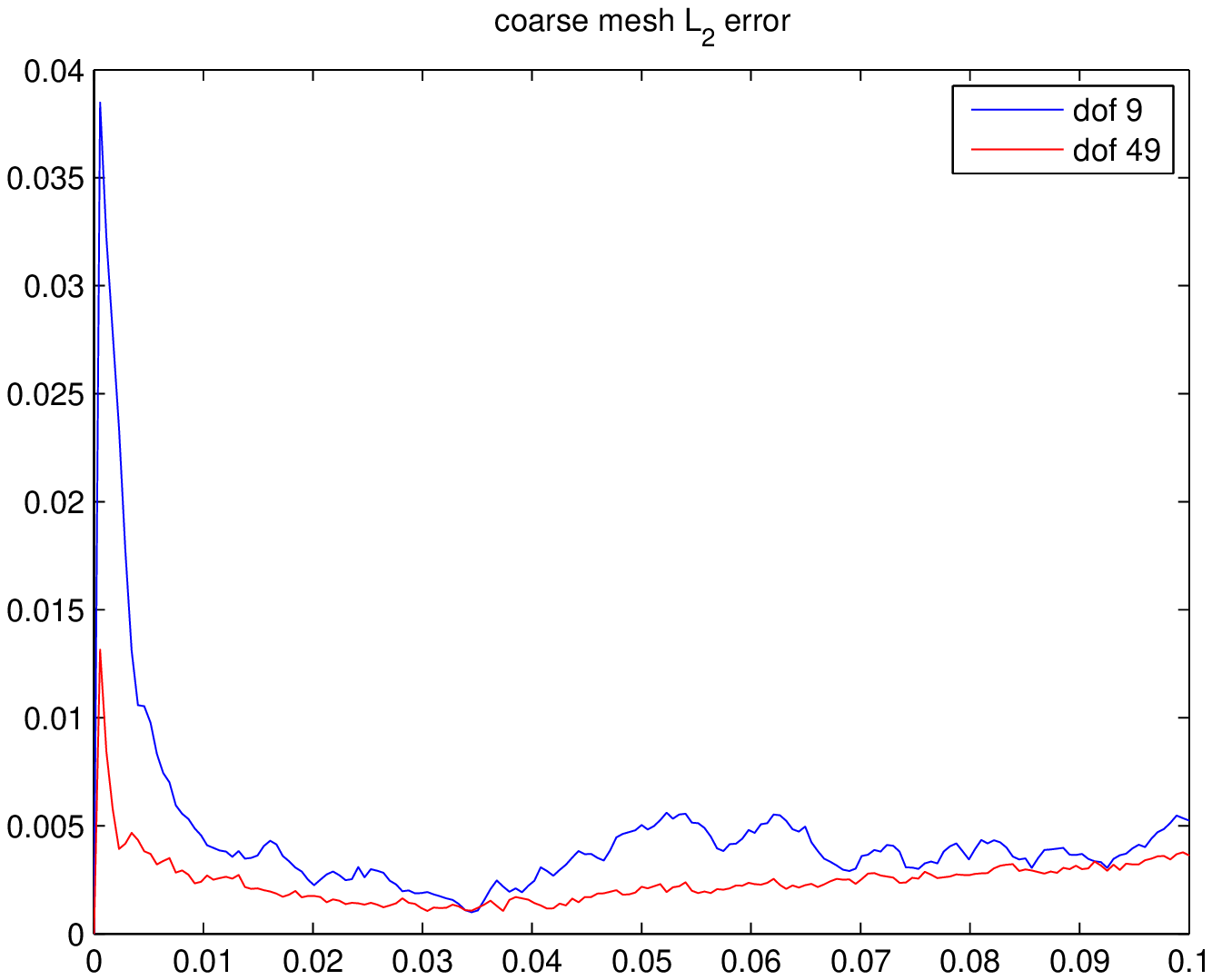}}
    \goodgap
    \subfigure[Fine Mesh $L^2$ error.]
    {\includegraphics[width=0.35\textwidth,height= 0.3\textwidth]{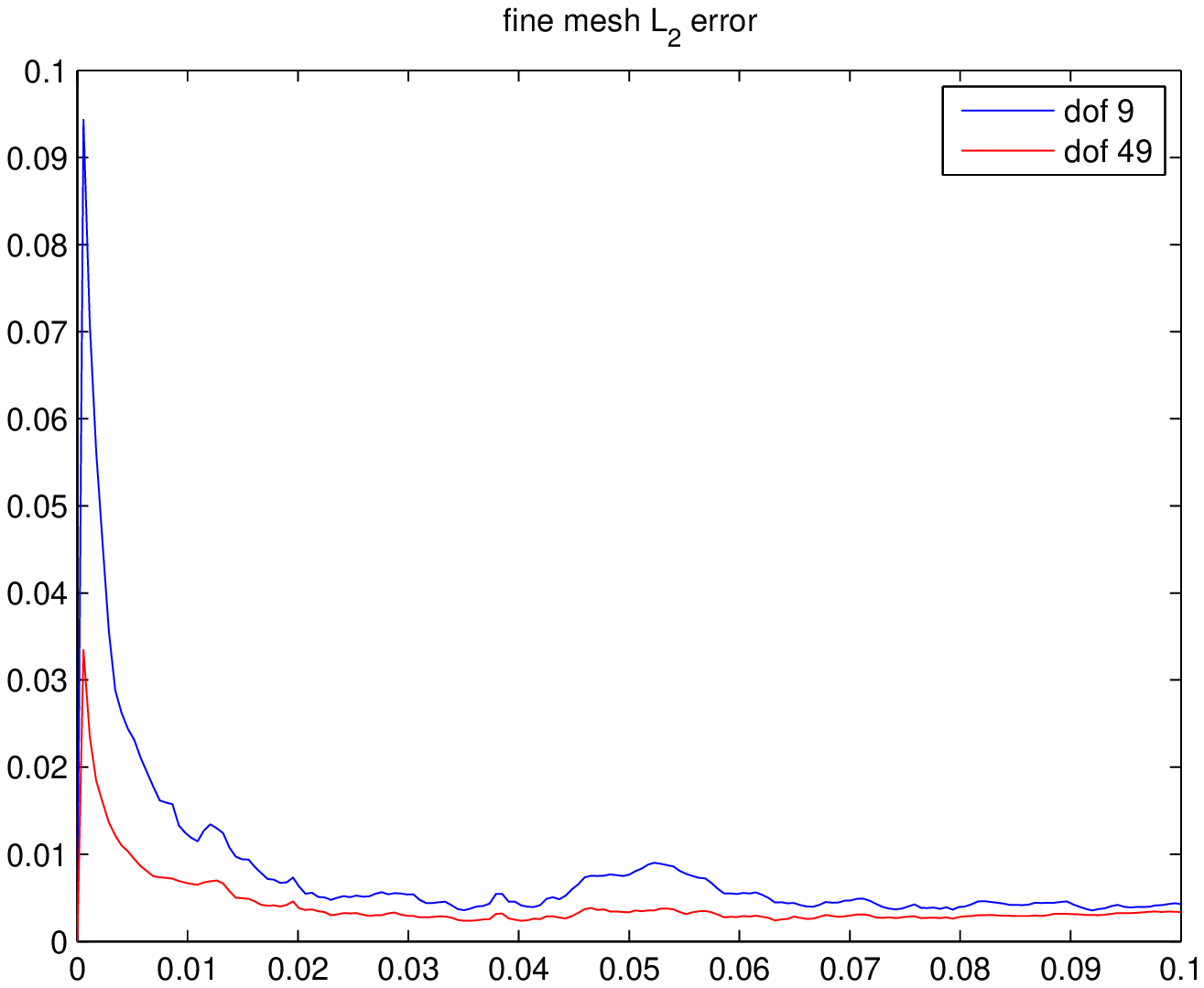}}\\
    \caption{$L^2$ error for the time dependent random fractal medium at $t=0.1$..}
    \label{errL2t10tdp6}
\end{center}
\end{figure}

\begin{figure}[httb]
  \begin{center}
    \subfigure[coarse mesh $L_{\infty}$ error.]
    {\includegraphics[width=0.35\textwidth,height= 0.3\textwidth]{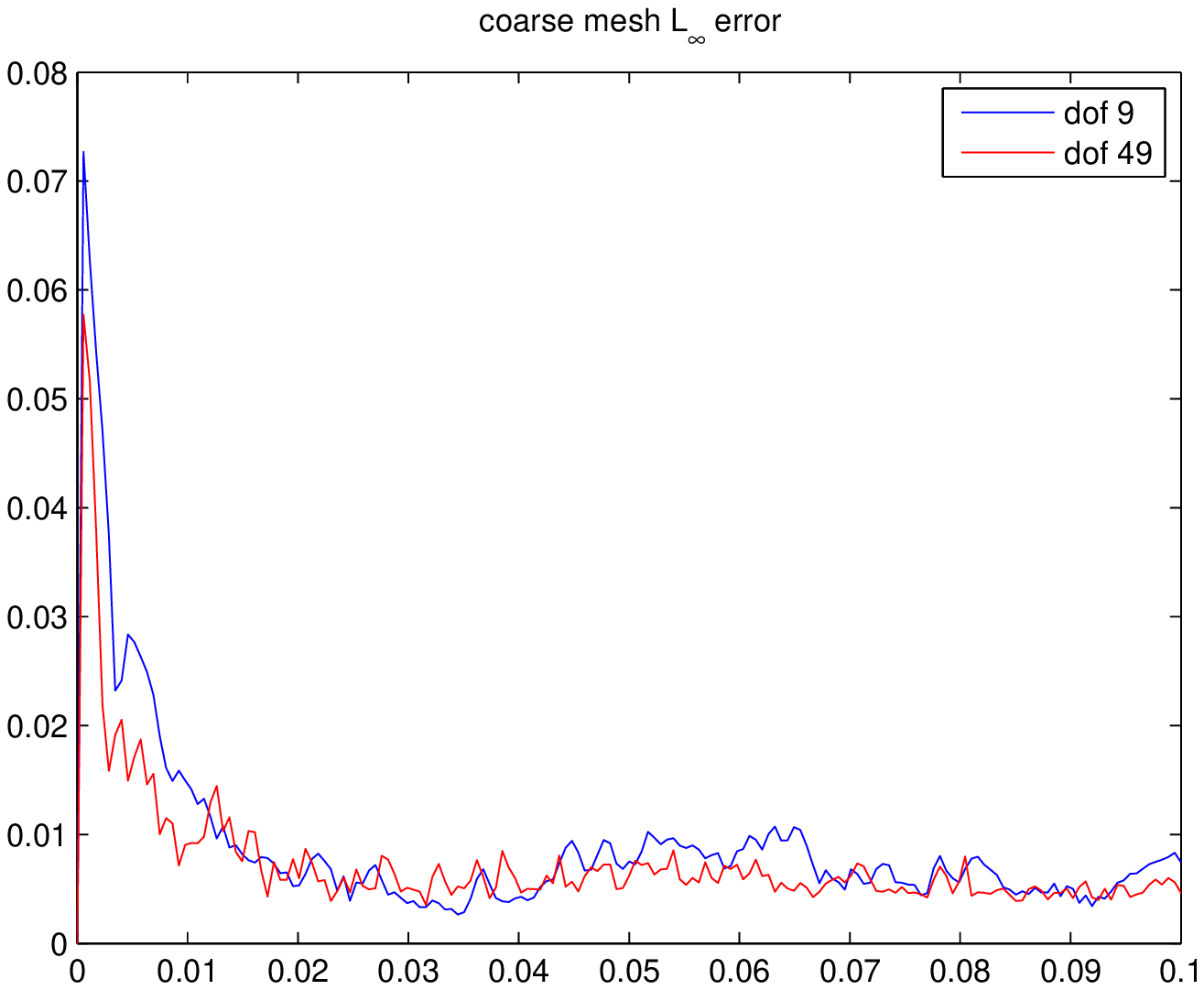}}
    \goodgap
    \subfigure[fine mesh $L_{\infty}$ error.]
    {\includegraphics[width=0.35\textwidth,height= 0.3\textwidth]{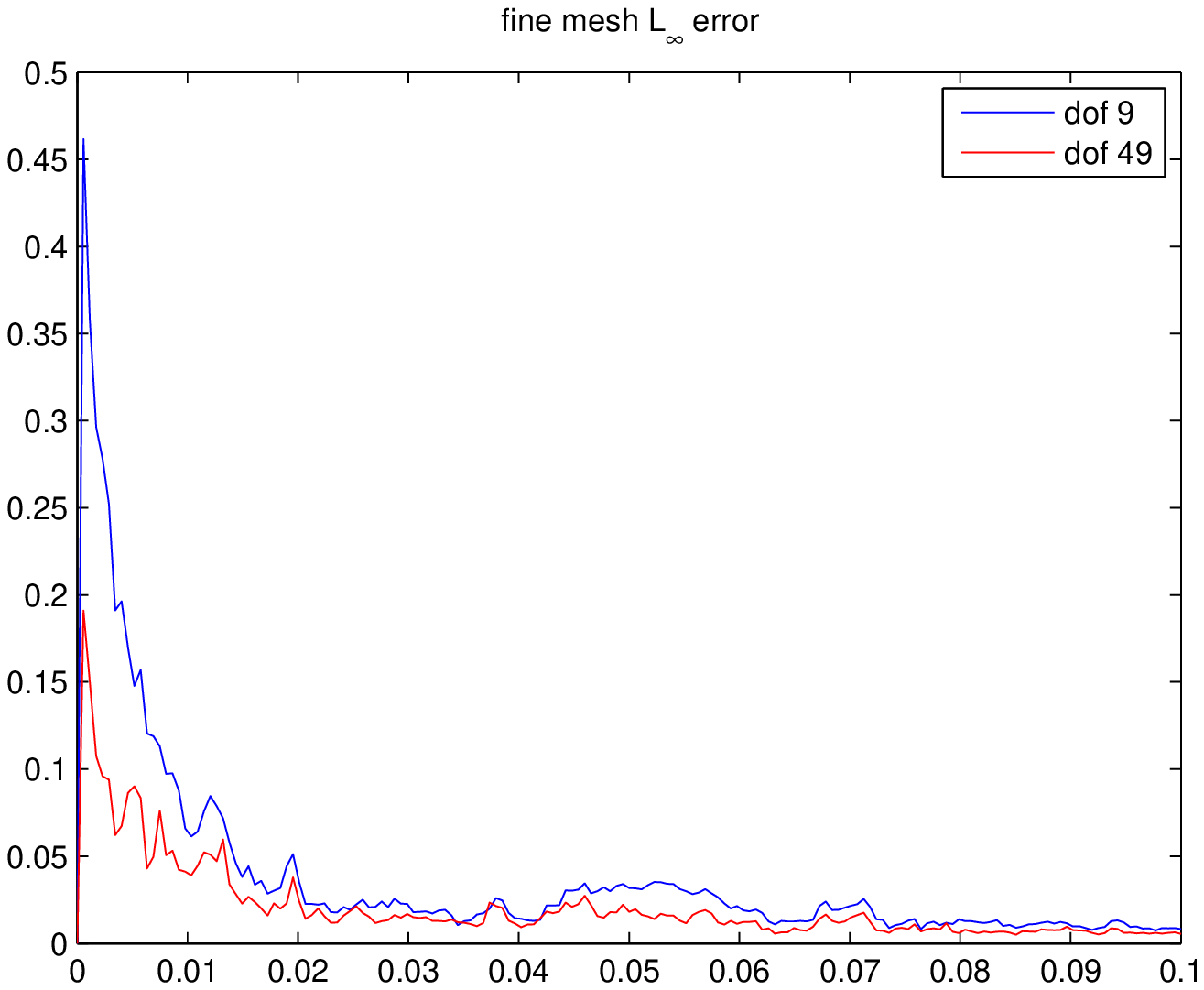}}\\
    \caption{$L^{\infty}$ error for the time dependent random fractal medium at $t=0.1$..}
    \label{errLit10tdp6}
\end{center}
\end{figure}

\begin{figure}[httb]
  \begin{center}
    \subfigure[coarse Mesh $H1$ error.]
    {\includegraphics[width=0.35\textwidth,height= 0.3\textwidth]{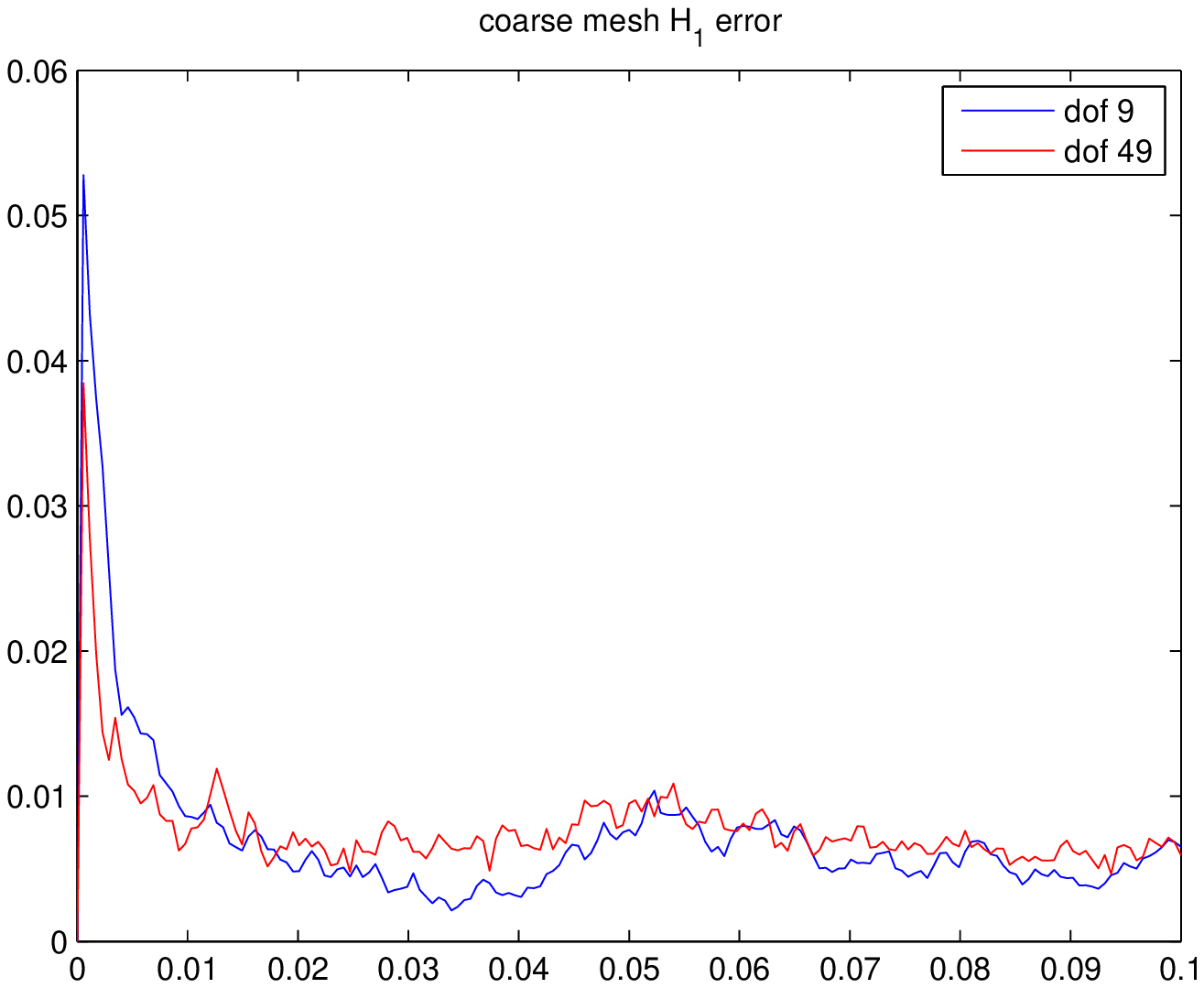}}
    \goodgap
    \subfigure[fine Mesh $H^1$ error.]
    {\includegraphics[width=0.35\textwidth,height= 0.3\textwidth]{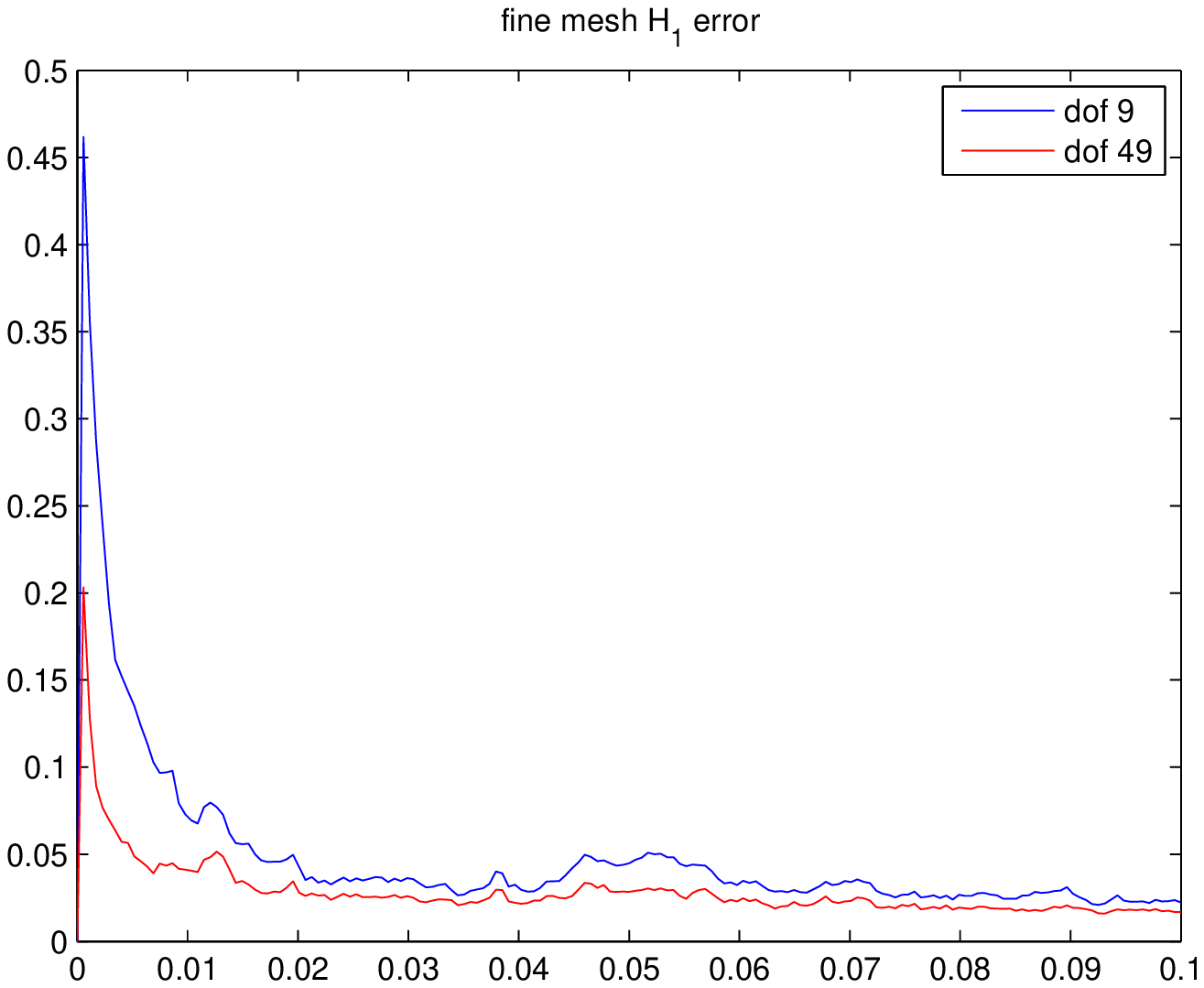}}\\
    \caption{$H^1$ error for the time dependent random fractal medium at $t=0.1$..}
    \label{errH1t10tdp6}
\end{center}
\end{figure}

\begin{example}
\label{exa:tdep_randomfractal}Time Dependent Random Fractal
\end{example}
In this example, $a$ is given by a product of discontinuous
functions oscillating randomly at different space and time scales.
Namely $a(x,t):=a_{1}(x,t)a_{2}(x,t)\cdots a_{n}(x,t)$ with $n=6$
and $a_{i}(x,t)=c_{pq}^i(t)$ for
$x\in[\frac{p}{2^{i}},\frac{p+1}{2^{i}})\times[\frac{q}{2^{i}},\frac{q+1}{2^{i}})$.
The coefficients $c_{pq}^i(t)$ are chosen at random with uniform law
 in $[\frac{1}{\gamma},\gamma]$ with $\gamma=0.7$ and independently
 in subdivision in space and in time, thus they are assumed to be
 constant in each time interval
 $0.1\times[\frac{k}{4^i},\frac{k+1}{4^i})$.  In this example we have $\frac{\lambda_{\max}(a)}{\lambda_{\min}(a)}=160.3295$. Although the number of
 fine
time steps to solve \eref{ghjh52} is $3482$, only $175$ coarse time
steps have been used to solve the homogenized equation which
corresponds to a reduction of the complexity of the scheme by a
factor of $35000$ in the case of the coarse tessellation with $9$
interior nodes. $a$ and the map $(F_1,F_2)$ are drawn in figure
\ref{mediatdp6}.  $L^1$, $L^2$, $L_{\infty}$ and $H_1$ errors are
given in figure \ref{errL1t10tdp6} to \ref{errH1t10tdp6}. Coarse and
fine mesh errors are given in table \ref{cerrtdp6sp} and
\ref{ferrtdp6sp} at time $t=0.1$. We have chosen $g=1$ in this
numerical experiment, one obtains similar results by choosing
$g=\sin(2.4 x-1.8y+2 \pi t)$.

\begin{table}
\begin{center}
\caption{Coarse Mesh Error for the time dependent random fractal
medium with spline element} \label{cerrtdp6sp}
\begin{tabular}{|c|c|c|c|c|}
\hline dof& $L^{1}$& $L^{\infty}$& $L^{2}$& $H^{1}$\tabularnewline
\hline 9& 0.0046& 0.0074& 0.0052& 0.0065\tabularnewline \hline 49&
0.0036& 0.0046& 0.0036& 0.0059\tabularnewline \hline
\end{tabular}
\end{center}
\end{table}

\begin{table}
\begin{center}
\caption{Fine Mesh Error for the time dependent random fractal
medium with spline element} \label{ferrtdp6sp}
\begin{tabular}{|c|c|c|c|c|}
\hline dof& $L^{1}$& $L^{\infty}$& $L^{2}$& $H^{1}$\tabularnewline
\hline 9& 0.0039& 0.0082& 0.0043& 0.0222\tabularnewline \hline 49&
0.0033& 0.0054& 0.0034& 0.0168\tabularnewline \hline
\end{tabular}
\end{center}
\end{table}

\clearpage

\def\polhk#1{\setbox0=\hbox{#1}{\ooalign{\hidewidth
  \lower1.5ex\hbox{`}\hidewidth\crcr\unhbox0}}} \def\cprime{$'$}
  \def\cprime{$'$}

\end{document}